\documentclass[a4paper,10pt]{article}
\usepackage{psfrag,
graphicx, harvard,
epsf,mathrsfs,
euscript,amsfonts,amsmath,latexsym,amssymb,mathrsfs
}
\usepackage[active]{srcltx}
\usepackage{fullpage}
\usepackage{color}

\newtheorem{lemma}{Lemma}[section]
\newtheorem{theorem}{Theorem}[section]
\newtheorem{proposition}{Proposition}[section]

\newtheorem{remark}{Remark}[section]
\newtheorem{definition}{Definition}[section]
\newtheorem{corollary}{Corollary}[section]

\newcommand{\cB}{{\cal B}}

\newcommand{\cL}{{\cal L}}
\newcommand{\cM}{{\cal M}}
\newcommand{\cN}{{\cal N}}

\newcommand{\cP}{{\cal P}}

\newcommand{\cR}{{\cal R}}
\newcommand{\cS}{{\cal S}}

\newcommand{\bC}{\mathbb C}

\newcommand{\bE}{\mathbb E}

\newcommand{\bL}{{\mathbb L}}

\newcommand{\bN}{{\mathbb N}}

\newcommand{\bP}{{\mathbb P}}

\newcommand{\bR}{{\mathbb R}}

\newcommand{\bZ}{{\mathbb Z}}
\newcommand{\sF}{{\mathscr F}}

\newcommand{\sH}{{\mathscr H}}

\newcommand{\sN}{{\mathscr N}}

\newcommand{\sS}{{\mathscr S}}

\newcommand{\rd}{\mathrm{d}}

\newcommand{\pr} {\par \noindent{\bf Proof\,:~}}
\newcommand{\epr}{\hfill\hbox{\hskip 4pt
                \vrule width 5pt height 6pt depth 1.5pt}\vspace{0.5cm}\par}
\begin{document}
\numberwithin{equation}{section}
\title{Smoluchowski processes and   
nonparametric estimation of functionals of
particle displacement distributions
\\
from count data 
\thanks{The work
is supported by the Israel Science Foundation (ISF) research grant.}}
\author{
A. Goldenshluger\thanks{Department of Statistics, University of Haifa, Haifa 31905, Israel.  e-mail:
{\em goldensh@stat.haifa.ac.il}. }
\and
\and R. Jacobovic\thanks{Department of Statistics, University of Haifa, Haifa 31905, Israel.  e-mail:
{\em royi.jacobovic@mail.huji.ac.il}.}
}
\maketitle
\begin{abstract}
Suppose that particles are randomly distributed in $\bR^d$, and they 
are subject to identical 
stochastic  motion independently of each other. 
The Smoluchowski process describes  fluctuations of the number of particles 
in an observation region over time.  
This paper studies properties of the Smoluchowski processes and considers related statistical problems. 
In the first part of the paper we revisit  probabilistic properties of the Smoluchowski process in a unified and principled way: explicit formulas for generating functionals
and moments are derived, conditions for stationarity and Gaussian approximation are discussed, and 
relations to other stochastic models are highlighted. 
The second part deals with statistics of the Smoluchowki processes. We consider two
different 
models of the particle displacement process: the undeviated uniform motion (when 
a particle moves with 
random constant velocity along  a straight line) and the Brownian motion displacement.
In the setting of the undeviated uniform motion we study the problems of estimating the mean speed
and the speed distribution, while for the Brownian displacement model the problem of estimating the diffusion 
coefficient is considered. In all these settings we develop estimators with provable accuracy guarantees.   
\end{abstract}
\vspace*{1em}
\par 
\vspace*{1em}

\section{Introduction}
Suppose that we have an infinite number of particles that are randomly distributed in $\bR^d$. 
The particles are subject to  the same stochastic movement independently of each other.  
We are interested in  characteristics of the stochastic process that
governs the particle movement, e.g., in  diffusion coefficients.
The natural assumption is that the particles are indistinguishable, and the  trajectory of a
single particle cannot be tracked over time.    
In this setting  \citeasnoun{Smoluchowski06} assumed that  particles perform 
Brownian motion, and suggested to measure the number of particles 
 (concentration)  in a fixed region 
over time in order to determine  unknown movement characteristics.
The classical works of Smoluchowski on  {\em concentration  fluctuations}
are of fundamental importance in statistical physics:
they are  
the part 
of the celebrated Einstein--Smoluchowski theory  that 
provided a molecular--kinetic  explanation of the Brownian movement.   
The exposition of probabilistic 
aspects of Smoluchowski's theory is given in classical 
surveys by \citeasnoun{Chandrasekhar} and \citeasnoun[Chapter~III, Sections~21-28]{Kac}, and, 
more recently, in \citeasnoun{Bingham};  we refer also to \citeasnoun{Mazo} for 
historical background and additional information.
\paragraph{The model.}
Consider the following model of moving particles in $\bR^d$.
Let   $\Xi:= \sum_{j\in \bZ} \varepsilon_{\xi_j}$ be a homogeneous Poisson process on 
$(\mathbb{R}^d,\mathcal{B}(\mathbb{R}^d))$ with rate $\lambda\in(0,\infty)$; here $\varepsilon_x$ is the 
Dirac measure at $x\in \bR^d$.
The location of $j$-th particle at time $t$ is determined by the equation
\begin{equation}\label{eq:X-t}
 X_t^{(j)}=\xi_j + Y_t^{(j)},\;\;\;j\in \bZ,\;\;t\geq 0,
\end{equation}
where 
$\{Y_t^{(j)}, t\geq 0\}$, $j\in \bZ$ are independent 
copies of  random process $\{Y_t, t\geq 0\}$ in $\bR^d$ with $Y_0=0$. 
In what follows process $\{Y_t, t\geq 0\}$
is referred to as 
{\em the  displacement process}. 
We assume 
that Poisson process~$\Xi$ describing the initial positions of particles and displacement processes  
$\{Y_t^{(j)}, t\geq 0\}$, $j\in \bZ$   are 
independent.
We also denote
$\{X_t, t\geq 0\}$ to be a generic location process with
$\{X_t^{(j)}, t\geq 0\}$, $j\in \bZ$  being  
the independent copies of $\{X_t, t\geq 0\}$. 
\par 
Let $B$ be a compact set with non--empty interior  
in $\bR^d$ representing {\em the observation region},  and define  
\begin{equation}\label{eq:N-t}
 N(t):= \sum_{j\in \bZ} {\bf 1}_{B}\big(X^{(j)}_t\big)= \sum_{j\in \bZ} {\bf 1}\big\{X^{(j)}_t \in B\big\},\;\;\;t\geq 0.
\end{equation}
By definition, process $\{N(t), t\geq 0\}$ counts   
the number of particles in the observation region $B$  over time. 
In all what follows we refer to random process $\{N(t), t\geq 0\}$  as 
{\em the Smoluchowski process}. 
\par\medskip  
Different particle displacement processes $\{Y_t, t\geq 0\}$ give rise to different versions
of the Smoluchowski process.
The following   two models of displacement  are of  particular interest. 
\begin{itemize}
 \item[(a)] {\em Undeviated uniform motion.} 
 Let 
$(v_j)_{j\in\mathbb{Z}}$ be  a sequence of iid random vectors in $\mathbb{R}^d$ with 
common distribution function $G$, independent of 
$\Xi$. Assume that 
\begin{equation}\label{eq:undeviated}
 Y_t^{(j)}= v_j t,\;\;\;j\in\bZ,\;\;t\geq 0.
\end{equation}
In this model particles move along  straight lines  with constant velocity
which varies from particle to particle according to probability distribution $G$.
We refer to $v_j\in \bR^d$ as the $j$-th {\em particle velocity}, while the Euclidean norm $\|v_j\|$ of $v_j$
is the corresponding  {\em particle speed}.
\item[(b)] {\em Brownian displacement model.} Assume that 
\begin{equation}\label{eq:brownian}
Y_t^{(j)}= \sigma W_t^{(j)},\;\;\;j\in \bZ,\;\;t\geq 0, 
\end{equation}
 where $\{W_t^{(j)},\, t\geq 0\}$, $j\in \bZ$ are independent standard 
$d$--dimensional 
Brownian motions in $\bR^d$, and $\sigma>0$ is the diffusion coefficient.  Under this setting
 particles  perform Brownian motion. 
\end{itemize}
\par 
The Smoluchowski process $\{N(t), t\geq 0\}$ associated 
with Brownian displacement~(\ref{eq:brownian}) was originally considered 
by \citeasnoun{Smoluchowski06} in the context of  his studies on colloidal suspensions.
The model of undeviated uniform motion~(\ref{eq:undeviated}) goes back to the work 
of F\"urth who applied Smoluchowski's methods 
to estimate the average speed of pedestrians from counts 
of the number of pedestrians in a fixed section of a road [see, \citeasnoun[Chapter~III, Section~26]{Kac}]. 
The term {\em undeviated uniform motion}
was coined  by  \citeasnoun{Lindley}.
\paragraph{Related literature.}
\citeasnoun[Chapter~III]{Chandrasekhar} discusses probabilistic properties of the Smoluchowski process. 
The presentation there does not explicitly state the probabilistic model, and 
it is   tacitly assumed that $\{N(t), t\geq0\}$ is Markovian even though this is not in general
true.
The review of  probabilistic properties of $\{N(t), t\geq 0\}$ that is most relevant to our work is 
\citeasnoun{Kac}.
In this work   expressions 
for the moment generating functional and finite dimensional distributions of the Smoluchowski process are 
derived. \citeasnoun{Kac} also discusses  the issue of Markovianity of the Smoluchowski process, 
considers the undeviated uniform motion model, and explains significance  
of the probabilistic results for problems of estimating physical entities from count data. 
Motivated by the F\"urth traffic problem and by the work of \citeasnoun{Rothschild}
on measuring motility of spermatozoa,
\citeasnoun{Lindley} considers  the problem of estimating expected  
velocity of 
the undeviated uniform motion on the real  line. This paper relates 
the expected velocity to the 
one--sided derivative at zero  
of the covariance function of the Smoluchowski process, and uses
this relation for constructing an estimator.  
\citeasnoun{Ruben64} suggests a generalization of the Smolluchowski process assuming that
the number of 
particles can be counted in several disjoint observation regions. Properties of  
the multivariate Smoluchowski process constructed in this way are studied in the aforementioned work, and 
it is argued that the use of such multivariate data improves accuracy of 
estimation procedures.  
Parameter estimation problems for 
the multivariate Smoluchowski processes are considered in 
\citeasnoun{McDunnough79b}.   
These papers promote the idea  of using  standard 
asymptotic results for estimating covariance functions of discrete time series
in conjunction with the delta method in order to derive
asymptotic distributions of parameter estimators. 
\par 
Another strand of research deals with 
stochastic models that are 
closely related to the Smoluchowski processes:
the $M/G/\infty$ queueing model and the branching processes with immigration. We postpone the detailed 
discussion of 
the connection between these models and the Smoluchowski processes to Section~\ref{sec:related};
here we restrict ourselves with brief description and related references.
\par 
\citeasnoun{Bingham} discuss the use of the $M/G/\infty$ queueing model in some problems of    
statistical  estimation for  Smoluchowski processes. 
Because the Smoluchowski process is not in general Markovian,   \citeasnoun{Bingham} suggest that 
there are two ways to handle the mathematical difficulties: (a)~to assume that
the service time distribution~$G$  is exponential which  
leads to the Markovian $M/M/\infty$ model; (b)~to work with reduced data when 
the observations of the busy and idle periods are available, i.e., 
the process $\{{\bf 1}(N(t)>0), t\geq 0\}$ is observed.    In this paper 
we demonstrate that problems
of statistical estimation for Smoluchowski processes can be solved in the original setting
without resorting to restrictive assumptions like~(a) or~(b).
\par 
A Markovian model that can be used as an approximation for the Smoluchowski process is
the branching process with immigration. In this context 
the model has been studied 
in \citeasnoun{Ruben62} and  \citeasnoun{Ruben64}.  
The relationship between branching processes 
with immigration and Smoluchowski processes  is discussed in detail 
in \citeasnoun{McDunnough78} where necessary and sufficient conditions for  Markovianity of 
the Smoluchowski 
processes are derived. 
Some statistical estimation problems for branching processes with immigration
are considered in \citeasnoun{Heyde-Seneta} and \citeasnoun{Heyde-Seneta-1};
we also refer to \citeasnoun{Wei-Winnicki} and \citeasnoun{Winnicki} for related results. 
\par 
We  note that there exists a considerable body of work dealing with 
statistical inference for diffusion processes and infinite particle
systems under assumption that trajectories of  particles are directly observable.
The problems of  parametric and nonparametric estimation  of drift and diffusion 
coefficients in these models
were intensively studied; 
we refer, 
e.g., to \citeasnoun{Jacod1}, \citeasnoun{Jacod2},  \citeasnoun{Hoffmann}, \citeasnoun{Gobet} and 
\citeasnoun{Kutoyants}
where further references can be found. However, as it is shown in the present paper, 
estimation settings based on  the count data lead to statistical inverse problems that require
completely different techniques and tools. One of the goals of this paper is to develop such techniques. 
\par 
Finally, it is also worth mentioning that although  Smoluchowski's theory originated in 
 statistical physics, 
it found numerous applications in  
diverse areas, e.g.,  in biology \cite{Rothschild}, spectroscopy~\cite{Brenner},
medicine~\cite{Aebersold} and geology~\cite{Culling}.
\paragraph{The paper contribution.}
Our goal in  this paper is two--fold. First, we introduce a general model of moving particles in
$\bR^d$ and revisit probability properties of the Smoluchowski processes  using original proof techniques. 
In the existing literature these properties are discussed or mentioned in passing 
in different sources under  
disparate and sometimes not fully specified assumptions on the probabilistic model. 
In contrast, the framework taken in this paper  enables us to investigate  
probability properties of Smoluchowski processes associated with 
arbitrary displacement models
in a unified and principled way.
\par 
Second, we are interested in estimation of some  
functionals of the particle displacement
distributions  from  continuous time observation  
$\sN_T:=\{N(t), 0\leq t \leq T\}$ of the Smoluchowski process. 
These estimation problems are reduced  to 
estimating functionals 
of the correlation function of the Smoluchowski process from indirect observations.
The resulting statistical ill--posed inverse problems depend on the  geometry of the 
observation region and require the use of special transform methods for 
constructing the estimators.  
Although the developed methods are applicable for observation regions  of arbitrary shape, 
in all what follows we focus
on the case when $B$ is a Euclidean ball in $\bR^d$.
In the setting of the undeviated uniform motion 
we consider problems of estimating the mean speed and the speed distribution, 
while in the setting of the Brownian displacement our focus is 
on estimating the diffusion coefficient.  
In all aforementioned settings  we develop estimators with provable accuracy guarantees. 
\paragraph{Organization of the paper.} 
The rest of the paper is structured as follows. Section~\ref{sec:properties} discusses probabilistic properties 
of the Smoluchowski process. In particular, we derive explicit expressions for the 
moment generating functional and mixed moments, present conditions for  
stationarity of the Smoluchowski processes, establish a Gaussian approximation 
and discuss connections to the $M/G/\infty$ queue and branching processes
with immigration. In Section~\ref{sec:covariance} we present results on estimation of the covariance function 
of the Smoluchowski process; these results play an important role   in 
all subsequent developments. 
Section~\ref{sec:undeviated} deals with estimation of the expected speed and speed distribution 
for the Smoluchowski process driven by the  undeviated uniform motion.
Section~\ref{sec:Brownian}
considers the Brownian displacement model, and studies estimation of the diffusion coefficient. 
Proofs of results of Sections~\ref{sec:covariance}, \ref{sec:undeviated} and \ref{sec:Brownian}
are given in Sections~\ref{sec:proofs-cov}, \ref{sec:proofs-undeviated} and \ref{sec:proofs-Brownian} respectively.
\paragraph{Notation.} The following notation  is used throughout the paper.
The $d$--dimensional volume (the Lebesgue content) 
of a set in $\bR^d$ is denoted 
${\rm vol}\{\cdot\}$.  For a set $C\subset \bR^d$ and point $x\in \bR^d$ we denote  
\[ 
C(x):=\{y\in \bR^d: y = z-x,\,z\in C\}. 
\]
{\em The covariogram}
of a compact set $C\subset \bR^d$ 
 is the function $g_C:\bR^d \to \bR$ defined by 
\begin{equation}\label{eq:covariogram}
 g_C(x):= {\rm vol}\{ C \cap C(x)\},\;\;x\in \bR^d;
\end{equation}
see, e.g., \citeasnoun[Section~4.3]{Matheron}. 
The Euclidean norm on $\bR^d$ is denoted 
$\|\cdot\|$. 
%
\section{Properties of the Smoluchowski process}\label{sec:properties}
In this section we discuss probabilistic properties of Smoluchowski processes.
Some of the presented results  appeared 
in different forms and for different settings, e.g., in
\citeasnoun{Lindley}, \citeasnoun{Kac}, 
and \citeasnoun{Ruben64}.
Our approach to  derivation of these results is purely analytic, and our proofs 
differ from those in the existing literature.
It is also worth to emphasize 
that the results of this section hold for any Smoluchowski's process, 
independently  of the displacement model.
\subsection{Moment generating functional}
Let $\Pi_n$ denote the set of all {\em non--empty} subsets of 
$\{1,2,\ldots, n\}$, and for $\pi \in \Pi_n$  write 
$\pi^c:= \{1, 2, \ldots, n\} \setminus \pi$. The following theorem establishes the moment generating functional
of process $\{N(t), t\geq 0\}$ defined in~(\ref{eq:N-t}).
\begin{theorem}\label{th:Eexp}
 For any $\theta=(\theta_1, \ldots, \theta_n)\in \bR^n$ and $t_1\leq t_2\leq \cdots\leq t_n$ one has 
 \begin{equation}\label{eq:formula}
  \ln \bE \exp\Big\{\sum_{k=1}^n \theta_k N(t_k)\Big\}=\lambda \sum_{\pi\in \Pi_n}
 \big[e^{\sum_{k\in \pi}\theta_k} -1\big] \bar{Q}_\pi (t_1,\ldots, t_n), 
 \end{equation}
 where 
 \begin{eqnarray}
 &&\bar{Q}_\pi (t_1, \ldots, t_n) := \int_{\bR^d} Q_\pi (x; t_1, \ldots, t_n) \rd x,
 \label{eq:Q-pi0}
 \\*[2mm]
 &&Q_{\pi}(x; t_1, \ldots, t_n) := \bP\big\{ \big (x+Y_{t_k} \in B, k\in \pi\big)\,\cap\, 
 \big (x+Y_{t_k}\notin B,\;
 k\in \pi^c \big) \big\}.
 \label{eq:Q-pi}
 \end{eqnarray}
\end{theorem}
\begin{remark}
 Note that function $\bar{Q}_\pi(t_1, \ldots, t_n)$ is well defined because the 
 integral in (\ref{eq:Q-pi0}) is bounded from above by $\bE\, {\rm vol}\big\{\cap_{k\in \pi} B(Y_{t_k})\big\} \leq  {\rm vol}(B)$
 for any $\pi \in \Pi_n$.
\end{remark}
\pr
For  $j\in \bZ$ and  $\pi\in \Pi_n$  define events
\[
 A_\pi^{(j)}:=\big\{X_{t_k}^{(j)}\in B, k\in \pi\big\} \cap \big\{ X^{(j)}_{t_k}\notin B, k\in\pi^c\big\}.
\]
Event $A_{\pi}^{(j)}$ states that $j$-th particle is inside  $B$ at time instances 
$\{t_k, k\in \pi\}$ and outside $B$ otherwise. 
Obviously, for fixed $j$ events $A_{\pi}^{(j)}$, $\pi\in \Pi_n$ are disjoint. We let 
$A^{(j)}:= \cup_{\pi\in \Pi_n} A_\pi^{(j)}$, and  $\bar{A}^{(j)}$  the event complimentary to 
$A^{(j)}$, i.e., $\bar{A}^{(j)}=\big\{X_{t_k}^{(j)}\notin B, \,\forall k\in \{1, \ldots, n\}\big\}$.
Since  
\[
 \sum_{k=1}^n \theta_k N(t_k)= \sum_{j\in \bZ} \sum_{k=1}^n \theta_k {\bf 1}\{X_{t_k}^{(j)}\in B\},
\]
with the introduced notation we have 
for any $j\in \bZ$ 
\[
 \sum_{k=1}^n \theta_k {\bf 1}\big\{X_{t_k}^{(j)} \in B\big\} = \sum_{\pi\in \Pi_n} 
 {\bf 1}\big\{A_\pi^{(j)}\big\} 
 \sum_{k\in \pi} \theta_k,
\]
and therefore 
\begin{align}
 \exp\Big\{\sum_{k=1}^n \theta_k {\bf 1}\big(X_{t_k}^{(j)}\in B\big)\Big\}=
 \sum_{\pi\in \Pi_n} {\bf 1}\big\{A_{\pi}^{(j)}\big\} e^{\sum_{k\in \pi}\theta_k} + 
 {\bf 1}\big\{\bar{A}^{(j)}\big\}
 \nonumber
 \\
 = \sum_{\pi\in \Pi_n} {\bf 1}\big\{A_{\pi}^{(j)}\big\} \big[e^{\sum_{k\in \pi}\theta_k}-1\big] + 1. 
\label{eq:exp}
 \end{align}
Then  in view of (\ref{eq:exp})
\begin{align*}
 \bE \bigg[\exp \Big\{\sum_{k=1}^n \theta_k N(t_k)\Big\}\Big| \,\Xi\bigg] &=
 \bE \bigg[ \exp \Big\{ \sum_{j\in \bZ} \sum_{k=1}^n \theta_k {\bf 1}\big(X_{t_k}^{(j)}\in B\big)
 \Big\} \Big| \,\Xi\bigg]
 \\
 &= \prod_{j\in \bZ} \bE \bigg[ \exp \Big\{\sum_{k=1}^n \theta_k {\bf 1}\big(X_{t_k}^{(j)}\in B\big)
 \Big\} \Big| \,\Xi \bigg]
 \\
 &= 
 \prod_{j\in \bZ} \bigg(1+
 \sum_{\pi\in \Pi_n} \bP\big\{A_\pi^{(j)}\, | \,\Xi\big\}
  \big[e^{\sum_{k\in \pi}\theta_k}-1\big]\bigg) 
  = 
  \exp\Big\{\sum_{j\in \bZ} f(\xi_j)\Big\},
\end{align*}
where 
\[
 f(\xi_j)=\ln \Big(1+
 \sum_{\pi\in \Pi_n} \bP\big\{A_\pi^{(j)}\, | \,\Xi\big\}
  \big[e^{\sum_{k\in \pi}\theta_k}-1\big]\Big).
\]
Using Campbell's formula we obtain 
\begin{align*}
 \ln \bE  \Big[\exp \Big\{\sum_{k=1}^n \theta_k N(t_k)\Big\}\Big]= \lambda \sum_{\pi\in \Pi_n}
 \big[e^{\sum_{k\in \pi}\theta_k} -1\big] \int_{\bR^d}  Q_\pi(x) \rd x,
\end{align*}
where function $Q_\pi(x)=Q_\pi(x;t_1, \ldots, t_n)$ is defined in (\ref{eq:Q-pi}).
\epr

\begin{remark}
 Theorem~\ref{th:Eexp} holds under very general assumptions on the displacement
 process $\{Y_t, t\geq 0\}$: only independence of $\{Y_t, t\geq 0\}$ and the initial
 position Poisson point process $\Xi$ is required. In particular, under these conditions 
 the Smoluchowski process $\{N(t), t\geq 0\}$ should not be stationary.
\end{remark}
\par 
Formula (\ref{eq:formula}) implies that for every $t\geq 0$ random 
variable $N(t)$ has Possion distribution
 with parameter 
 \[
 \rho:=\bE N(t) = \lambda\int_{\bR^d} \bP\{x+Y_t\in B\}\rd x = \lambda \bE\, {\rm vol}\{B(Y_t)\}=\lambda {\rm vol}\{B\}.
 \]
Thus the  
expectation of 
$N(t)$ does not depend on properties of displacement process  $\{Y_t, t\geq 0\}$.
This fact is consistent with the well known result on random displacement of Poisson point process
which implies that 
$\sum_{j\in \bZ} \varepsilon_{X_t^{(j)}}$ is   homogeneous Poisson process  of intensity $\lambda$ 
for any $t\geq 0$
[see, e.g., \citeasnoun[Section~5, Chapter~VIII]{Doob}].
Since the expectation of $N(t)$ does not bring any information
on   the particle displacement process, statistical inference for the 
displacement process  should be based on 
moments of higher order.
\par 
It is  instructive to specialize formula (\ref{eq:formula}) for the case  $n=2$:
\begin{align*}
&\frac{1}{\lambda} 
\ln \bE \exp\{\theta_1 N(t_1)+\theta_2 N(t_2)\} 
\\
&\;\;\;= 
(e^{\theta_1}-1) \bar{Q}_{\{1\}}(t_1, t_2)+ (e^{\theta_2}-1) \bar{Q}_{\{2\}}(t_1,t_2)
+ (e^{\theta_1+\theta_2}-1) \bar{Q}_{\{1,2\}}(t_1, t_2),
\end{align*}
where
\begin{align*}
& \bar{Q}_{\{1\}}(t_1, t_2) = \bE \int_{\bR^d} {\bf 1}\Big\{ x\in B(Y_{t_1}) \cap B^c(Y_{t_2})\Big\} \rd x
 =  {\rm vol}\{B\} - \bE {\rm vol} \big\{  B(Y_{t_1}) \cap B(Y_{t_2})\big\},
 \\
& \bar{Q}_{\{2\}}(t_1, t_2)= {\rm vol}\{B\} - \bE {\rm vol} \big\{  B(Y_{t_1}) \cap B(Y_{t_2})\big\},\;\;\;
 \bar{Q}_{\{1,2\}}(t_1, t_2) =  \bE {\rm vol} \big\{  B(Y_{t_1}) \cap B(Y_{t_2})\big\}.
\end{align*}
Thus 
\begin{align*}
& \frac{1}{\lambda} 
\ln \bE \exp\big\{\theta_1 N(t_1)+\theta_2 N(t_2)\big\} 
\\
&\;\;\;= {\rm vol}\{B\}
\big[ (e^{\theta_1}-1)+ (e^{\theta_2}-1)\big]
+ \bE {\rm vol}\big\{B(Y_{t_1})\cap B(Y_{t_2})\big\} 
\big(e^{\theta_1}-1\big)\big(e^{\theta_2}-1\big).
\end{align*}
This formula should be compared with~(III.22.14) in \citeasnoun{Kac}
which was derived using different considerations.
\par 
Another interesting consequence of Theorem~\ref{th:Eexp} is stated in the next corollary. Denote 
\[
 N(t_1, \ldots, t_n) := \sum_{j\in \bZ} {\bf 1}\big\{X^{(j)}_{t_1}\in B, \ldots, X^{(j)}_{t_n}\in B\big\}. 
\]
In words, $N(t_1, \ldots, t_n)$ is the number of particles that were in $B$ at all
time instances 
$t_1, \ldots, t_n$. 
\begin{corollary}\label{cor:Nt123}
 For any $\theta\in \bR$, $n$, and $t_1, \ldots, t_n$ one has 
 \[
  \ln \bE \exp\{\theta N(t_1, \ldots, t_n)\} = \lambda (e^\theta-1) \int_{\bR^d} \bP\{x+Y_{t_1}\in B,\ldots, 
  x+Y_{t_n}\in B\} \rd x
 \]
\end{corollary}
\pr The proof is basically coincides with the one of Theorem~\ref{th:Eexp}. 
For $j\in \bZ$ define the event $A^{(j)}=\{X_{t_k}^{(j)}\in B, k=1, \ldots, n\}$, and let $\bar{A}^{(j)}$ 
be the complementary event.
We have 
\[
 \theta N(t_1, \ldots, t_n) =  \theta \sum_{j\in \bZ} {\bf 1}\big\{X_{t_k}\in B, k=1, \ldots, n\} = 
 \sum_{j\in \bZ} \theta {\bf 1}\{A^{(j)}\}. 
\]
Therefore 
\begin{align*}
 \bE \Big\{\exp\{\theta N(t_1, \ldots, t_n)\}\big|\, \Xi\Big\}=
 \bE \Big[\exp\Big\{\theta \sum_{j\in \bZ} {\bf 1}(A^{(j)}) \Big\} \Big|\, \Xi\Big]
 =
 \prod_{j\in \bZ} \bE \Big[ e^{\theta  {\bf 1}(A^{(j)})}  \Big|\, \Xi\Big] 
 \\
 = \prod_{j\in \bZ}
 \Big[(e^\theta-1) \bP\big\{A^{(j)}|\,\Xi\big\} + 1\Big]. 
\end{align*}
Then application of Campbell's formula completes the proof. 
\epr 
Corollary~\ref{cor:Nt123} shows that number of particles that are in $B$ at $n$ 
time instances $t_1,\ldots, t_n$
is a Poisson random variable with expectation $\lambda \bar{Q}_{\{1, \ldots, n\}}(t_1, \ldots, t_n)=\lambda\int_{\bR^d} \bP\{x+Y_{t_k}\in B, k=1, \ldots, n\}\rd x$. In particular, 
for any $t_1, t_2$
\begin{align*}
N(t_1, t_2) \sim {\rm Poisson}\Big(\lambda\bE {\rm vol}\{B(Y_{t_1})\cap B(Y_{t_2})\}\Big),
\;\;\bE N(t_1, t_2) = \lambda\bE {\rm vol}\{B(Y_{t_1})\cap B(Y_{t_2})\}.
\end{align*}
\subsection{Stationarity}
It follows from Theorem~\ref{th:Eexp} that Smoluchowski process $\{N(t), t\geq 0\}$ is strictly stationary 
if and only if for all $n$, $t_1,\ldots, t_n$, and $\pi\in \Pi_n$ 
 one has
 \begin{equation}\label{eq:stationary}
 \bar{Q}_\pi(t_1, \ldots, t_n):=\int_{\bR^d} Q_{\pi}(x;t_1,\ldots,t_n) \rd x =\bar{Q}_\pi (t_1+\tau, \ldots, 
 t_n+\tau),\;\;\;\forall \tau,
\end{equation}
i.e., function $\bar{Q}_\pi(t_1, \ldots, t_n)$ is invariant with respect to the shift in its arguments: for all $n$, $t_1, \ldots, t_n$ and $\tau$ 
\[
 \int_{\bR^d} \bP\big\{Y_{t_1}\in B(x), \ldots, Y_{t_n}\in B(x)\big\} \rd x =  
 \int_{\bR^d} \bP\big\{Y_{t_1+\tau}\in B(x), \ldots, Y_{t_n+\tau}\in B(x)\big\} \rd x. 
\]
An equivalent  condition of stationarity is given in 
\citeasnoun{McDunnough78}.
\par 
Condition  (\ref{eq:stationary}) holds  
for a wide class of particle displacement  processes $\{Y_t, t\geq 0\}$. The following  two 
examples are particularly important. 
\begin{itemize}
\item[(i)] Let $\{Y_t, t\geq 0\}$  be a  strictly stationary process; then 
$Q_\pi(x; t_1, \ldots, t_n)=Q_\pi(x; t_1+\tau, \ldots, t_n+\tau)$ for all $\tau\in \bR$, $x\in \bR^d$, 
$\pi\in \Pi_n$, and (\ref{eq:stationary}) holds trivially. Thus  $\{N(t), t\geq 0\}$ is  strictly
stationary.
\item[(ii)]
Let $\{Y_t, t\geq 0\}$ be a time homogeneous Markov process with transition function $P_t(x, A)$ on 
$[0, \infty)\times \bR^d\times \cB(\bR^d)$; then 
\begin{align*}
 Q_\pi (x; t_1, \ldots,  t_n)
 = \int_{B_1}\bigg[\int_{B_2}\cdots \int_{B_{n}} 
 \prod_{j=2}^n P_{t_j-t_{j-1}}(y_{j-1}, \rd y_j)\bigg] P_{t_1}( x, \rd y_1),
\end{align*}
where $B_k$ stands for $B$ if $k\in \pi$ and for $\bR^d\setminus B$ for $k\in \pi^c$.
If the integral   
\begin{equation}\label{eq:barP}
\int_{\bR^d} P_t(x, A)\rd x=:\bar{P}(A),\;\;A\in \cB(\bR^d)
\end{equation}
defines measure $\bar{P}$ independent of $t$ then  
\begin{align*}
 \bar{Q}_\pi(x)= \int_{B_1}\bigg[\int_{B_2}\cdots \int_{B_{n}} 
 \prod_{j=2}^n P_{t_j-t_{j-1}}(y_{j-1}, \rd y_j)\bigg] \bar{P}(\rd y_1),
\end{align*}
and process $\{N(t), t\geq 0\}$ is strictly stationary. 
For instance, if $P_t(x, A)=\int_A q_t(x-y)\rd y$  for some transition probability density $q_t$ 
then
(\ref{eq:barP}) holds with  
 $\bar{P}$  being the Lebesgue measure. 
 An important specific case of the discussed setting is when 
  $\{Y_t, t\geq 0\}$ is a process with independent and stationary increments. 
\end{itemize}
%
\subsection{Moments and  covariance function} 
Theorem~\ref{th:Eexp} allows us to calculate  the covariance function and the moments of the finite dimensional distributions 
of $\{N(t), t\geq 0\}$. These formulas are repeatedly used in the sequel.
\par 
For every fixed $\pi\in \Pi_n$ let  $\Pi_n(\pi):=\{\pi^\prime\in \Pi_n: \pi^\prime \supseteq \pi\}$ be the set of all supersets of $\pi$ in $\Pi_n$. Define
\begin{equation}\label{eq:H-def}
 U_\pi (t_1, \ldots, t_n) := \frac{1}{{\rm vol}(B)} \sum_{\pi^\prime\in \Pi_n(\pi)} \bar{Q}_{\pi^\prime} (t_1, \ldots, t_n),\;\;\;\pi\in\Pi_n,
\end{equation}
where $\bar{Q}_\pi (t_1, \ldots, t_n)$ is given in (\ref{eq:Q-pi0})--(\ref{eq:Q-pi}). 
We obviously have 
\begin{align}
U_\pi(t_1, &\ldots,t_n)=\frac{1}{\text{vol}\{B\}}
\sum_{\pi^\prime \in\Pi_n(\pi)}\int_{\mathbb{R}^d}Q_{\pi^\prime} (x;t_1,\ldots,t_n)\rd x
\nonumber
\\
&=\frac{1}{\text{vol}\{B\}}\int_{\mathbb{R}^d}\sum_{\pi^\prime\in\Pi_n(\pi)}
\bP\Big\{ \Big (x+Y_{t_j} \in B, j\in \pi^\prime \Big)\,\cap\, 
\Big (x+Y_{t_j}\notin B,\;
j\in (\pi^\prime)^c \Big) \Big\}\rd x\nonumber
\\
&=\frac{1}{\text{vol}\{B\}}\int_{\mathbb{R}^d}
\bP\Big\{x\in\cap_{j\in \pi} B (Y_{t_j})\Big\}\rd x
=\frac{1}{\text{vol}\{B\}}\bE\,\text{vol}\,\Big\{\cap_{j\in \pi} B(Y_{t_j})\Big\},
\label{eq:H-pi}
\end{align}
where the third equality follows 
because for every $\pi\in \Pi_n$  events 
\[
 \big\{x+Y_{t_j} \in B, j\in \pi^\prime\setminus \pi\big\}\cap 
 \big\{x+Y_{t_j}\not\in B, \;j\in (\pi^\prime)^c\big\}, \;\;\;\pi^\prime\in \Pi_n(\pi)  
\]
are disjoint and constitute a partition of the sample space. 
\par
Although $U_\pi (t_1, \ldots, t_n)$ is formally defined  
as a function of $n$ variables $t_1, \ldots, t_n$, it follows from (\ref{eq:H-pi}) that in fact 
$U_\pi(t_1, \ldots, t_n)$ 
is a function of 
$\{t_j, j\in \pi\}$ only; here   $\pi$ indicates the subset of variables on which function
$U_\pi (t_1, \ldots, t_n)$ actually depends. Therefore  with slight abuse of notation when appropriate 
we will drop subscript $\pi$ and indicate the corresponding variables explicitly, e.g., 
\[
U(t_i, t_j, t_k) := U_{\{i,j,k\}}(t_1, \ldots, t_n),\;\;\;i,j,k\in \{1, \ldots, n\}.
\]
Several remarks on the properties of function $U_\pi (t_1, \ldots, t_n)$ defined in (\ref{eq:H-def}) 
are in order.
\begin{remark}\mbox{}
\begin{itemize} 
\item[{\rm (a)}] If $\pi$ is a singleton, $\pi=\{k\}$, then it follows from (\ref{eq:H-pi}) that
$U_{\{k\}}(t_1, \ldots, t_n)=U(t_k)=1$ for all $k\in \{1, \ldots, n\}$.
\item[{\rm (b)}] $U_\pi(t_1, \ldots, t_n)$ is 
monotone non--increasing in the following sense:
for all 
$t_1, \ldots, t_n$
one has 
\[
 U_\pi (t_1, \ldots, t_n) \leq U_{\pi^\prime}(t_1, \ldots, t_n),\;\;\;\forall \pi\supseteq \pi^\prime,
 \;\;\pi, \pi^\prime \in \Pi_n.
\]
\item[{\rm (c)}] If $\{N(t), t\geq 0\}$ is strictly stationary then for all $n, t_1,\ldots, t_n$, and $\pi\in\Pi_n$
\[
 U_\pi (t_1, \ldots, t_n) = U_\pi (t_1+\tau,\ldots, t_n+\tau),\;\;\forall \tau.
\]

\end{itemize}
\end{remark}
\begin{proposition}\label{prop:moments}
 For any $t_1\leq t_2\leq \cdots \leq t_n$ one has
 \begin{equation}\label{eq:moments}
  \bE \Big[\prod_{j=1}^n N(t_j)\Big] =  \sum_{k=1}^n \rho^k \sum_{\cP(n, k)} \prod_{j=1}^k U_{\pi_j}(t_1, \ldots, t_n),
 \end{equation}
 where $\cP(n, k)$ is the set of all partitions $\pi=(\pi_1, \ldots, \pi_k)$ of $\{1, \ldots, n\}$
in $k$ subsets. In particular, for $n=2$
\begin{equation}\label{eq:cov}
 \bE \big[N(t_1)N(t_2)\big] = \rho^2 + \rho U(t_1, t_2),\;\;\;{\rm cov}\{N(t_1), N(t_2)\} = 
 \rho U(t_1, t_2).
\end{equation}
\end{proposition}
\pr The proof follows  from Theorem~\ref{th:Eexp} and  
a multivariate version of the Fa\'a~di~Bruno formula
[see, e.g., \citeasnoun{Hardy}]. Letting 
\[ 
y(\theta_1, \ldots, \theta_n)= \lambda \sum_{\pi^\prime\in \Pi_n}
 \big[e^{\sum_{k\in \pi^\prime}\theta_k} -1\big] \bar{Q}_{\pi^\prime} (t_1,\ldots, t_n)
\]
we have 
$ \psi(\theta_1, \ldots, \theta_n) :=  \bE e^{\sum_{k=1}^n \theta_k N(t_k)} = e^{y(\theta_1, \ldots,\theta_n)}$.
Therefore by the Fa\'a~di~Bruno formula 
\[
 \frac{\partial^n \psi (\theta_1, \ldots, \theta_n)}{\partial \theta_1\cdots \partial \theta_n}
 = e^{y(\theta_1, \ldots, \theta_n)} \sum_{\cP} \prod_{\pi\in \cP} \frac{\partial^{|\pi|} y(\theta_1, \ldots \theta_n)}{\prod_{l\in \pi}\partial \theta_l}, 
\]
where the sum is over all partitions $\cP$ of set $\{1, \ldots, n\}$, the product is 
over all blocks $\pi$
of partition $\cP$, and $|\pi|$ stands for the cardinality of block $\pi$.
We obviously have 
\[
 \frac{\partial^{|\pi|} y(\theta_1, \ldots \theta_n)}{\prod_{l\in \pi}\partial \theta_l} = 
 \lambda \sum_{\pi^\prime\in \Pi_n: \pi^\prime \supseteq \pi} 
 e^{\sum_{k\in \pi^\prime}\theta_k} \bar{Q}_{\pi^\prime} (t_1,\ldots, t_n),
\]
and therefore 
\begin{align*}
 \bE \prod_{j=1}^n N(t_j) = \frac{\partial^n \psi (\theta_1, \ldots, \theta_n)}{\partial \theta_1\cdots \partial \theta_n}\bigg|_{\theta_1=\cdots=\theta_n=0} = 
 \sum_{\cP} \prod_{\pi \in \cP} \rho U_\pi(t_1, \ldots, t_n) 
\\
 = \sum_{k=1}^n \rho^k \sum_{\cP(n, k)} \prod_{j=1}^k U_{\pi_j}(t_1, \ldots, t_n),
 \end{align*}
as claimed.
\epr

\begin{remark}
 Note that process $\{N(t), t\geq 0\}$ is weakly stationary if and only if 
 $U(t_1, t_2)$ depends on $t_1$ and $t_2$ 
 via difference $t_2-t_1$. 
 Because  
 ${\rm vol}\big\{B(Y_{t_1})\cap B(Y_{t_2})\big\}={\rm vol}\big\{B\cap B(Y_{t_2}-Y_{t_1})\big\}$, 
 in view of 
 (\ref{eq:H-pi}), $\{N(t), t\geq 0\}$ is weakly stationary if process 
 $\{Y_t, t\geq 0\}$ has stationary increments, i.e., $Y_{t_2}-Y_{t_1} \stackrel{d}={}Y_{t_2-t_1}$ for all 
 $0\leq t_1\leq t_2$.
\end{remark}
\begin{remark}\label{rem:covariogram}
If $\{N(t), t\geq 0\}$ is weakly stationary then 
by definition 
\[
 U(t_1, t_2)= \frac{1}{{\rm vol}(B)} \bE \,{\rm vol}\{B(Y_{t_1})\cap B(Y_{t_2})\}=
 \frac{1}{{\rm vol}(B)} \bE \,{\rm vol}\{B\cap B(Y_{t_1-t_2})\}=: H(t_1-t_2),
\]
and it follows from    (\ref{eq:cov}) that $H(t)$ is the correlation function of $\{N(t), t\geq 0\}$,
while $R(t)=\rho H(t)={\rm cov}\{N(s), N(t+s)\}$ is the covariance function. 
It is also worth noting that 
\[
R(t) = \rho H(t)= {\rm cov}\{N(s), N(t+s)\}=\lambda\, \bE g_B(Y_t),\;\;\;\;t\geq 0,
\]
where $g_B$ is the covariogram of  set $B$ [cf.~(\ref{eq:covariogram})].
\end{remark}
\par 
For weakly stationary $\{N(t), t\geq 0\}$ we have 
 \[
  \bE |N(s+t)-N(s)|^2 = 2\rho [1- H(t)].
 \]
This equation was derived by \citeasnoun{Smoluchowski14}; in his terminology  function 
$1-H(t)$ is called the {\em probability after--effect}
(Wahrscheinlichkeitnachwirkung)
[see \citeasnoun{Chandrasekhar} and \citeasnoun{Kac}]. 
In addition to  relationship to the correlation function of the Smoluchowski process, 
the probability after--effect $1-H(t)$ has another interpretation. It is 
loosely referred in~\citeasnoun{Chandrasekhar}
as  {\em a probability that a particle somewhere inside $B$ will have emerged from it during time $t$.}  
In fact 
\[
 H(t) =
 \frac{\int_{\bR^d} \bP\{x+Y_{s}\in B, x+Y_{s+t} \in  B\} \rd x}{\int_{\bR^d} 
 \bP\{x+Y_s\in B\}\rd x} = \frac{\bE N(s, s+t)}{\lambda {\rm vol}\{B\}},
\]
where $N(s, s+t)$ is the number of particles that are in $B$ at both time instances $s$ and $s+t$.
\subsection{Probability generating functional}
Theorem~\ref{th:Eexp} gives the probability generating functional of 
$\{N(t), t\geq 0\}$:
\begin{equation*}
 \bE \prod_{j=1}^n z_j^{N(t_j)}= \exp\Big\{\lambda \sum_{\pi\in \Pi_n} \Big[\prod_{k\in \pi} z_k -1\Big]
 \bar{Q}_\pi (t_1, \ldots, t_n) \Big\}.
\end{equation*}
This expression allows us to calculate finite dimensional  distributions of the Smoluchowski process.
In particular, for strictly stationary  Smoluchowski process $\{N(t), t\geq 0\}$  
\begin{align*}
\psi(z_1, z_2):=\bE \big[z_1^{N(s)} z_2^{N(s+t)}\big] = 
\exp\big\{ \rho [(z_1-1)+ (z_2-1)] + \rho H(t)(z_1-1)(z_2-1)\big\},
\end{align*}
and routine differentiation yields: for $m\geq 0$ and $k\geq 0$
\begin{align}
\bP\{N(s)=m,\, & N(s+t)=m+k\} = \frac{1}{m! (m+k)!}\Big\{\frac{\partial^m}{\partial z_1^m } 
\frac{\partial^{m+k}}{\partial z_2^{m+k}} \psi(z_1, z_2)\Big\}\Big|_{z_1=z_2=0}
\nonumber
\\*[2mm]
&= \frac{e^{-\rho} \rho^m}{m!} 
\sum_{j=0}^m \binom{m}{j} [H(t)]^j [1-H(t)]^{m-j} \frac{e^{-\rho[1-H(t)]} 
[\rho(1-H(t))]^{m+k-j}}{(m+k-j)!},
\label{eq:m+k}
\end{align}
and for $m\geq 0$ and  $0\leq k\leq m$
\begin{align}
\bP\{N(s)=m,\, & N(s+t)=m-k\} = \frac{1}{m! (m-k)!}\Big\{\frac{\partial^m}{\partial z_1^m } 
\frac{\partial^{m-k}}{\partial z_2^{m-k}} \psi(z_1, z_2)\Big\}\Big|_{z_1=z_2=0}
\nonumber
\\*[2mm]
&= \frac{e^{-\rho} \rho^m}{m!} 
\sum_{j=k}^{m} \binom{m}{j} [H(t)]^j [1-H(t)]^{m-j} \frac{e^{-\rho[1-H(t)]} 
[\rho(1-H(t))]^{j-k}}{(j-k)!}.
\label{eq:m-k}
\end{align}
Since $N(t)$ is a Poisson random variable with parameter $\rho$ for every $t$
we have 
\begin{align}
&\bP\{N(s+t)=m+k | N(s)= m\} 
\nonumber
\\
&\;\;\;=
\sum_{j=0}^m \binom{m}{j} [H(t)]^j [1-H(t)]^{m-j} \frac{e^{-\rho[1-H(t)]} 
[\rho(1-H(t))]^{m+k-j}}{(m+k-j)!},\;\;\;k\geq 0,
\label{eq:m+k-cond}
\end{align}
and 
\begin{align}
&\bP\{N(s+t)=m-k| N(s)=m\}
\nonumber
\\
&\;\;\;= 
\sum_{j=k}^{m} \binom{m}{j} [H(t)]^j [1-H(t)]^{m-j} \frac{e^{-\rho[1-H(t)]} 
[\rho(1-H(t))]^{j-k}}{(j-k)!},\;\;\;0\leq k\leq m.
\label{eq:m-k-cond}
\end{align}
\par 
Formulas (\ref{eq:m+k}) and (\ref{eq:m-k}) have been derived in \citeasnoun{Smoluchowski14} using combinatorial arguments; see also 
\citeasnoun{Chandrasekhar} and \citeasnoun{Kac}. It is worth noting that the 
expressions on the right hand sides of (\ref{eq:m+k-cond}) and (\ref{eq:m-k-cond})
are discrete convolutions of binomial and Poisson distributions. 
\subsection{Gaussian approximation}
The next statement establishes a Gaussian approximation for the finite dimensional distributions of process
$\{N(t), t\geq 0\}$.  This result is an easy consequence of  Theorem~\ref{th:Eexp}.
To state the result we need the following notation. 
We consider a family of processes $\{N_\rho(t), t\geq 0\}$
indexed by parameter $\rho=\lambda {\rm vol}(B)>0$, and for fixed $(t_1, \ldots, t_n)\in \bR^n$ let 
\[
 N^n_\rho:= \big(N_\rho(t_1), \ldots, N_\rho(t_n)\big),\;\;\; e_n=(1, \ldots, 1)\in \bR^n.
\]

\begin{proposition}\label{prop:Normal}
 For any $n\geq 1$ and $(t_1, \ldots, t_n)\in \bR^n$ one has 
 \[
  \frac{N^n_\rho - \rho e_n}{\sqrt{\rho}} \;\stackrel{d}{\to}\; \cN_n\big(0, \Sigma_H\big),\;\;\;
  \rho\to\infty,
 \]
where  $\Sigma_H$ is the $n\times n$ matrix with elements
\[
 \big[ \Sigma_H]_{i,j} := U_{\{i,j\}}(t_1, \ldots, t_n)=U(t_i, t_j),\;\;i,j=1, \ldots, n.
\]

\end{proposition}
\pr  The proof is standard; it is  based on application of Theorem~\ref{th:Eexp}.
It follows from this theorem that for any  $\theta=(\theta_1, \ldots, \theta_n)\in \bR^n$
one has 
\begin{align}\label{eq:sum-theta}
\ln\bE 
\exp\left\{\sum_{k=1}^n\theta_k\frac{N_\rho(t_k)-\rho}{\sqrt{\rho}}\right\} &=
-\sqrt{\rho}\sum_{k=1}^n\theta_k+\ln\bE \exp\left\{\sum_{k=1}^n N_\rho(t_k)\frac{\theta_k}{\sqrt{\rho}}\right\}\nonumber
\\&=-\sqrt{\rho}\sum_{k=1}^n\theta_k+ \lambda \sum_{\pi\in \Pi_n}
\big[e^{\sum_{k\in \pi}\frac{\theta_k}{\sqrt{\rho}}} -1\big] \bar{Q}_\pi(t_1, \ldots, t_n) \nonumber
\\ &= -\sqrt{\rho}\sum_{k=1}^n\theta_k 
+ \lambda \sum_{l=1}^\infty \frac{1}{l!} \sum_{\pi\in \Pi_n} 
\bigg(\sum_{k\in \pi} \frac{\theta_k}{\sqrt{\rho}}\bigg)^l \,\bar{Q}_\pi (t_1, \ldots, t_n).
\end{align}
In the second sum on the right hand side of (\ref{eq:sum-theta}) 
the term corresponding to $l=1$ is 
\begin{align}
 \lambda \sum_{\pi\in \Pi_n} \sum_{k\in \pi} \frac{\theta_k}{\sqrt{\rho}} \bar{Q}_\pi(t_1, \ldots, t_n)
 = \frac{\lambda {\rm vol}(B)}{\sqrt{\rho}} \sum_{k=1}^n \theta_k  \frac{1}{{\rm vol}(B)}\sum_{\pi \in \Pi_n(\{k\})} 
 \bar{Q}_\pi (t_1, \ldots, t_n)
 \nonumber
 \\
 = \sqrt{\rho} \sum_{k=1}^n \theta_k U(t_k) = \sqrt{\rho} \sum_{k=1}^n \theta_k,
\label{eq:l=1}
 \end{align}
where, we remind,  $\Pi_n(\pi)=\{\pi^\prime\in \Pi_n: \pi^\prime \supseteq \pi\}$, and we took into account 
the definition of $U_\pi(t_1, \ldots, t_n)$ [cf. (\ref{eq:H-pi})].
 For $l=2$ we have
 \begin{align}
  \frac{\lambda}{2} \sum_{\pi\in \Pi_n} \sum_{k_1\in \pi}\sum_{k_2\in \pi} \frac{\theta_{k_1}\theta_{k_2}}{\rho}
  \bar{Q}_\pi (t_1, \ldots, t_n) = \frac{\lambda}{2\rho} \sum_{k_1=1}^n 
  \sum_{k_2=1}^n \theta_{k_1}\theta_{k_2} \sum_{\Pi_n(\{k_1, k_2\})} \bar{Q}_\pi (t_1, \ldots, t_n)
  \nonumber 
  \\
  = \frac{1}{2} \sum_{k_1=1}^n\sum_{k_2=1}^n \theta_{k_1}\theta_{k_1} U(t_{k_1}, t_{k_2}),
  \label{eq:l=2}
 \end{align}
and similarly for $l>2$ 
\begin{align*}
  \frac{\lambda}{l!\rho^{l/2}} \sum_{\pi\in \Pi_n} 
&\sum_{k_1\in \pi}\cdots \sum_{k_l\in \pi} \theta_{k_1}\cdots \theta_{k_l}  \,\bar{Q}_\pi (t_1, \ldots, t_n)
\\
&= \frac{1}{l!\rho^{\frac{1}{2}l-1}} \sum_{k_1=1}^n \cdots \sum_{k_l=1}^n \theta_{k_1}\cdots \theta_{k_l}
U(t_{k_1}, \ldots, t_{k_l}).
\end{align*}
Taking into account that the terms with $l>2$ tend to zero as $\rho\to\infty$ and combining 
(\ref{eq:l=2}), (\ref{eq:l=1}) and (\ref{eq:sum-theta}) we obtain 
\[
 \lim_{\rho\to\infty} \bE \exp\left\{\sum_{k=1}^n\theta_k\frac{N_\rho(t_k)-\rho}{\sqrt{\rho}}\right\}=
 \exp\big\{ \tfrac{1}{2} \theta^T \Sigma_H \theta\big\}
\]
for all $t_1, \ldots, t_n$. The statement of the proposition follows. 
\epr
\begin{remark}
 Proposition~\ref{prop:Normal} demonstrates that convergence of the finite dimensional distributions 
 of Smoluchowski process $\{N(t), t\geq 0\}$ 
 to the multivariate normal distribution takes place not only if 
 the rate parameter $\lambda$
 tends to infinity, but also if the volume of the  observation region ${\rm vol}(B)$ increases without bound.  
\end{remark}

\subsection{Related stochastic models}\label{sec:related}
In this section we discuss some stochastic models  that are intimately connected with  the 
Smoluchowski processes. 

\paragraph{$M/G/\infty$ queue.}
The $M/G/\infty$ queue is 
one of the most well studied and well understood models 
in Queueing Theory [see, e.g., \citeasnoun{Takacs}]. 
In this model 
there is an infinite number of servers, 
customers arrive at time epochs $(\tau_j)_{j\in \bZ}$ of the homogeneous Poisson process of 
intensity~$\lambda$, obtain service upon arrival and leave the system after the service completion.
The service times $(s_j)_{j\in \bZ}$ 
are independent identically distributed random variables with common distribution function $S$; they are assumed to be  independent of the 
arrival process. 
\par 
If we assume that the system operates infinite time (it is in a stationary  regime) then 
the number of busy servers in the system 
at time $t$ is given by the formula 
\[
 N(t)= \sum_{j\in \bZ} {\bf 1} \{ \tau_j \leq t, \tau_j+ s_j>t\},\;\;\;t\geq 0.
\]
It has been shown in \citeasnoun{Goldenshluger16} 
that 
\begin{align*}
 \frac{1}{\rho}\ln \bE \exp\Big\{\sum_{k=1}^n \theta_k N(t_k)\Big\}= 
 \sum_{k=1}^n (e^{\theta_k}-1) + \sum_{k=1}^{n-1} H(t_k) \sum_{m=k}^{n-1}(e^{\theta_{m-k+1}}-1)
 e^{\sum_{i=m-k+2}^n\theta_i} (e^{\theta_{k+1}}-1),
\end{align*}
where 
\[
\rho:= \frac{\lambda}{\int_0^\infty [1-S(x)]\rd x},\;\;
 H(t):= \frac{\int_t^\infty [1-S(x)]\rd x}{\int_0^\infty [1-S(x)]\rd x}. 
\]
This formula can be rewritten as 
\begin{align}\label{eq:MGinfty}
 \frac{1}{\rho}\ln \bE \exp\Big\{\sum_{k=1}^n \theta_k N(t_k)\Big\}= 
 \sum_{k=1}^n  \sum_{l=k}^n (e^{\sum_{j=l}^m \theta_j}-1) \tilde{H}_{l, k},
\end{align}
where 
\[
 \tilde{H}_{l, k}:= H(|t_l-t_k|)+ H(|t_{l-1}-t_{k+1}|)- H(|t_{l-1}-t_k|)- H(|t_{l}-t_{k+1}|),\;\;
 1\leq l, k\leq n,
\]
and the term $H(|t_l-t_k|)$ is interpreted as zero if $l\not\in \{1, \ldots, n\}$ or 
$k\notin \{1, \ldots, n\}$.
It follows from (\ref{eq:MGinfty}) that 
structure of the moment generating function 
coincides with that in (\ref{eq:formula}).  Therefore the process of   
the number of busy servers in the $M/G/\infty$ queue is a
Smoluchowski process. In fact, it is a version of the Smoluchowski process
with displacement determined by 
the undeviated uniform motion (this process is studied in detail in Section~\ref{sec:undeviated}).
\par 
A particular case of the $M/G/\infty$ model is the $M/M/\infty$ queue 
when  the service time distribution is exponential. The $M/M/\infty$ model is Markovian; it  provides
an example of the Smoluchowski process with Markov property. We refer to \citeasnoun{Bingham} for further 
details on the connection between $M/G/\infty$ queue and the Smoluchowski process.
\paragraph{Bernoulli--Poisson branching process with immigration.}
Branching process with immigration is a model for  the size of a population evolving over time.  
It is a sequence $\{N_{t}, t=0,1,2,\ldots\}$ of non-negative integer 
random variables  defined by equation
\begin{equation}\label{eq:branching}
 N_{t+1}= \sum_{j=1}^{N_t} Z_j^{(t+1)} + I_{t+1},\;\;\;t=0,1,\ldots 
\end{equation}
where 
$N_t$ stands for the size of the population at time $t$, 
$Z_j^{(t+1)}$ is the number of offsprings of the $j$th individual existing in the population 
at time $t$, and 
$I_{t+1}$ is the number of immigrants joining the population at time $t+1$. 
The non--negative integer random variables 
$\{Z_j^{(t)}, j\geq 1\}$ and $\{I_t\}$, $t=1,2,\ldots$ are 
independent of each other; they are assumed to be    
 sequences of independent 
identically distributed random variables. 
\par 
The following special case of the  branching process with immigration is 
closely related to the Smoluchowski process. 
If the offspring distribution is Bernoulli with parameter $1-p$ (each individual
replaces itself with probability $1-p$ or ``dies'' with probability $p$), and 
if the immigration distribution is Poisson with parameter $\lambda$
then~(\ref{eq:branching}) implies that conditionally on $N_t=i$, random variable $N_{t+1}$
is a sum of two independent random variables: binomial random variable with parameters $i$ and $1-p$ and Poisson 
random variable with parameter $\lambda$. Therefore 
\begin{align*}
 \bP\{N_{t+1}= j | N_t=i\}= \sum_{l=0}^{\min\{i,j\}} \binom{j}{l} (1-p)^l p^{j-l} 
 \frac{e^{-\lambda}\lambda^{j-l}}{(j-l)!},
 \end{align*}
which coincides with (\ref{eq:m+k-cond}) and (\ref{eq:m-k-cond}) for a particular 
choice of parameters $\lambda$ and $p$.  This fact shows the relationship between  
the Bernoulli--Poisson branching processes with immigration and the Smoluchowski processes.  
\par
The Bernoulli--Poisson branching process with immigration is Markovian while the Smoluchowski process is in general not.
We refer to \citeasnoun{McDunnough78} where, using connections between these two 
models,  conditions for Markovianity of the Smoluchowski processes are established. 
\section{Covariance function estimation}\label{sec:covariance}
In this section and from now on we assume that   process $\{N(t), t\geq 0\}$ is strictly stationary.
We consider the problem of estimating 
the covariance function of the Smoluchowski process
from continuous time observations $\{N(t), 0\leq t\leq T\}$. 
The results established in this section 
are repeatedly used in the sequel. 
\par 
Recall that covariance and correlation functions   of $\{N(t), t\geq 0\}$ are 
given by 
\begin{align}
 R(t) & =  {\rm cov}\{N(s), N(t+s)\} = \rho H(t),
 \label{eq:covariance-0}
 \\
 H(t) &  = U(t+s, s) =
 \frac{1}{\text{vol}\{B\}}\bE\,\text{vol}\,\big\{ B \cap B(Y_{t})\big\} = \frac{\bE
 g_B(Y_{t})}{\text{vol}\{B\}},
 \label{eq:correlation-0}
\end{align}
and for any $k, t_1, \ldots, t_k$
\[
 U(t_1, \ldots, t_k)= \frac{1}{{\rm vol}(B)} \bE \,{\rm vol}\big\{ \cap_{j=1}^k B(Y_{t_j})\big\}.
\]
The covariance and correlation functions are extended to the entire real line 
by symmetry: $R(-t)=R(t)$, and $H(-t)=H(t)$ for every $t\in \bR$.
\par 
The standard unbiased estimators of $R(t)$ and $H(t)$ 
are 
\begin{equation}\label{eq:R-hat}
 \hat{R}(t) := \frac{1}{T-t} \int_0^{T-t} (N(s)-\rho) (N(t+s)-\rho)\rd s,\;\;\;\hat{H}(t):= \hat{R}(t)/\rho.
\end{equation}
The next statement establishes an upper bound on the 
estimation accuracy of~$\hat{R}(t)$.
\begin{theorem}\label{lem:R-Rhat}
For any $0\leq   t<T$ and $0\leq s<T$  one has 
\begin{align}
&(T-s)(T-t)\; \bE [\hat{R}(t)-R(t)][\hat{R}(s)-R(s)] 
 \nonumber
 \\
& = \rho^2 \int_0^{T-s}\int_0^{T-t} \big[ H(\tau_1-\tau_2)H(\tau_1-\tau_2+t-s)+
 H(\tau_1-\tau_2-s)H(\tau_1-\tau_2+t)\big]\rd \tau_1\rd \tau_2
\nonumber
 \\
& \hspace{60mm}\;\;\;+ \rho \int_0^{T-s}\int_0^{T-t} U(\tau_1, \tau_1+t, \tau_2, \tau_2+s)
 \rd \tau_1\rd \tau_2.
\label{eq:R-Rhat}
 \end{align}
In particular, 
\begin{align}\label{eq:cov-err}
 \bE |\hat{R}(t)-R(t)|^2 
 \leq  \frac{C(\rho^2+\rho)}{T-t} \int_0^T H(y)\rd y,
\end{align}
where $C$ is an absolute constant.
\end{theorem}
\begin{remark}
 Theorem~\ref{lem:R-Rhat} is valid for any Smoluchowski process, 
 independently of the displacement model.
\end{remark}
\par 
It is instructive to compare Theorem~\ref{lem:R-Rhat} with results on estimation of covariance 
functions of stationary Gaussian processes. It is well known that  for 
a stationary Gaussian process with correlation function $H(t)$ the standard estimator
of    $H(t)$ 
is consistent in the mean square sense if and only if  $T^{-1} \int_0^T H^2(t)\rd t\to 0$ as $T\to\infty$, and 
the mean squared estimation error admits an upper bound which is proportional to $T^{-1} \int_0^T H^2(t)\rd t$;
see, e.g., \citeasnoun[Section~17]{Yaglom}. As Theorem~\ref{lem:R-Rhat} shows,  
this is not the case for the Smoluchowski process: here the mean squared error is proportional to $T^{-1}\int_0^T H(t)\rd t$. 
The difference is due to the presence of the last term on the right hand side
of (\ref{eq:R-Rhat}) which,  in general, cannot be bounded in terms of the $\bL_2$--norm of 
the correlation function.
\section{Undeviated uniform motion}\label{sec:undeviated}
In this section  
we consider  estimation problems for  
the Smoluchowski process 
(\ref{eq:X-t})--(\ref{eq:N-t})   
associated with the undeviated uniform motion; see~(\ref{eq:undeviated}).
In this setting function 
$Q_\pi (t_1, \ldots, t_n)$ in 
(\ref{eq:Q-pi})
is non--zero only for the subsets 
$\pi$ consisting  of consecutive numbers 
from $\{1, \ldots, n\}$, i.e., for the subsets  such that   
$\pi = \{l, l+1, \ldots, m\}$ for some  
 $1\leq l\leq m\leq n$. Under these circumstances
\begin{align*}
 Q_\pi(x) &=\bP \big\{ (x+ v t_k \in B, k\in\{l,\ldots, m\}) \cap ( x+vt_k\notin B,\;
 k\notin \{l, \ldots, m\}\big\}
 \\
 &= \bP\big\{ x\in B(vt_l)\cap B(vt_m), x\notin B(vt_{l-1}), x\notin B(vt_{m+1})\big\}
 \\
 &= \bP\big\{ x\in B(vt_{l-1})\cap B(v t_{m+1})\big\} + 
 \bP\big\{ x\in B(vt_{l})\cap B(v t_{m})\big\}
 \\
 &\;\;\;\;\;\;\;- \bP\big\{ x\in B(vt_{l-1})\cap B(v t_{m})\big\}
 - \bP\big\{ x\in B(vt_{l})\cap B(v t_{m+1})\big\},
\end{align*}
where $B(vt_k)$   is interpreted as the empty set if  $k\notin \{1, \ldots, n\}$.
The resulting Smoluchowski process is  strictly stationary. Indeed, 
since
$${\rm vol}\{B(vt)\cap B(v\tau)\}= {\rm vol}\{B \cap B(v(t-\tau))\}= g_B(v(t-\tau)), \;\;\;\forall t, \tau, $$
we have for \mbox{$\pi=\{l, l+1, \ldots, m\}$} that   
\begin{align}
 \bar{Q}_\pi (t_1, &\ldots,  t_n) 
 \nonumber
 \\ 
 & = 
  \bE \Big[g_B(v(t_{m+1}-t_{l-1})) + 
 g_B(v(t_m-t_l))
 - g_B(v(t_m-t_{l-1}))
 - g_B(v(t_{m+1}-t_{l}))\Big], 
\label{eq:Q-undeviated}
 \end{align}
where $g_B(v(t_j-t_i))=0$ whenever $i\notin \{1, \ldots, n\}$ or 
$j\not\in \{1, \ldots, n\}$. 
It follows that 
\[
 U_\pi (t_1, \ldots, t_n)= \frac{1}{{\rm vol}(B)} \bE\, {\rm vol}\big\{ 
 B\big(v \max_{j\in \pi} t_j\big) \cap  B\big(v \min_{j\in \pi} t_j\big) \big\}
 = \frac{\bE g_B\big(v (\max_{j\in \pi} t_j - \min_{j\in\pi}t_j)\big)}{{\rm vol}(B)},
\]
and the covariance function of $\{N(t), t\geq 0\}$ is 
\begin{equation}\label{eq:cov-undeviated}
R(t)=\rho H(t),\;\;\;
 H(t) = \frac{1}{{\rm vol}(B)}\, \bE \,{\rm vol}\{ B\cap B(vt)\} = \frac{1}{{\rm vol}(B)}
 \int_{\bR^d} g_B(vt) \,\rd G(v),
\end{equation}
where, we recall,  $G$ is the velocity distribution function.  
\par 
The covariance function in (\ref{eq:cov-undeviated}) depends on the geometry of 
the observation region via covariogram~$g_B$. In this section and from now on we will assume that 
$B$ is the closed Euclidean ball of radius $r>0$ centered at the origin,  
\[
B=\{x: \|x\|\leq r\},\;\;\;r>0.
\]
Then the covariogram is 
 \begin{equation}\label{eq:sph-cap}
 g_{B}(x)=  {\rm vol}\{B\cap B(x)\}= {\rm vol}\{B\} 
 I\Big( \tfrac{d+1}{2}, \tfrac{1}{2}; 1-\tfrac{\|x\|^2}{4r^2}\big) {\bf 1}\{\|x\|\leq 2r\},
 \end{equation}
where  $I(x;a,b)$ is the regularized incomplete beta function
\[
 I(a,b;x):= \frac{1}{B(a,b)} \int_0^x t^{a-1} (1-t)^{b-1}\rd t,\;\;a>0, b>0,\;\;0\leq x\leq 1,
\]
and $B(a,b)$ is the beta function. 
Formula (\ref{eq:sph-cap}) is a consequence of the 
 well known expression for the volume of the spherical cap; see, e.g., \citeasnoun{Li}.
 \par
The particular form of covariogam $g_B$ in (\ref{eq:sph-cap}) has immediate implications on identifiability of certain functionals of $G$
from the Smoluchowski process data.
 In view of (\ref{eq:sph-cap}), $g_B(x)$ depends on $x$ via the Euclidean norm $\|x\|$ only, 
i.e., $g_{B}(x)$ is a ridge function, $g_B(x)=\tilde{g}_{B}(\|x\|)$.  Therefore 
formula (\ref{eq:Q-undeviated}) together with Theorem~\ref{th:Eexp}
imply that all finite dimensional distributions of the corresponding Smoluchowski process 
$\{N(t), t\geq 0\}$ depend on $v$ via  
 $\|v\|$ only. Thus, if $B$ is a Euclidean ball  then  only  speed 
distribution (the distribution of $\|v\|$) is identifiable from the count data. 
\par 
Therefore  we deal with 
the following two estimation problems: given  
 continuous time observations $\sN_T=\{N(t), 0\leq t\leq T\}$ of the 
Smoluchowski process we want to estimate: 
(a)~the mean speed $\mu :=\bE \|v\|= \int_0^\infty x\rd F(x)$; and (b)~the value 
$F(x_0)$ of 
the speed
distribution function $F$ at given point~$x_0$.
The problem of estimating the mean speed $\mu$ from discrete time observations was discussed
by \citeasnoun{Lindley} in the setting  of undeviated uniform motion on the real line. 
This paper established the relationship between the one--sided 
derivative of the correlation function at zero and the mean speed, and discussed construction of an estimator based on discrete time data,
but did not present rigorous analysis of its  accuracy. 
To the best of our knowledge, the problem of estimating the speed distribution has 
not been studied in the literature. 
\par
We adopt the minimax framework for measuring estimation accuracy. 
Let $\psi=\psi(F)$ be a functional of the speed distribution $F$ such 
as  $\mu=\bE\|v\|$ or $F(x_0)$. By an estimator $\hat{\psi}$ of $\psi(F)$
we mean any measurable function of observation $\sN_T$, and accuracy of $\hat{\psi}$ is measured
by the maximal root mean squared error 
\[
 \cR_T[\hat{\psi}; \sF]:=\sup_{F\in \sF} \Big\{\bE |\hat{\psi}-\psi(F)|^2\Big\}^{1/2}
\]
on a natural class $\sF$ of speed distribution functions $F$.
The minimax risk is defined by 
$\cR_T^*[\sF]:=\inf_{\hat{\psi}} \cR_T[\hat{\psi}; \sF]$, where $\inf$ is taken over all possible 
estimators of $\psi(F)$. The goal is to develop {\em a rate--optimal estimator} $\hat{\psi}_*$ 
of $\psi(F)$ such that 
$\cR_T[\hat{\psi}_*; \sF]\asymp \cR_T^*[\sF]$ as $T\to\infty$.
\subsection{Correlation function and its properties}
In view of (\ref{eq:cov-undeviated}) and (\ref{eq:sph-cap}), 
covariance function of the Smoluchwski process associated with the undeviated uniform motion 
is given by
\begin{align}
 R(t) = \rho H(t) &=\frac{\rho}{B(\frac{d+1}{2}, \frac{1}{2})} \int_0^\infty  \int_0^{1-\frac{t^2 x^2}{4r^2}}
 y^{\frac{1}{2}(d-1)} (1-y)^{-1/2}  {\bf 1} \{xt\leq 2r\}\rd y\,\rd F(x)
 \nonumber
 \\
 & = \frac{1}{B(\frac{d+1}{2}, \frac{1}{2})} \int_0^1 
 F\Big(\tfrac{2r\sqrt{y}}{t}\Big) (1-y)^{\frac{1}{2}(d-1)} y^{-1/2}\rd y.
\label{eq:R(t)}
 \end{align}
 The next statement establishes some useful properties of the correlation function.
\begin{lemma}\label{lem:H}
One has 
 \[
  \int_0^\infty  H(t)\rd t = \frac{2rB(\frac{d+1}{2},1)}{B(\frac{d+1}{2}, \frac{1}{2})} \int_0^\infty 
  x^{-1} \rd F(x).
 \]
Moreover, if $F$ is absolutely continuous with density $f$, and  $f(x)\leq M x^\alpha$
 for some $\alpha>-1$ and $0\leq x\leq \delta$ then  
\begin{equation}\label{eq:H(t)<=}
 H(t)\leq \frac{M}{1+\alpha} \Big(\frac{2r}{t}\Big)^{1+\alpha},
 \;\;\;\forall t \geq 2r/\delta,
\end{equation}
 and for any $T\geq 2r/\delta$
 \begin{equation}\label{eq:int-H}
  \int_0^T H(t)\rd t \leq  c_0\big(r +  M r^{1\vee (1+\alpha)}\eta_T\big),\;\;\;\;\;\;
  \eta_T:= \left\{\begin{array}{ll}
 T^{-\alpha}, & -1<\alpha<0,\\
\ln T, & \alpha=0,\\
1, & \alpha>0,
                  \end{array}
\right.
 \end{equation}
 where constant $c_0$ depends on $\alpha$ and $\delta$ only. 
\end{lemma}
%
\par 
Lemma~\ref{lem:H} demonstrates that the Smoluchowski process associated with the undeviated uniform motion 
exhibits  {\em short range dependence}
if and only if $\int_0^\infty x^{-1}\rd F(x)<\infty$. 
By short range dependence (or short memory)  
we mean  the property of integrability of the covariance function, while
absence of this property corresponds to {\em the long range dependence} (or long memory) 
[see, e.g., \citeasnoun[Chapter~6]{Samorodnitsky}].
The Smoluchowski process has long memory 
when $\int_0^\infty x^{-1} \rd F(x)=\infty$; this is the case, e.g., if 
 $F$ is absolutely continuous with respect to the Lebesgue measure
with density $f$, and $f(0)>0$. 
In particular, if $f$ is bounded and $f(0)>0$  
then 
$\int_0^T H(t)\rd t = O(\ln T)$ as $T\to\infty$. In general, 
the rate of decay of the correlation function 
$H(t)$ as $t\to\infty$ is determined by the local behavior of $F$ near zero. 
\subsection{Mean speed estimation}
Now we are in a position to define an estimator of the mean speed $\mu=\bE\|v\|$.
As mentioned above,  the covariogram $g_{B}(x)$ depends on $x$ via $\|x\|$
only, $g_{B}(x)=\tilde{g}_B(\|x\|)$, where  
\[
 \tilde{g}_B(y):= {\rm vol}\{B\}
 I\Big( \tfrac{d+1}{2}, \tfrac{1}{2}; 1-\tfrac{y^2}{4r^2}\Big) {\bf 1}\{y\leq 2r\};
\]
see (\ref{eq:sph-cap}).
Function $\tilde{g}_B(y)$ is continuous for all $y\geq 0$, monotone decreasing and 
has continuous 
derivative 
in interval $(0, 2r)$ which is 
\[
 \tilde{g}_B^\prime(y) = -\frac{{\rm vol}\{B\}}{r B(\frac{d+1}{2}, \frac{1}{2})}\Big(1-\frac{y^2}{4r^2}
 \Big)^{\frac{1}{2}(d-1)},\;\;\;y\in (0,2r).
\]
In view of (\ref{eq:cov-undeviated}) and (\ref{eq:sph-cap}) we have 
$R(t) = (\rho/{\rm vol}\{B\}) \int_0^\infty \tilde{g}(xt) \rd F(x)$, $t\geq 0$,
so that
\begin{align}\label{eq:Rprime-0}
 R^\prime (t) = - \frac{\rho}{r B(\frac{d+1}{2}, \frac{1}{2})} \int_0^{2r/t}
 \Big(1-\frac{t^2x^2}{4r^2}
 \Big)^{\frac{1}{2}(d-1)}  x\, \rd F(x),
\end{align}
and 
\begin{align*}
R^\prime(0+)=\lim_{t\downarrow 0} R^\prime(t) =  - \frac{\rho}{r B(\frac{d+1}{2}, \frac{1}{2})}\, 
\bE\|v\|.
\end{align*}
The last equation can be used as a  basis 
for constructing an estimator of the mean speed $\mu = \bE \|v\|$. 
The main idea is to estimate 
the one--sided derivative of the covariance function 
at zero $R^\prime(0+)$ using available data;
then the estimator of $\mu=\bE\|v\|$
is easily obtained from (\ref{eq:Rprime-0}). 
For one--dimensional case and discrete observations this idea 
has been discussed in \citeasnoun{Lindley}; see also \citeasnoun{Bingham}.
\par 
Let  
$K:[0,1]\to \bR$ be a kernel satisfying the following conditions
\begin{equation}\label{eq:K-Kernel}
 \int_0^1 K(x)\rd x=0,\;\;
 \int_0^1 xK(x)\rd x=1,\;\;\int_0^1 x^2 |K(x)|\rd x=:C_K <\infty.
\end{equation}
Fix real number $h>0$, 
and consider the following estimator of $\mu$:
\begin{equation}\label{eq:psi-hat}
 \hat{\mu}_{h} = - \frac{r B(\frac{d+1}{2}, \frac{1}{2})}{\rho h^2} 
 \int_0^{h} K\Big(\frac{t}{h}\Big) \hat{R}(t)\rd t,
\end{equation}
where $\hat{R}(t)$ is defined in (\ref{eq:R-hat}). The bandwidth $h$ is a design parameter
of the estimator; it will be specified in the sequel.
\begin{definition}
 Let $L>0$; we say that distribution function $F$ on $[0,\infty)$ 
 belongs to the class $\sF(L)$ if it is absolutely
 continuous with differentiable density $f$, and 
\[
  \sup_{x>0} \Big\{(1+ x^{l+2}) |f^{(l)}(x)|\Big\} \leq L,\;\;\;\;
 \int_0^\infty x^{l+2} \big|f^{(l)}(x)\big|\rd x \leq L,\;\;\;l=0,1.
\]
\end{definition}

\begin{theorem}\label{th:expected-value}
 Let $\hat{\mu}_{h_*}$ be  the estimator of $\mu=\bE \|v\|$ defined in (\ref{eq:psi-hat}) and associated with bandwidth $h_*$ satisfying
$(\ln T)^2/T \leq h_*\leq 1/(\sqrt{T}\ln T)$;
 then 
 \[
  \limsup_{T\to\infty}  \Big\{\sqrt{T}\, \cR_T[\hat{\mu}_{h_*}; \sF(L)]\Big\} \leq 
 C\Big(1+\frac{1}{\rho}\Big)^{1/2} [L(1+L)r]^{1/2},
 \]
%
where $C$ is a constant depending on $d$ only.
\end{theorem}
\begin{remark}
 Theorem~\ref{th:expected-value} demonstrates that the mean speed $\mu=\bE\|v\|$  can be estimated from count 
 data with the parametric rate: the maximal  root mean squared error of the proposed estimator over  class 
 $\sF(L)$ converges to zero  at the rate $1/\sqrt{T}$ as $T\to\infty$. 
 Therefore $\hat{\mu}_{h_*}$ is a rate--optimal estimator.  
 The definition of $\sF(L)$ requires
 existence and boundedness of the first derivative of $f$ (second derivative of~$F$) along with mild  tail
 conditions. The condition $F\in\sF(L)$ implies 
boundedness of the first and second derivatives of the correlation function $H$ and  integrability 
of the first  derivative of $H$; these properties are   essential in the theorem proof.
\end{remark}
\subsection{Estimation of the speed distribution}
In this section we deal with the problem of estimating the value $F(x_0)$ of 
 the speed distribution $F$ at  given point 
 $x_0>0$.  
 Our construction  uses 
 formula (\ref{eq:R(t)}); it will be convenient to rewrite it in the following form:
 \begin{align}\label{eq:w0}
 H(t) & = \int_0^\infty w(tx) \rd F(x),\;\;\;\;
 \\
 w(t) &:= \frac{ {\bf 1}\{t\leq 2r\}}{B(\frac{d+1}{2}, \frac{1}{2})}
 \int_0^{1- \frac{t^2}{4r^2}} y^{\frac{1}{2}(d-1)} (1-y)^{-\frac{1}{2}} \rd y,\;\;\;t\geq 0.
 \label{eq:w}
 \end{align}
 The correlation function on the left hand side can be estimated from the data 
$\sN_T=\{N(t), 0\leq t\leq T\}$. The
 distribution function $F$  is related to the correlation function via 
 integral operator   
 in (\ref{eq:w0}) that should be inverted. 
To construct the estimator we use a method based on the Laplace and Mellin transforms. 
 \paragraph{Preliminaries.}
First we  introduce notation and  recall  some 
 standard facts about the Laplace and Mellin
transforms; for details we refer to \citeasnoun{Widder}.
\par 
For generic function $g$ on $\bR$ the bilateral Laplace transform of $g$ is defined 
by
\[
 \widehat{g}(z):=  \cL[g; z] =\int_{-\infty}^\infty g(t) e^{-zt} \rd t,\;\;\;
\]
and $\widehat{g}(z)$ is an analytic function in 
the region where the integral converges.   In general, the convergence region 
is  a vertical strip in the complex plane, say,  
$\widehat{\Sigma}_g:=\{z\in \bC: \widehat{\sigma}_g^-< {\rm Re}(z)< \widehat{\sigma}_g^+\}$ for 
some $-\infty \leq \widehat{\sigma}_g^-< \widehat{\sigma}_g^+\leq \infty$.
The inverse Laplace transform is given by 
\[
 g(t)= \frac{1}{2\pi i} \int_{s-i\infty}^{s+i\infty} \widehat{g}(z) e^{zt}\rd z=
 \frac{1}{2\pi}\int_{-\infty}^\infty \widehat{g}(s+i\omega) e^{(s+i\omega)t}\rd \omega,
 \;\;\widehat{\sigma}_g^-<s<\widehat{\sigma}_g^+,
\]
where the integration is performed over any vertical line in the convergence region.
The Mellin transform of a function $g$ on $[0,\infty)$ is defined by the integral
\[
 \widetilde{g}(z) := \cM[g; z] = \int_0^\infty t^{z-1} g(t)\rd t
\]
with convergence region $\widetilde{\Sigma}_g:= 
\{z\in \bC: \widetilde{\sigma}_g^-< {\rm Re}(z)< \widetilde{\sigma}_g^+\}$. 
Then the inversion formula for the 
Mellin transform is  
\[
 g(t) = \frac{1}{2\pi i}\int_{s-i\infty}^{s+i\infty} t^{-z} \widetilde{g}(z)\rd z=
 \frac{1}{2\pi}\int_{-\infty}^\infty \widetilde{g}(s+i\omega) t^{-s-i\omega} \rd \omega,\;\;\;
 \widetilde{\sigma}_g^-< s<\widetilde{\sigma}_g^+.
\]
\par 
The following standard facts about the Mellin transform are repeatedly used in the sequel.
If $g_1$ and $g_2$ are two functions such that the integral 
$\int_0^\infty g_1(x)g_2(d)\rd x$ exists, and if the Mellin transforms $\widetilde{g}_1(1-z)$
and $\widetilde{g}_2(z)$ have a common strip of analyticity  then 
for any line $\{z: {\rm Re}(z)=c\}$ in this strip 
\[
 \int_0^\infty g_1(x) g_2(x)\rd x = \frac{1}{2\pi i}\int_{c-i\infty}^{c+i\infty} 
 \widetilde{g}_1(1-z) \widetilde{g}_2(z) \rd z. 
\]
In addition, the Parseval identity for the Mellin transform reads as 
\begin{equation}\label{eq:Parseval-Mellin}
 \int_0^\infty g^2(x) x^{2s-1}\rd x = \frac{1}{2\pi} \int_{-\infty}^\infty |\widetilde{g}(s+i\omega)|^2
 \rd \omega. 
\end{equation}
\paragraph{Estimator construction.}
Now we proceed with construction of the estimator.  
Let $K$ be a kernel satisfying the following condition. 
\begin{itemize}
\item[(K)]  Function~$K:\bR\to\bR$ be an infinitely differentiable bounded function on $\bR$ such that 
\[
 {\rm supp}(K)=[0,1],\;\;\;\int_0^1 K(y)\rd y=1, 
\]
and for given positive integer $m$
\[
 \int_0^1 K(y) y^j \rd y =0,\;\;\;j=1, \ldots, m.
\]
\end{itemize}
Condition~(K) is  standard in nonparametric estimation with kernel methods.   
\par For $0<h<1/2$ and $x_0 > 2h$ define function
\begin{equation}\label{eq:varphi}
 \varphi_{x_0,h}(t):= \int_0^t \Big[\frac{1}{h} K\Big(\frac{x}{h}\Big) - \frac{1}{xh}K\Big(\frac{\ln(x/x_0)}{h}\Big) \Big] \rd x,\;\;\;t\geq 0.
\end{equation}
The next statement demonstrates that for small $h$ function $\varphi_{x_0,h}$ 
is a smooth approximation
to the indicator function ${\bf 1}_{[0, x_0]}(\cdot)$.
\begin{lemma}\label{lem:properties-varphi}
 Let $K$ be a kernel satisfying condition~(K); then 
 $\varphi_{x_0, h}$ possesses the following properties:
 \[ 
 {\rm supp}(\varphi_{x_0,h})= [0, x_0e^h],\;\;\; \varphi_{x_0}(t)=1,\;\; h\leq t \leq x_0.
 \]
 In addition, $\widetilde{\varphi}_{x_0}$ is an entire function, and 
 \[
  \widetilde{\varphi}_{x_0,h}(z)= \frac{1}{z} \big[x_0^z \widehat{K}(-zh)-
  h^z\widetilde{K}(z+1)\big],\;\;
  \forall z\in \bC.
 \]
\end{lemma}
\par 
The next step in our construction is to define  
\begin{align}
 \psi_{x_0, h}(t) &:= \frac{1}{2\pi i} \int_{s-i\infty}^{s+i\infty} 
 \frac{\widetilde{\varphi}_{x_0,h}(1-z)}
 {\widetilde{w}(1-z)} t^{-z} \rd z
 \nonumber
 \\*[2mm]
 &=
 \frac{1}{2\pi i}
 \int_{s-i\infty}^{s+i\infty}  \frac{t^{-z} B(\frac{d+1}{2}, \frac{1}{2})}{(2r)^{1-z}
 B(\frac{d+1}{2}, 1-\frac{z}{2})} \Big[x_0^{1-z} \widehat{K}((z-1)h) - h^{1-z} \widetilde{K}(2-z)\Big]
 \rd z,\;\;\; s <1.
 \label{eq:psi}
\end{align}
The special form of the integrand 
in the second line  
on the right hand side of (\ref{eq:psi}) is a consequence  of 
the 
following formula for the Mellin transform of function $w$ defined in (\ref{eq:w}):
 \[
  \widetilde{w}(z)=\frac{(2r)^z}{z} \frac{B(\frac{d+1}{2}, \frac{z+1}{2})}{B(\frac{d+1}{2}, \frac{1}{2})},\;\;\; {\rm Re}(z)>0.
 \]
This formula is stated in Lemma~\ref{lem:Mellin-w} in Section~\ref{sec:proof-of-distr}
and proved there.
Note that the integrand in (\ref{eq:psi}) is an analytic function  in $\{z: {\rm Re}(z)<1\}$ so that 
the integration can be performed over any vertical line in this region.
It is also seen that $\psi_{x_0,h}$  is a function on $[0,\infty)$ defined by the 
inversion formula for the Mellin transform.
\par 
Our construction of estimator of $F(x_0)$ utilizes special properties of 
function $\psi_{x_0, h}(t)$ established  in the following lemma.
\begin{lemma}\label{lem:psi}
 Under assumption (K) one has  
 \begin{equation}\label{eq:abs-conv}
  \int_{-\infty}^\infty \bigg|\frac{\widetilde{\varphi}_{x_0, h}(1-s-i\omega)}{\widetilde{w}(1-s-i\omega)}\bigg|
  \rd \omega <\infty,\;\;\;\forall s<1.
 \end{equation}
 Moreover, 
if $\int_0^\infty |\psi_{x_0,h}(t)| H(t)\rd t <\infty$
then 
\begin{equation}\label{eq:lin-strategy}
 \int_0^\infty \psi_{x_0,h}(t) H(t)\rd t = \int_0^\infty \varphi_{x_0,h}(x) \rd F(x),
\end{equation}
where $\varphi_{x_0,h}$ is given by   (\ref{eq:varphi}).
\end{lemma}
\par 
The property (\ref{eq:lin-strategy}) is of crucial importance for our purposes. 
In view of Lemma~\ref{lem:properties-varphi}, the integral on the right hand side of (\ref{eq:lin-strategy}) 
approximates the value
$F(x_0)$ to be estimated. 
Therefore the main idea is to estimate
the left hand side of (\ref{eq:lin-strategy}) 
by plugging in the estimator of the correlation function. 
Specifically, define
\begin{equation}\label{eq:estimator}
 \hat{F}_h(x_0) := \int_0^{T/2} \psi_{x_0,h}(t) \hat{H}(t) \rd t,
\end{equation}
where $\hat{H}(t)$ is  an estimator of $H(t)$ defined in  (\ref{eq:R-hat}).
\paragraph{Upper bound on the risk.} 
Our current goal is to study the risk of 
the constructed estimator $\hat{F}_h(x_0)$. For this purpose we first 
define
the functional class of speed distribution functions on which 
 the risk of $\hat{F}_h(x_0)$ is assessed. 
 \begin{definition}
 Let $A>0$, $\beta>0$ be fixed real numbers. We say that 
 distribution function $F$ on $[0, \infty)$ belongs to the functional class 
$\sH_\beta (A)$ if 
 $F$ is $\ell :=\lfloor\beta\rfloor =\min\{k\in \bN \cup \{0\}: k<\beta\}$ times continuously differentiable 
 and 
 \[
\max_{k=1,\ldots, \ell}|F^{(k)}(x)|\leq A,\;\;\;\;|F^{(\ell)}(x) - F^{(\ell)}(x^\prime)|\leq A |x-x^\prime|^{\beta-\ell}, \;\;\;\;\forall x, x^\prime \in [0,\infty).
 \]
\end{definition}
In words, $\sH_\beta(A)$ is the class of all distribution functions on $[0,\infty)$ satisfying  
H\"older's condition of order~$\beta$. 
\begin{definition}
 Let $\alpha>-1$, $M>0$ and $\delta>0$. We say that distribution function $F$ on $[0,\infty)$ belongs to the functional class 
 $\sS_\alpha (M)$ if $F$ is absolutely continuous with density $f$ and 
 \[
  f(x)\leq M x^\alpha,\;\;\;\forall 0\leq x\leq \delta.
 \]
\end{definition}
\par 
The functional class $\sS_\alpha(M)$ imposes restrictions of the local behavior 
of the speed density $f$ near the origin. According to Lemma~\ref{lem:H}, this behavior is responsible
for the long/short range dependence of the Smoluchowski process associated with the undeviated uniform motion. 
Note that the definition of the functional class $\sS_\alpha(M)$ also 
involves parameter $\delta$, but we do  not 
indicate it in the notation.
\par 
Define also 
\[
 \sF_{\alpha, \beta}(A, M):= \sH_\beta(A)\cap \sS_\alpha(M).
\]
\begin{theorem}\label{th:distribution}
 Let $\hat{F}_{h_*}(x_0)$ be the estimator defined in   (\ref{eq:estimator}) and associated with
 kernel $K$ that satisfies assumption~(K) with $m>\beta+1$, and with bandwidth $h$ that 
 is set to be  
 \begin{equation}\label{eq:h-*}
 h_*= \bigg[A^{-2}(x_0^\beta+1)^{-2}\Big(1+\frac{1}{\rho}\Big) \frac{\tilde{\eta}_T \ln T}{T}\bigg]^{1/(2\beta+d+2)},
 \;\;\;\tilde{\eta}_T:= (r+ Mr^{1\vee (1+\alpha)})\eta_T,\;\;\;
\end{equation}
where $\eta_T$ is defined in (\ref{eq:int-H}).
Let 
\[
 \phi_T:= \Big[A^2(x_0^\beta+1)^2\Big]^{\frac{\beta+(d+1)/2}{2\beta+d+2}}
 \Big[\Big(1+\frac{1}{\rho}\Big) \frac{\tilde{\eta}_T\ln T}{T}\Big]^{\frac{\beta}{2\beta+d+2}}~;
\]
then 
\[
 \limsup_{T\to\infty} \Big\{\phi_T^{-1} \, \cR_T\big[\hat{F}_{h_*}(x_0); \sF_{\alpha, \beta}(A, M)\big]
 \Big\} \leq C, 
\]
%
%
where  $C$ may depend on $d, \beta$,  and $\alpha$  only. 
\end{theorem}
\begin{remark}
 Theorem~\ref{th:distribution} shows 
 that the rates at which the risk  of $\hat{F}_{h_*}(x_0)$ converges to zero are the following: 
 $(\ln T/T)^{\beta/(2\beta+d+2)}$ if $\alpha>0$;  $(\ln^2 T/T)^{\beta/(2\beta+d+2)}$ if $\alpha=0$; and 
 $(\ln T/T^{1+\alpha})^{\beta/(2\beta+d+2)}$ if $-1<\alpha <0$.  The existence of these three regimes 
 in the rate of convergence is explained by the short/long range dependence of the Smoluchowski process: 
 $\alpha>0$
 corresponds to the short memory, while $-1<\alpha\leq 0$ results in the long memory 
 with $\alpha=0$ being the boundary case.
 \end{remark}
 \par 
 We do not have a formal proof that $\hat{F}_{h_*}(x_0)$ is nearly rate--optimal
 up to a logarithmic factor; however  we conjecture that this is so. 
 Our conjecture is based on the connection  to the results obtained recently in
 a closely related statistical inverse problem. \citeasnoun{BG} considered 
 the problem  of 
 density estimation from observations with multiplicative measurement errors. This is 
 a statistical  inverse problem with integral operator 
 of  type (\ref{eq:w0}). It was shown there that 
 the achievable estimation accuracy in such problems is 
 determined by smoothness of the function to be estimated and by the ill--posedness index  $\gamma$ characterizing 
 the 
 rate of decay of the Mellin transform of function $w$ along vertical lines in the convergence
 region.  
 In particular, if the function to be estimated satisfies  the 
 H\"older condition with index $\beta$, and 
 the rate of decay of the Mellin transform is  $|\omega|^{-\gamma}$ as $|\omega|\to\infty$ 
 then  the minimax pointwise risk converges to zero at rate  
 $n^{-\beta/(2\beta+2\gamma+1)}$, where  $n$ is the sample size in
 the setting of \citeasnoun{BG}.
 As it is shown in the proof of Theorem~\ref{th:distribution},  the ill--posedness index of the inverse 
 problem in (\ref{eq:w0})--(\ref{eq:w}) is $\gamma=(d+1)/2$. With this definition of $\gamma$, 
 the rate of convergence 
 established in Theorem~\ref{th:distribution}, $(\ln T/T)^{\beta/(2\beta+ 2\gamma+1)}$,
 matches the one in \citeasnoun{BG} up to a logarithmic factor. 
 That is why  we conjecture that $\hat{F}_{h_*}(x_0)$ is nearly rate--optimal for $F(x_0)$. 

\section{Brownian displacement model}\label{sec:Brownian}
In this section we consider the Smoluchowski process associated with 
the displacement governed by the Brownian motion, see (\ref{eq:brownian}). 
Recall that in this setting
$ X_t^{(j)}=\xi_j+ Y_t^{(j)}$, $Y_t^{(j)}=\sigma W_t^{(j)}$, for $j\in \bZ$, $t\geq 0$,
where $W_t^{(j)}$ is 
the standard Brownian motion in $\bR^d$, and $\sigma>0$ is the diffusion coefficient. 
\par 
With this displacement process Theorem~\ref{th:Eexp} holds with
\begin{align*}
 Q_\pi (x) = \int_{B_1}\int_{B_2}\cdots \int_{B_p}
 \varphi(t_1; x, y_1) \varphi(t_2-t_1; y_1, y_2) \cdots \varphi(t_p-t_{p-1}; y_{p-1}, y_p)
 \rd y_1 \rd y_2 \cdots \rd y_p,
\end{align*}
where  $\pi\in \Pi_n$, $\varphi(t; x, y):=(\sqrt{2\pi \sigma^2 t})^{-1} 
\exp\{-\|x-y\|^2/(2\sigma^2 t)\}$ is the 
Gaussian kernel, and 
\[
 B_k:= \left\{\begin{array}{ll}
              B, & k\in \pi,\\
              B^c:=\bR^d\setminus B, & k\in \pi^c,
             \end{array}
\right.\;\;\;\;k\in \{1, \ldots, p\}.
\]
Therefore  function $\bar{Q}_\pi$ in (\ref{eq:Q-pi0}) is  given by 
\begin{align*}
 \bar{Q}_\pi(t_1, \ldots, t_p) = \int_{B_1}\int_{B_2}\cdots \int_{B_p}
 \varphi(t_2-t_1; y_1, y_2) \cdots \varphi(t_p-t_{p-1}; y_{p-1}, y_p)
 \rd y_1 \rd y_2 \cdots \rd y_p.
\end{align*}
This formula shows that 
the Smoluchowski process  associated with the Brownian displacement 
is strictly stationary.  As before, in this section we assume that the observation region 
is the Euclidean ball of radius $r$ centered in the origin. 
\subsection{Correlation function and its properties}
The correlation function of the Smoluchowski process 
governed by  Brownian displacement  is 
easily calculated using (\ref{eq:covariance-0}), (\ref{eq:correlation-0}) and (\ref{eq:sph-cap}): 
\begin{align*}
H(t) &= \frac{1}{{\rm vol}\{B\}} \bE {\rm vol}\, \{ B\cap B(\sigma W_t)\}
 = \bE\, I\Big( \tfrac{d+1}{2}, \tfrac{1}{2}; 1-\tfrac{\sigma^2\|W_t\|^2}{4r^2} \big)
 \\
  & = \frac{1}{B(\tfrac{d+1}{2}, \tfrac{1}{2})} \,\bE
 \int_0^1 {\bf 1}\Big\{ y\leq 1 - \frac{\sigma^2\|W_t\|^2}{4r^2}\Big\}
 y^{(d-1)/2} (1-y)^{-1/2}\rd y
 \\
 &=\frac{1}{B(\tfrac{d+1}{2}, \tfrac{1}{2})}
 \int_0^1 \bP\Big\{\eta \leq \frac{4r^2 y}{\sigma^2 t} \Big\} (1-y)^{(d-1)/2} y^{-1/2}\rd y,
\end{align*}
where $\eta$ is a random variable distributed $\chi^2$  with $d$ degrees of freedom.
If $\Gamma(s;x):= \int_x^\infty t^{s-1} e^{-t} \rd t$ 
is the upper incomplete Gamma function [see \cite[Chapter~6]{Abr-Ste}] then 
\[
 \bP\big\{\eta > x \big\} = \frac{1}{2^{d/2} \Gamma(d/2)} \int_x^\infty t^{d/2-1} e^{-t/2} \rd t =
 \frac{\Gamma(\tfrac{d}{2}; \tfrac{x}{2})}{\Gamma(\tfrac{d}{2})},
\]
and  we obtain
\begin{align}
 H(t) = 1- \frac{1}{B(\tfrac{d+1}{2}, \tfrac{1}{2}) \Gamma(\tfrac{d}{2})}
 \int_0^1  \Gamma\bigg(\frac{d}{2}; \frac{2r^2 y}{\sigma^2 t}\bigg) (1-y)^{(d-1)/2} y^{-1/2}\rd y
 \nonumber
 \\
  = \frac{1}{B(\tfrac{d+1}{2}, \tfrac{1}{2}) \Gamma(\tfrac{d}{2})}
 \int_0^1  \gamma\bigg(\frac{d}{2}; \frac{2r^2 y}{\sigma^2 t}\bigg) (1-y)^{(d-1)/2} y^{-1/2}\rd y,
\label{eq:covariance}
 \end{align}
where $\gamma(s;x):=\Gamma(x)-\Gamma(s;x)$ is the lower incomplete Gamma function.
\par 
In the next statement 
we summarize some properties of the correlation function that are useful for our purposes.
\begin{lemma}\label{lem:cov-properties}
\begin{itemize} 
\item[{\rm (a)}] The following asymptotic relationships hold:
  \begin{eqnarray*}
   1-H(t)&\sim& C(d)  \frac{\sigma\sqrt{t}}{2\pi r},\;\;C(d):=\left\{\begin{array}{ll}
                                                                        1, & d=1,\\
                                                                        d/(d-1), & d\geq 2,
                                                                       \end{array}
                                                                       \right.
 \;\;\;\;\;t\to 0;
 \\
 H(t)  &\sim & \frac{\Gamma(\frac{d+1}{2})}{\Gamma(d+1)  \Gamma(\frac{1}{2})} \bigg(\frac{\sqrt{2}r}{\sigma\sqrt{t}}\bigg)^{d},\;\;\;t\to\infty,
  \end{eqnarray*}
where $a\sim b$ means that $\lim(a/b)=1$.
\item[{\rm (b)}]  For every $t>0$  
\begin{align}\label{eq:H(t)leq}
 H(t) \leq  \frac{2^{d/2}\Gamma(\frac{d+1}{2})}{\Gamma(d+1) \Gamma(\frac{1}{2})}\bigg(\frac{r^2}{\sigma^2 t}\bigg)^{d/2}.
\end{align}
In addition, for $T>0$
\begin{align}\label{eq:H-int}
 \int_0^T H(t) \rd t \leq 
 \left\{\begin{array}{ll}
         \frac{4r^2}{\sigma^2(d^2-4)}, & d>2,\\*[2mm]
         1+ \frac{r^2}{2\sigma^2} \ln T, & d=2,\\*[2mm]
         \frac{2\sqrt{2}r}{\sqrt{\pi}\sigma} \sqrt{T}, & d=1.
        \end{array}
\right. 
\end{align}
\end{itemize}
\end{lemma}
\begin{remark}
 Lemma~\ref{lem:cov-properties} demonstrates  
 that  the one--sided first derivative of correlation function $H$ at zero is infinite, i.e., $H^\prime(0+)=\infty$. 
 Also, the rate of decay of $H$ at infinity depends on dimension $d$: we have 
 $H(t) = O(t^{-d/2})$ as
 $t\to\infty$.
 This implies that  for $d=1,2$ 
 the Smoluchowski
 process has long memory: in these cases $\int_0^\infty H(t)\rd t=\infty$. 
\end{remark}
\subsection{Estimation of the diffusion coefficient}
In this section we consider the problem of estimating 
diffusion coefficient $\sigma$ from  the data
$\sN_T=\{N(t), 0\leq t\leq T\}$. Although the problem of estimating the diffusion coefficient 
from observations of the Smoluchowski process was discussed in the literature [see, e.g.,  
in \citeasnoun{Kac} and \citeasnoun{Bingham}], we are not aware of specific estimators with provable 
accuracy guarantees. The goal of this section is to develop such an estimator.  
\par 
The proposed estimator is based on the following simple idea. 
Define
\[
J(t):=\frac{1}{B(\tfrac{d+1}{2}, \tfrac{1}{2}) \Gamma(\tfrac{d}{2})}
 \int_0^1  \gamma\Big(\frac{d}{2}; \frac{2r^2y}{t}\Big) (1-y)^{(d-1)/2} y^{-1/2}\rd y,\;\;\;t\geq 0.
\]
Note that $J$ is a fixed known function; it is completely determined by known parameters 
$r$ and $d$ and can be computed at any point. Observe that 
$H(t)=J(\sigma^2t)$, i.e.,  
the problem of estimating $\sigma^2$ is the problem of estimating the scale parameter of 
function $J$. 
\par 
Let $\alpha\in (0,\tfrac{1}{2})$ be a parameter to be  specified, 
and consider  the following functional of the correlation function
\[
 \Psi_\alpha:=\Psi_\alpha(H) = \int_0^\infty \frac{H(t)}{t^{1-\alpha}} \rd t.
\]
It follows from (\ref{eq:H(t)leq}) that 
the integral on the right hand side is finite for all $d\in \bN$. Note also that 
\begin{align*}
 &\Psi_\alpha= \int_0^\infty \frac{J(\sigma^2 t)}{t^{1-\alpha}}\rd t = \sigma^{-2\alpha} 
 J_\alpha,\;\;
 \end{align*}
 where 
 \begin{align*}
 J_\alpha:=\Psi_\alpha(J):= \int_0^\infty \frac{J(t)}{t^{1-\alpha}}\rd t= 
 \frac{(2r^2)^\alpha\Gamma(\frac{d}{2}-\alpha)B(\frac{1}{2}+\alpha, \frac{d+1}{2})}
 {\alpha\Gamma(\frac{d}{2})B(\frac{d+1}{2},\frac{1}{2})}~.
\end{align*}
With the introduced  notation
$\sigma^{2}= [J_{\alpha}/\Psi_\alpha]^{1/\alpha}$, and  
if  $\hat{\Psi}_\alpha$ is an estimator of $\Psi_\alpha$ then a natural estimator 
of $\sigma^{2}$ can be defined as 
$\hat{\sigma}^2:= [J_{\alpha}/\hat{\Psi}_\alpha]^{1/\alpha}$.
\par 
We consider the following estimator of $\Psi_\alpha$. 
Recall that  $\rho =\lambda {\rm vol}(B)$, and let  
\begin{align*}
\hat{H}(t) := \tfrac{1}{\rho} \{\hat{R}(t)\}_+,\;\;\;\;
 \hat{R}(t) := \frac{1}{T-t}
 \int_0^{T-t} [N(s)- \rho] [N(s+t)- \rho] \rd s,
\end{align*}
where $\{\cdot\}_+=\max\{\cdot, 0\}$.
For parameter $b>0$ to be specified define 
\begin{equation}\label{eq:Psi-est}
 \hat{\Psi}_{\alpha, b}:= \int_0^b \frac{\hat{H}(t)}{t^{1-\alpha}} \rd t,
\end{equation}
and the corresponding estimator of $\sigma^2$ is 
\[
 \hat{\sigma}^2_{\alpha, b} := \bigg(\frac{J_{\alpha}}{\hat{\Psi}_{\alpha, b}}\bigg)^{1/\alpha}.
\]
The proposed estimator $\hat{\sigma}^2_{\alpha, b}$ depends on two design 
parameters $\alpha$ and $b$ that 
are specified in the sequel.
\par 
Now we are in a position to present a result on 
accuracy of the proposed estimator of $\sigma^2$.
For $x, y>0$ we let
$\Delta(x,y):= |x-y|/(x+y)$ and note that $\Delta (\cdot, \cdot)$ defines 
a distance on 
$\{x\in \bR: x>0\}$. We use $\Delta (\cdot, \cdot)$ as a loss function in the problem of estimating 
$\sigma^2$; it measures the relative estimation accuracy. 
\par 
\begin{theorem}\label{th:sigma-1}
 Let $\hat{\sigma}^2_*= [J_{\alpha_*}/\hat{\Psi}_{\alpha_*,b_*}]^{1/\alpha_*}$
 be the estimator of $\sigma^2$, where $\hat{\Psi}_{\alpha_*, b_*}$ is the estimator in  (\ref{eq:Psi-est})
 associated with 
 \[
  \alpha_*:=\frac{1}{\ln T},\;\;\;\;\; b_*:=\left\{\begin{array}{ll}
                                             (T/\ln^2 T)^{1/d}, & d>2,\\
                                             (T/\ln^3 T)^{1/2}, & d=2,\\
                                             \sqrt{T}/\ln^2T, & d=1,
                                            \end{array}\right.                                                                             
 \]
Then   
\[
 \limsup_{T\to\infty}\sup_{\sigma^2: \sigma^2\leq r^2T} \bigg\{  
 \Big[\Big(1+\frac{1}{\rho}\Big)\Big(\frac{r}{\sigma}\Big)^{d\wedge 2} + \Big(\frac{r}{\sigma}
\Big)^{2d}\Big]^{-1}
 \phi_T^{-1}\;
\bE \big[\Delta (\hat{\sigma}_*^2, \sigma^2)\big]^2  \bigg\}
\leq C ,\;\;\;
\]
where $C$ is a constant depending on $d$ only, and 
\[
 \phi_T:=\left\{\begin{array}{ll}
                                             \ln^2 T/T, & d>2,\\
                                             \ln^3 T/T, & d=2,\\
                                             \ln^2T/\sqrt{T}, & d=1.
                                            \end{array}\right. 
\]
\end{theorem}
\begin{remark}
 Theorem~\ref{th:sigma-1} establishes an upper bound on the relative risk of
 $\hat{\sigma}^2$.
 The main advantage in using the loss function $\Delta$   
 is that the derived upper bound holds for a broad range of 
 possible values of~$\sigma$, and the estimator does 
 not require any prior information on~$\sigma$.  
 \end{remark}
 \par 
 Theorem~\ref{th:sigma-1} demonstrates that  the rate of convergence of the squared relative 
 risk coincides with the best achievable rate in estimation of the 
 correlation function $H$ up 
 to a logarithmic factor. 
In the case $d\geq 2$ the risk converges at the rate that is within logarithmic factors of the 
parametric rate; thus the estimator is nearly rate--optimal.
 The slower rate of convergence  for $d=1$ is a consequence of 
 the long range dependence. In fact, the correlation
 function is  non--integrable in the case $d=2$ too; here, however, this leads to
  an extra 
 logarithmic factor in the upper bound. 
 We conjecture that  logarithmic factors in $T$ 
 appearing in the upper bounds in the cases $d\ne 2$ 
 can be eliminated, and in the case $d=2$ the degree of the logarithmic factor can be improved. 

\par 
In the specific case of $d=1$ we are able to develop another estimator whose squared 
 risk converges to zero at the rate $\ln T/\sqrt{T}$, but this estimator 
 requires prior information on $\sigma^2$. 
The estimator construction is based on 
the first order approximation
of the correlation function $H(t)$ near zero that is established in the following lemma.
\begin{lemma}\label{lem:d=1}
 Let $d=1$ then for any $t\leq \sigma^2/(2r^2)$ one has 
\[
 1-H(t)= \frac{\sigma\sqrt{t}}{\sqrt{2\pi} r}\Big[1- \exp\Big\{-\frac{2r^2}{\sigma^2 t}\Big\} \Big] + 
 \delta(t),\;\;\;
|\delta(t)| \leq  \exp\Big\{-\frac{2r^2}{\sigma^2 t}\Big\}.
\]
\end{lemma}
\par 
Lemma~\ref{lem:d=1} suggests the following construction of  an estimator of 
 $\sigma^2$. 
Let $\tau>0$ be a parameter to be specified and define
\begin{equation}\label{eq:sigma-tau}
 \hat{\sigma}_\tau := \frac{\sqrt{2\pi} r}{\sqrt{\tau}} [1-\hat{H}(\tau)],
\end{equation}
where $\hat{H}$ is the standard  estimator of the correlations function.
\begin{theorem}\label{th:sigma-2}
 Let $\hat{\sigma}_{\tau_*}$ be the estimator (\ref{eq:sigma-tau}) associated with 
 $\tau=\tau_*= 4r^2/(\bar{\sigma}^2\ln T)$,
where $\bar{\sigma}$ is a constant. Then for every $\sigma\leq \bar{\sigma}$ and sufficiently large $T$ one has  
\begin{equation}\label{eq:sigma-tau-upper-bound}
 \bE \big|\hat{\sigma}_{\tau_*} - \sigma\big|^2
 \leq 
 c\Big(1+\frac{1}{\rho}\Big)\Big(\frac{r}{\sigma}\Big) \frac{\bar{\sigma}^2\ln T}{\sqrt{T}},
\end{equation}
where $c$ is an absolute constant.
\end{theorem}

\begin{remark}\mbox{}
\begin{itemize}
\item[{\rm (i)}] If the upper bound $\bar{\sigma}$ on parameter $\sigma$ is known then the squared 
risk of estimator $\hat{\sigma}_{\tau_*}$ converges to zero at the rate $\ln T/\sqrt{T}$
as $T\to\infty$. This can be compared with the result of Theorem~\ref{th:sigma-1} which establishes 
the rate $(\ln T)^2/\sqrt{T}$. The derived upper bound 
(\ref{eq:sigma-tau-upper-bound}) depends on $\bar{\sigma}$, and the accuracy may be poor under 
conservative choice of $\bar{\sigma}$. In contrast, the estimator $\hat{\sigma}_*^2$ does not require
any prior information on $\sigma^2$.
\item[{\rm (ii)}]
 In the case $d\geq 2$ the first order approximation of $1-H(t)$ 
 near zero is much less accurate, and  the risk 
 of the corresponding estimator  
 is much worse than the one established in Theorem~\ref{th:sigma-1}.
\end{itemize}
 \end{remark}
\section{Proofs for Section~\ref{sec:covariance}}
\label{sec:proofs-cov}
\subsection{Proof of Theorem~\ref{lem:R-Rhat}}\label{sec:proof-of-cov}
Recall that
\[
 \hat{R}(t)= \hat{r}(t) - \frac{\rho}{T-t}\int_0^{T-t} \big[N(\tau)+N(\tau+t)\big]\rd \tau +
 \rho^2,\;\;\;
 R(t)=r(t)-\rho^2=\rho H(t),
\]
where we denoted   $r(t):= \bE [N(s)N(s+t)]$, and 
 $\hat{r}(t) := (1/(T-t)) \int_0^{T-t} N(s) N(t+s)\rd s$.
We have 
\begin{align*}
 \bE \hat{R}(t)\hat{R}(s) 
= & \bE [\hat{r}(t)\hat{r}(s)] 
\\
&\;- \bE \hat{r}(t)\bigg[\frac{\rho}{T-s}\int_0^{T-s} [N(\tau_2)+
N(\tau_2+s)]\rd \tau_2 -\rho^2\bigg]
\\
&\; - \bE \hat{r}(s)\bigg[\frac{\rho}{T-t}\int_0^{T-t} [N(\tau_1)+
N(\tau_1+t)]\rd \tau_1 -\rho^2\bigg]
\\
&\;+   \bE \bigg[\frac{\rho}{T-s}\int_0^{T-s} [N(\tau_2)+
N(\tau_2+s)]\rd \tau_2 -\rho^2\bigg]
\bigg[\frac{\rho}{T-t}\int_0^{T-t} [N(\tau_1)+
N(\tau_1+t)]\rd \tau_1 -\rho^2\bigg]
\\
=: & \;J_1 - J_2 -J_3 + J_4.
 \end{align*}
Our current goal is to calculate the terms on the right hand side of the previous formula. 
\par
First we compute $J_1$.  Using (\ref{eq:moments}) with $n=4$ we have  
\begin{align*}
 &J_1:=\bE \hat{r}(t)\hat{r}(s) = \frac{1}{(T-t)(T-s)}\iint
\bE[ N(\tau_1) N(\tau_1+t)N(\tau_2) N(\tau_2+s)] \rd \tau_1\rd \tau_2
\\
&=\rho^4  + \rho^3[H(t)+H(s)] + \rho^2 H(t) H(s)
\\
&\;\;\;+ \frac{\rho^3}{(T-t)(T-s)}\iint \big[H(\tau_1-\tau_2)+H(\tau_1-\tau_2-s)+H(\tau_1-\tau_2+t)+
H(\tau_1-\tau_2+t-s)\big]\rd \tau_1 \rd \tau_2
\\
&\;\;\;+ \frac{\rho^2}{(T-s)(T-t)}
\iint \Big[U(\tau_1, \tau_1+t, \tau_2)+ U(\tau_1, \tau_1+t, \tau_2+s)
+ U(\tau_1, \tau_2, \tau_2+s)
\\
& \hspace{10mm}\;\;
+ U(\tau_1+t, \tau_2, \tau_2+s) +  H(\tau_1-\tau_2)H(\tau_1-\tau_2+t-s)+
H(\tau_1-\tau_2-s)H(\tau_1-\tau_2+t)\Big]\rd \tau_1 \rd \tau_2
\\
&\;\;\;+ \frac{\rho}{(T-t)(T-s)}\iint U(\tau_1, \tau_1+t, \tau_2, \tau_2+s) \rd \tau_1\rd \tau_2,
\end{align*}
where in the above expression 
the integrals over $\tau_1$ are taken from $0$ to $T-t$ and integrals over $\tau_2$ 
are taken from $0$ to $T-s$.
\par
Now we compute $J_2$ and $J_3$. Using 
(\ref{eq:moments}) with $n=3$  
we obtain  
\begin{align*}
 &\bE [\hat{r}(t) N(\tau_2)]= \frac{1}{T-t}\int_0^{T-t} \bE [N(\tau_1) N(\tau_1+t) N(\tau_2)] \rd \tau_1
 \\
 &= \rho^3 + \rho^2 H(t) + \frac{1}{T-t}\int_0^{T-t}\Big\{\rho^2 H(\tau_1-\tau_2)+
 \rho^2 H(\tau_1+t-\tau_2) +\rho U(\tau_1, \tau_1+t, \tau_2)\Big\}\rd \tau_1,
 \\
 & \bE [\hat{r}(t) N(\tau_2+s)]=
 \frac{1}{T-t}\int_0^{T-t} \bE [N(\tau_1) N(\tau_1+t) N(\tau_2+s)] \rd \tau_1
 \\
 &= \rho^3 + \rho^2 H(t) + \frac{1}{T-t}\int_0^{T-t} \Big\{\rho^2H(\tau_1-\tau_2-s) +
 \rho^2 H (\tau_1+t-\tau_2-s) + \rho U(\tau_1, \tau_1+t, \tau_2+s)\Big\}\rd \tau_1,
 \end{align*}
and similarly
\begin{align*}
 &\bE [\hat{r}(s) N(\tau_1)]= \frac{1}{T-s}\int_0^{T-s} \bE [N(\tau_2) N(\tau_2+s) N(\tau_1)] \rd \tau_2
 \\
 &= \rho^3 + \rho^2 H(s) + \frac{1}{T-s}\int_0^{T-s}\Big\{\rho^2 H(\tau_2-\tau_1)+
 \rho^2 H(\tau_2+s-\tau_1) +\rho U(\tau_2, \tau_2+s, \tau_1)\Big\}\rd \tau_2,
 \\
 & \bE [\hat{r}(s) N(\tau_1+t)]=
 \frac{1}{T-s}\int_0^{T-s} \bE [N(\tau_2) N(\tau_2+s) N(\tau_1+t)] \rd \tau_2
 \\
 &= \rho^3 + \rho^2 H(s) + \frac{1}{T-s}\int_0^{T-s} \Big\{\rho^2H(\tau_2-\tau_1-t) +
 \rho^2 H (\tau_2+s-\tau_1-t) + \rho U(\tau_2, \tau_2+s, \tau_1+t)\Big\}\rd \tau_2.
 \end{align*}
Therefore
\begin{align*}
 &J_2 =  \rho^4+ \rho^3 H(t) 
 \\
 &+ 
 \frac{\rho^3}{(T-t)(T-s)} \iint \big[ H(\tau_1-\tau_2)+H(\tau_1-\tau_2-s) + 
 H(\tau_1-\tau_2+t)+H(\tau_1-\tau_2+t-s)\big]\rd \tau_1\rd \tau_2
 \\
 &+ \frac{\rho^2}{(T-t)(T-s)}\iint 
 \big[U(\tau_1,\tau_1+t, \tau_2) + U(\tau_1, \tau_1+t, \tau_2+s)\big]\Big\}\rd \tau_1
 \rd \tau_2,   
\end{align*}
and 
\begin{align*}
 &J_3 =  \rho^4+ \rho^3 H(s) 
 \\
 &+ 
 \frac{\rho^3}{(T-t)(T-s)} \iint \big[ H(\tau_2-\tau_1)+H(\tau_2-\tau_1+s) + 
 H(\tau_2-\tau_1-t)+H(\tau_2-\tau_1+s-t)\big]\rd \tau_1\rd \tau_2
 \\
 &+ \frac{\rho^2}{(T-t)(T-s)}\iint 
 \big[U(\tau_2,\tau_2+s, \tau_1) + U(\tau_2, \tau_2+s, \tau_1+t)\big]\Big\}\rd \tau_1
 \rd \tau_2.  
\end{align*}
 \par 
 Now we compute $J_4$:
 \begin{align*}
  J_4 & = \frac{\rho^2}{(T-s)(T-t)}\iint \bE\big[ (N(\tau_2)+N(\tau_2+s))(N(\tau_1)+N(\tau_1+t))\big]
  \rd \tau_1\rd \tau_2 -  3\rho^4
  \\
  & = \rho^4 + \frac{\rho^3}{(T-t)(T-s)}\iint \Big[ H(\tau_1-\tau_2)+ H(\tau_1-\tau_2+t)+
  H(\tau_1-\tau_2-s) + H(\tau_1-\tau_2+t-s)\Big]\rd \tau_1 \rd \tau_2.
 \end{align*}
Now we combine expressions for $J_1, J_2, J_3$ and $J_4$ to get 
\begin{align*}
 &\bE\hat{R}(t) \hat{R}(s) = J_1-J_2-J_3+J_4
 \\
 &= \rho^2 H(t)H(s) + \frac{\rho}{(T-t)(T-s)}\iint U (\tau_1, \tau_1+t, \tau_2, \tau_2+s) \rd \tau_1\rd \tau_2
 \\
 & \;\;\;\;+ 
 \frac{\rho^2}{(T-s)(T-t)}\iint 
 \big[H(\tau_1-\tau_2) H(\tau_1-\tau_2+t-s) + H(\tau_1-\tau_2-s) H(\tau_1-\tau_2+t)\big]\rd \tau_1
 \rd \tau_2, 
\end{align*}
and finally taking into account that $R(t)R(s)=\rho^2 H(t) H(s)$ we obtain
\begin{align*}
 &\bE [\hat{R}(t)-R(t)] [\hat{R}(s)- R(s)] 
 \\
 &=  \frac{\rho}{(T-t)(T-s)}\int_0^{T-t}\int_0^{T-s} U (\tau_1, \tau_1+t, \tau_2, \tau_2+s) \rd \tau_1\rd \tau_2
 \\
 & \;\;\;\;+ 
 \frac{\rho^2}{(T-s)(T-t)}\int_0^{T-t} \int_0^{T-s} 
 \big[H(\tau_1-\tau_2) H(\tau_1-\tau_2+t-s) + H(\tau_1-\tau_2-s) H(\tau_1-\tau_2+t)\big]\rd \tau_1
 \rd \tau_2.
\end{align*}
This completes the proof.
\epr

\section{Proofs for Section~\ref{sec:undeviated}}
\label{sec:proofs-undeviated}
\subsection{Proof of Lemma~\ref{lem:H}}
We have 
\begin{align}
 H(t) & =\frac{1}{B(\frac{d+1}{2},\frac{1}{2})}\int_0^\infty 
 \int_0^1 {\bf 1}\Big\{t \leq \frac{2r\sqrt{1-y}}{x}\Big\} y^{\frac{1}{2}(d-1)}(1-y)^{-1/2}\rd y\rd F(x),
\label{eq:1}
 \\
&=
 \frac{1}{B(\frac{d+1}{2},\frac{1}{2})} 
 \int_0^1 F\Big(\frac{2r\sqrt{1-y}}{t}\Big) y^{\frac{1}{2}(d-1)}(1-y)^{-1/2}\rd y.
 \label{eq:2}
 \end{align}
It follows from (\ref{eq:1}) that 
\[
 \int_0^\infty H(t)\rd t = \frac{1}{B(\frac{d+1}{2},\frac{1}{2})}
 \int_0^\infty \int_0^1 \frac{2r}{x} y^{\frac{1}{2}(d-1)} \rd y\rd F(x)= 
 \frac{2rB(\frac{d+1}{2},1)}{B(\frac{d+1}{2}, \frac{1}{2})} \int_0^\infty \frac{\rd F(x)}{x}
\]
If $f(x)\leq Mx^\alpha$ for $0\leq x\leq \delta$ then for $t\geq 2r/\delta$ 
\begin{align*}
 H(t) \leq \frac{1}{B(\frac{d+1}{2}, \frac{1}{2})} \int_0^{2r/t} f(x) \int_0^{1} y^{(d-1)/2}(1-y)^{-1/2}
 \rd y \rd x \leq \frac{M}{1+\alpha} \Big(\frac{2r}{t}\Big)^{1+\alpha}
\end{align*}
Inequality (\ref{eq:int-H}) is obtained by integration of the above upper bound on $H(t)$.
\epr

\subsection{Proof of Theorem~\ref{th:expected-value}}

We begin with a lemma that establishes bounds on the derivatives of the correlation function.
\begin{lemma}\label{lem:H-diff}
Let $L>0$ be a real number, and assume that speed distribution $F$ is absolutely continuous with density~$f$.
\begin{itemize}
 \item[{\rm (i)}]  If $\sup_x f(x)\leq L$ and $\int_0^\infty xf(x)\rd x \leq L$  then for any $t>0$
 \begin{equation}\label{eq:Hprime-bound}
  |H^\prime(t)| \leq \frac{L}{rB(\frac{d+1}{2},\frac{1}{2})} \bigg\{1 \;\wedge\; 
  \frac{4r^2}{(d+1)t^2}\bigg\},
 \end{equation}
 and 
 \begin{equation}\label{eq:Hprime-int}
  \int_0^\infty |H^\prime(t)|\rd t \leq  \frac{L}{B(\frac{d+1}{2},\frac{1}{2})}
 \big[1+ 2B(\tfrac{d+1}{2},1)
 \big].
 \end{equation}
\item[{\rm (ii)}] If $f$ is differentiable and $\sup_x |f^\prime(x)| \vee \sup_x f(x)\leq L$ then 
\begin{align*}
  |H^{\prime\prime}(t)|\leq 
  4Lr \bigg(\frac{(d+2)r}{t^4} + \frac{1}{t^3}\bigg),\;\;\;\forall t>0.
 \end{align*}
\item[{\rm (iii)}]
 If $f$ is differentiable and 
 \begin{equation}\label{eq:L-bound}
 \sup_x x^3 |f^\prime(x)|\leq L,\;\; 
 \int_0^\infty x^3 |f^\prime(x)|\rd x \leq L, \;\;\int_0^\infty x^2 f(x) \rd x \leq L
 \end{equation}
 then
  \begin{align}\label{eq:Hprime-prime}
 |H^{\prime\prime}(t)|\leq \frac{3dL}{r^2 B(\frac{d+1}{2},\frac{1}{2})},\;\;\forall t>0,
\end{align}
and   
 \[
  \int_0^\infty |H^{\prime\prime}(t)|\rd t \leq \frac{L}{r}\bigg( \frac{3d}{B(\frac{d+1}{2}, \frac{1}{2})} + 
  2(2d+3)\bigg).
 \]
\end{itemize}
 \end{lemma}
\pr 
(i).~Differentiating (\ref{eq:R(t)}) we obtain 
\begin{eqnarray*}
  H^\prime(t) &=& - \frac{2r}{B(\frac{d+1}{2}, \frac{1}{2})} \int_0^1 
 f\Big(\frac{2r\sqrt{y}}{t}\Big) \frac{1}{t^2} (1-y)^{\frac{1}{2}(d-1)} \rd y
 \\
 &=&
- \frac{1}{rB(\frac{d+1}{2}, \frac{1}{2})} \int_0^{2r/t} xf(x) 
\Big(1-\frac{t^2x^2}{4r^2}\Big)^{\frac{1}{2}(d-1)}\rd x.
\end{eqnarray*}
If $\sup_x f(x)\leq L$ then the first equality yields for every $d=1,2,\ldots$
\[
 |H^\prime (t)| \leq \frac{2rL}{B(\frac{d+1}{2},\frac{1}{2})t^2} \int_0^1(1-y)^{\frac{1}{2}(d-1)}\rd y=
 \frac{2rL B(\frac{d+1}{2},1)}{B(\frac{d+1}{2},\frac{1}{2}) t^2},\;\;\;\forall t>0.
\]
On the other hand, it follows from the second equality that 
\[
 |H^\prime (t)|\leq  \frac{1}{rB(\frac{d+1}{2}, \frac{1}{2})} \int_0^\infty x f(x)\rd x
 \leq \frac{L}{rB(\frac{d+1}{2}, \frac{1}{2})},\;\;\;\forall t>0. 
\]
Combining these two inequalities we come to (\ref{eq:Hprime-bound}).
Moreover, 
\begin{align*}
 \int_0^\infty |H^\prime(t)|\rd t \leq \frac{L}{B(\frac{d+1}{2},\frac{1}{2})} + 
 \frac{2rLB(\frac{d+1}{2},1)}{B(\frac{d+1}{2},\frac{1}{2})}\int_r^\infty\frac{\rd t}{t^2}
 \leq \frac{L}{B(\frac{d+1}{2},\frac{1}{2})}
 \bigg[1+\frac{2B(\frac{d+1}{2},1)}{B(\frac{d+1}{2},\frac{1}{2})}
 \bigg].
\end{align*}
\par\medskip
(ii). The second derivative of $H$ is 
\begin{eqnarray*}
  H^{\prime\prime}(t) &=&  \frac{2r}{B(\frac{d+1}{2}, \frac{1}{2})} \int_0^1 
 \bigg[f^\prime \Big(\frac{2r\sqrt{y}}{t}\Big) \frac{2r\sqrt{y}}{t^4} +
 f\Big(\frac{2r\sqrt{y}}{t}\Big)\frac{2}{t^3}
 \bigg]
 (1-y)^{\frac{1}{2}(d-1)} \rd y
\\
&=& \frac{1}{B(\frac{d+1}{2}, \frac{1}{2}) r t}
 \int_0^{2r/t} \big[x^2 f^\prime (x) + 2xf(x)\big]\Big(1- \frac{t^2x^2}{4r^2}\Big)^{\frac{1}{2}(d-1)}
 \rd x.
 \end{eqnarray*}
 From the first equality we have 
 \begin{align*}
  |H^{\prime\prime}(t)|\leq \frac{4r^2L}{t^4} \frac{B(\frac{d+1}{2},\frac{1}{2})}{B(\frac{d+1}{2},\frac{3}{2})}
  + \frac{4rL}{t^3} \frac{B(\frac{d+1}{2},1)}{B(\frac{d+1}{2},\frac{1}{2})}\leq 
  4Lr \bigg(\frac{(d+2)r}{t^4} + \frac{1}{t^3}\bigg),
 \end{align*}
 provided that $\sup_x f(x) \vee \sup_x |f^\prime(x)|\leq L$.
\par\medskip  
(iii). By the third inequality in 
(\ref{eq:L-bound}), 
$\lim_{x\to\infty} x^2 f(x)=0$ and  $\int_0^\infty [2xf(x)+ x^2f^\prime(x)]\rd x=0$; therefore  
\begin{align}
 |H^{\prime\prime}(t)| &\leq  \frac{1}{B(\frac{d+1}{2}, \frac{1}{2}) r t}\bigg\{
 \bigg|\int_0^{2r/t} \big[x^2 f^\prime (x) + 2xf(x)\big]\,
 \Big[\Big(1- \frac{t^2x^2}{4r^2}\Big)^{\frac{1}{2}(d-1)} -1\Big]
 \rd x\bigg|
\nonumber
 \\
 &\hspace{80mm} + \Big|\int_{2r/t}^\infty \big[2xf(x)+x^2f^\prime(x)\big] \rd x\Big|\bigg\}
\nonumber
 \\
&= \frac{1}{B(\frac{d+1}{2}, \frac{1}{2}) r t}\bigg\{\bigg|
 \int_0^{2r/t} \big[x^2 f^\prime (x) + 2xf(x)\big]\,
 \Big[\Big(1- \frac{t^2x^2}{4r^2}\Big)^{\frac{1}{2}(d-1)} -1\Big]
 \rd x \bigg| + \frac{4r^2}{t^2} f\Big(\frac{2r}{t}\Big)\bigg\}.
\label{eq:Hpp}
 \end{align}
If $d=1$ then 
\begin{align*}
 |H^{\prime\prime}(t)| &\leq  \frac{1}{B(\frac{d+1}{2}, \frac{1}{2}) r t}\frac{4r^2}{t^2} 
 f\Big(\frac{2r}{t}\Big) \leq \frac{L}{2r^2 B(\frac{d+1}{2},\frac{1}{2})},
\end{align*}
where we have used the first inequality in (\ref{eq:L-bound}). If $d=2$ then  by the elementary inequality
$1- (1-a)^{1/2}\leq a$ for $a\in [0,1]$ and by (\ref{eq:L-bound})
 \begin{align*}
  |H^{\prime\prime}(t)| &\leq \frac{1}{B(\frac{d+1}{2}, \frac{1}{2}) r t}
  \bigg\{\int_0^{2r/t} \big[x^2|f^\prime(x)| + 2xf(x)\big]\frac{t^2x^2}{4r^2}\rd x + \frac{4r^2}{t^2} 
 f\Big(\frac{2r}{t}\Big)\bigg\}
 \\
 &\leq \frac{1}{2B(\frac{d+1}{2}, \frac{1}{2}) r^2} \int_0^{2r/t} \big[x^3 |f^\prime(x)| + 2x^2 f(x)\big]\rd x
 + \frac{L}{2r^2 B(\frac{d+1}{2},\frac{1}{2})} \leq \frac{2L}{r^2 B(\frac{d+1}{2},\frac{1}{2})}.
 \end{align*}
Finally, if $d\geq 3$ then expanding in Taylor's series in (\ref{eq:Hpp}) 
we obtain 
\begin{align*}
 |H^{\prime\prime}(t)| & \leq  \frac{1}{B(\frac{d+1}{2}, \frac{1}{2}) r t}
 \bigg\{\int_0^{2r/t} \big[ x^2 |f^\prime(x)|+ 2xf(x)\big] \frac{(d-1)x^2t^2}{2r^2} \rd x
 + \frac{4r^2}{t^2} f\Big(\frac{2r}{t}\Big)\bigg\}
 \\
 &\leq \frac{d-1}{B(\frac{d+1}{2}, \frac{1}{2}) r^2}\int_0^{2r/t} \big[ x^3 |f^\prime(x)|+ 2x^2f(x)\big]\rd x
 + \frac{L}{2r^2 B(\frac{d+1}{2},\frac{1}{2})} < \frac{3dL}{r^2 B(\frac{d+1}{2},\frac{1}{2})},
\end{align*}
where we took into account (\ref{eq:L-bound}).
This completes the proof.
\epr 
\par 
Now we proceed to the proof of  the theorem.
\paragraph{Proof of Theorem~\ref{th:expected-value}.}
In the subsequent proof $c_1,c_2,\ldots$ stand for positive constants that may depend 
on $d$  and characteristics of kernel $K$ only. These constants may be different
on different occasions. 
\par 
We have the following bias--variance decomposition of the mean squared error of $\hat{\psi}_{h}$:
\begin{align}
 \frac{\rho}{r B(\frac{d+1}{2}, \frac{1}{2})} \Big[\bE |\hat{\mu}_{h} - \mu|^2\Big]^{1/2}
 \leq 
 \bigg\{\bE \Big|\frac{1}{h^2}\int_{0}^{h} K\Big(\frac{t}{h}\Big)\big[\hat{R}(t) - R(t)\big]\rd t\Big|^2\bigg\}^{1/2} 
 \nonumber
 \\
 +
 \bigg|\frac{1}{h^2}\int_{0}^{h} K\Big(\frac{t}{h}\Big)R(t)\rd t - R^\prime (0+)\bigg|.
\label{eq:bias-variance}
 \end{align}
Our current goal is to bound from above the two terms on the right hand side of the above formula.
\par 
The bound on the  bias is immediate.
Expanding $R(t)$ is Taylor's series and using Lemma~\ref{lem:H-diff} we obtain
\begin{align*}
 \bigg|\frac{1}{h^2}\int_{0}^{h} K\Big(\frac{t}{h}\Big)R(t)\rd t - R^\prime (0+)\bigg| \leq 
 \tfrac{1}{2} h C_K \sup_{0< y\leq h} |R^{\prime\prime}(y)| \leq 
 \frac{3d C_K\rho L h}{2r^2 B(\frac{d+1}{2},\frac{1}{2})} =  c_1\rho Lh r^{-2}.
\end{align*}
We continue with bounding the variance.
\par\medskip 
{\em Bound on the variance.}
 We  have   
\begin{align}\label{eq:var}
V :=
& \bE \bigg|\frac{1}{h^2}\int_{0}^{h} K\Big(\frac{t}{h}\Big)\big[\hat{R}(t) - R(t)\big]\rd t\bigg|^2
 \nonumber
 \\
&\;\;\; = \frac{1}{h^4}\int_0^h\int_0^h K\Big(\frac{t}{h}\Big)K\Big(\frac{s}{h}\Big)
 \bE \big[\hat{R}(t) - R(t)\big]\big[\hat{R}(s) - R(s)\big] \rd t\,\rd s.
\end{align}
In our derivation of the upper bound on the variance we substitute formula (\ref{eq:R-Rhat}) given
in Theorem~\ref{lem:R-Rhat} in (\ref{eq:var}) and  bound the resulting terms. 
The proof proceeds in the following steps.
\par
1$^0$.
Denote for brevity 
\begin{align*}
 A_1(\tau_1-\tau_2, s, t) &:=  H(\tau_1-\tau_2)H(\tau_1-\tau_2+t-s)+
 H(\tau_1-\tau_2-s)H(\tau_1-\tau_2+t)
 \\
 A_2(\tau_1-\tau_2, s,t) &:= U(\tau_1, \tau_1+t, \tau_2, \tau_2+s)= U(\tau_1-\tau_2, 
 \tau_1-\tau_2+t, 0, s),
\end{align*}
and remind that $A_1$ and $A_2$ are non--negative functions.
Then  using Theorem~\ref{lem:R-Rhat} we can write 
\begin{align*}
 (T-t)(T-s)  \bE \big[ &\hat{R}(t)  - R(t)\big]\big[\hat{R}(s) - R(s)\big] 
\nonumber
 \\
 &= \int_0^{T-s}\int_0^{T-t} 
 \Big\{\rho^2 A_1(\tau_1-\tau_2, s, t) 
 + \rho  A_2(\tau_1-\tau_2, s, t)\Big\}
 \rd \tau_1\rd \tau_2
 \\
  &=  \int_0^{T}\int_0^{T} \Big\{\rho^2 A_1(\tau_1-\tau_2, s, t) 
 + \rho A_2(\tau_1-\tau_2, s, t) \Big\}
 \rd \tau_1\rd \tau_2 + E(t, s),
 \nonumber
 \end{align*}
where 
\begin{align*}
 E(t, s) :=  \int_0^{T-t}\int_{T-s}^{T} \Big\{\rho^2 A_1(\tau_1-\tau_2, t, s) 
 + \rho A_2(\tau_1-\tau_2, s, t)\Big\}\rd 
 \tau_1
 \rd \tau_2
 \nonumber
\\
 \;\;\;+ \int_{T-t}^{T} \int_0^{T} \Big\{\rho^2 A_1(\tau_1-\tau_2, t, s)\rd 
 \tau_1 + \rho A_2(\tau_1-\tau_2, s, t)\Big\}
\rd\tau_1  \rd \tau_2.
\end{align*}
Note that 
for $i=1,2$ and for $t, s\in [0,h]$ one has 
\begin{align*}
 \int_0^{T-t} \int_{T-s}^{T} 
 A_i(\tau_1-\tau_2, t, s) \rd \tau_1\rd \tau_2 \leq 2 \int_0^{T-t} 
 \int_{T-s}^{T}
 H(\tau_1-\tau_2)\rd \tau_1\rd \tau_2 
 \leq 4h \int_0^T H(y)\rd y,
\end{align*}
and similarly, 
\begin{align*}
 \int_{T-t}^{T} \int_{0}^{T} A_i(\tau_1-\tau_2, t, s) \rd \tau_1\rd \tau_2 \leq 2 
 \int_{T-t}^{T} \int_{0}^{T}
 H(\tau_1-\tau_2)\rd \tau_1\rd \tau_2 \leq 4h \int_0^T H(y)\rd y.
\end{align*}
Therefore $|E(t,s)|\leq 8h(\rho^2+\rho)\int_0^T H(y)\rd y$, and 
\begin{align*}
& \bigg|\frac{1}{h^4}\int_0^h\int_0^h K\Big(\frac{t}{h}\Big)K\Big(\frac{s}{h}\Big) 
 \frac{E(t,s)}{(T-s)(T-t)}\rd t\,\rd s\bigg| 
 \nonumber
 \\*[2mm]
 & \leq 
 \frac{c_2 (\rho^2+\rho)}{h^3 T^2}\int_0^h\int_0^h \bigg|K\Big(\frac{t}{h}\Big)K\Big(\frac{s}{h}\Big)\bigg|
 \rd t\,\rd s \int_0^T H(y)\rd y
\leq \frac{c_3(\rho^2+\rho)}{h T^2} \int_0^T H(y)\rd y.
 \end{align*}
Thus,
\begin{align}\label{eq:var-psi}
 V \leq  \tilde{V} + \frac{c_3(\rho^2+\rho)}{h T^2} \int_0^T H(y)\rd y,
\end{align}
where we denoted  
\begin{align*}
\tilde{V} & 
:= \frac{1}{h^4} \int_0^h\int_0^h \bigg[ K\Big(\frac{s}{h}\Big)K\Big(\frac{t}{h}\Big)
\frac{1}{(T-s)(T-t)}
\nonumber 
\\
  &\;\;\;\;\;\;\;\;\;\;\times 
  \int_0^{T}\int_0^{T} 
 \Big\{\rho^2A_1(\tau_1-\tau_2, t, s)+ \rho A_2(\tau_1-\tau_2, t,s)
 \Big\}\rd \tau_1 \rd \tau_2 \bigg]\rd t\,\rd s.
 \end{align*}
Denote 
\begin{align}
\bar{V} := \frac{1}{h^4} \int_0^h\int_0^h K\Big(\frac{s}{h}\Big)K\Big(\frac{t}{h}\Big)
\frac{1}{T^2}
\int_0^{T}\int_0^{T} 
 \Big\{\rho^2A_1(\tau_1-\tau_2, t, s)+ \rho A_2(\tau_1-\tau_2, t,s)
 \Big\}\rd \tau_1 \rd \tau_2\rd t\,\rd s.
\label{eq:Vbar}
 \end{align} 
Because 
 for $t, s\in [0,h]$ 
\begin{align*}
 \bigg|\frac{1}{(T-s)(T-t)} -\frac{1}{T^2}\bigg| 
 \int_0^{T}\int_0^{T} 
 \Big\{\rho^2A_1(\tau_1-\tau_2, t, s)+ \rho A_2(\tau_1-\tau_2, t,s)
 \Big\}\rd \tau_1 \rd \tau_2\rd t \rd s
\\
 \;\leq\; \frac{c_4h(\rho^2+\rho)}{(T-h)^2}  \int_{0}^T H(y)\rd y
\end{align*}
we have 
\begin{equation*}
 |\tilde{V} - \bar{V}| \leq \frac{c_5h(\rho^2+\rho)}{(T-h)^2}\int_0^T H(y)\rd y,
\end{equation*}
 Combining this inequality with (\ref{eq:var-psi}) we obtain 
 \begin{equation}\label{eq:var-psi-1}
  V \leq \bar{V}+  c_6(\rho^2+\rho)\bigg(\frac{1}{hT^2}+ 
  \frac{h}{(T-h)^2}\bigg)\int_0^T H(t)\rd t.
 \end{equation}
 Our current goal is to bound $\bar{V}$ from above [see (\ref{eq:Vbar})]. Write 
 \begin{align}
  \bar{V} &=: \rho^2 V_1 + \rho V_2
  \label{eq:V-bar}
  \\
  V_1 & :=\frac{1}{h^4} \int_0^h\int_0^h K\Big(\frac{s}{h}\Big)K\Big(\frac{t}{h}\Big)
\frac{1}{T^2}
\int_0^{T}\int_0^{T}  A_1(\tau_1-\tau_2, t, s) \rd\tau_1\rd\tau_2\rd t\rd s
\nonumber
\\
V_2 & := \frac{1}{h^4} \int_0^h\int_0^h K\Big(\frac{s}{h}\Big)K\Big(\frac{t}{h}\Big)
\frac{1}{T^2}
\int_0^{T}\int_0^{T}  A_2(\tau_1-\tau_2, t, s) \rd\tau_1\rd\tau_2\rd t\rd s.
\nonumber
\end{align}
\par 
2$^0$. 
Consider first $V_1$. 
Because 
\[
 A_1(\tau_1-\tau_2, t, s)=H^2(\tau_1-\tau_2)+ H(\tau_1-\tau_2-s)H(\tau_1-\tau_2+t)
\]
we have $V_1= V_{1,1}+V_{1,2}$ where 
\begin{align*}
V_{1,1} := &  \frac{1}{h^4} \int_0^h\int_0^h K\Big(\frac{s}{h}\Big)K\Big(\frac{t}{h}\Big)
 \frac{1}{T^2}\int_0^{T}\int_0^{T}   
 H(\tau_1-\tau_2-s) H(\tau_1-\tau_2+t)\rd \tau_1\rd\tau_2 \rd t\rd s,
\end{align*} 
and 
\begin{align*}
 V_{1,2}:= \frac{1}{h^4} \int_0^h\int_0^h K\Big(\frac{s}{h}\Big)K\Big(\frac{t}{h}\Big)
 \frac{1}{T^2}\int_0^{T}\int_0^{T}   
 H^2(\tau_1-\tau_2)\rd \tau_1\rd\tau_2 \rd t\rd s.
\end{align*}
By  (\ref{eq:K-Kernel}), $V_{1,2}=0$; therefore we need to bound $V_{1,1}$ only.
To this end, 
let us  introduce the following subsets of $\bR^2$
 \begin{align}
  &\cS_1:=\big\{(\tau_1, \tau_2): 0\leq \tau_1-\tau_2\leq h\big\}\cap [0, T]^2,\;\;\;\;
  \label{eq:setS1}
  \\
  & \cS_2:= \big\{(\tau_1, \tau_2): -h\leq  \tau_1-\tau_2\leq 0\big\}\cap [0, T]^2,\;\;\;
  \cS:=
  [0, T]^2 \setminus (\cS_1\cup \cS_2),
  \label{eq:setS2}
 \end{align}
and divide the double integral over $\tau_1$ and  $\tau_2$ in three integrals  corresponding to  these subsets. 
For the set $\cS_1$  we obtain 
\begin{align*}
V_{1,1}(\cS_1)& :=\frac{1}{h^4} \int_0^h\int_0^h K\Big(\frac{s}{h}\Big)K\Big(\frac{t}{h}\Big)
 \frac{1}{T^2}\iint_{\cS_1} 
 H(\tau_1-\tau_2-s) H(\tau_1-\tau_2+t)\rd \tau_1\rd\tau_2 \rd t\rd s
 \\
& = \frac{1}{h^2}\int_0^h K\Big(\frac{s}{h}\Big)  \frac{1}{T^2}\iint_{\cS_1}
H(\tau_1-\tau_2-s) \Big\{\frac{1}{h^2}\int_0^h K\Big(\frac{t}{h}\Big) H(\tau_1-\tau_2+t)\rd t\Big\}\rd \tau_1\rd\tau_2 \rd s.
\end{align*}
For $(\tau_1, \tau_2)\in \cS_1$ we have $\tau_1-\tau_2+t >0$ for all $t\in [0,h]$ so that 
$H(\tau_1-\tau_2+t)$ is smooth and can be expanded in Taylor's series around $\tau_1-\tau_2$:
\[
 \frac{1}{h^2}\int_0^h K\Big(\frac{t}{h}\Big) H(\tau_1-\tau_2+t)\rd t= H^\prime (\tau_1-\tau_2) + 
 \frac{h}{2!}\int_0^1 y^2 K(y) H^{\prime\prime}(\tau_1-\tau_2 + \vartheta yh)\rd y,\;\;\;
 \vartheta\in [0, 1].
\]
We have
\begin{align*}
 \bigg|\frac{1}{h^2}\int_0^h K\Big(\frac{t}{h}\Big) H(\tau_1-\tau_2+t)\rd t - H^\prime (\tau_1-\tau_2)\bigg|
\leq \frac{h}{2!}\int_0^1 y^2 |K(y)| |H^{\prime\prime}(\tau_1-\tau_2 + \vartheta yh)|\rd y
\leq c_6Lhr^{-2},
\end{align*}
where in the last inequality we have used (\ref{eq:Hprime-prime}) of Lemma~\ref{lem:H-diff}.
Therefore 
\begin{align*}
 |V_{1,1}(\cS_1)| &\leq \frac{1}{h^2}\int_0^h \Big|K\Big(\frac{s}{h}\Big)\Big|
 \frac{1}{T^2}\iint_{\cS_1} H(\tau_1-\tau_2-s) \big\{ |H^\prime (\tau_1-\tau_2)| 
 + c_6Lhr^{-2} \big\}\rd \tau_1\rd \tau_2 \rd s 
 \\
& \leq \frac{c_7}{hT^2}\iint_{\cS_1} |H^\prime (\tau_1-\tau_2)|\rd\tau_1\rd \tau_2  +  
 \frac{c_8L h}{r^2T} \leq \frac{c_9L}{rT}(1+hr^{-1}),
\end{align*}
where we again use Lemma~\ref{lem:H-diff}.
Similarly,
\begin{align*}
V_{1,1}(\cS_2)& :=\frac{1}{h^4} \int_0^h\int_0^h K\Big(\frac{s}{h}\Big)K\Big(\frac{t}{h}\Big)
 \frac{1}{T^2}\iint_{\cS_2} 
 H(\tau_1-\tau_2-s) H(\tau_1-\tau_2+t)\rd \tau_1\rd\tau_2 \rd t\rd s
 \\
 & = \frac{1}{h^2}\int_0^h K\Big(\frac{t}{h}\Big)  \frac{1}{T^2}\iint_{\cS_2}
H(\tau_1-\tau_2+t) \Big\{\frac{1}{h^2}\int_0^h K\Big(\frac{s}{h}\Big) H(\tau_1-\tau_2-s)\rd s\Big\}\rd \tau_1\rd\tau_2 \rd t.
\end{align*}
and by the same reasoning as above we obtain the same bound
\begin{align*}
 |V_{1,1}(\cS_2)| 
 \leq \frac{c_{10}L}{rT} (1+hr^{-1}). 
\end{align*}
\par 
Now we consider 
\begin{align*}
&V_{1,1}(\cS) :=\frac{1}{h^4} \int_0^h\int_0^h K\Big(\frac{s}{h}\Big)K\Big(\frac{t}{h}\Big)
 \frac{1}{T^2}\iint_{\cS} 
 H(\tau_1-\tau_2-s) H(\tau_1-\tau_2+t)\rd \tau_1\rd\tau_2 \rd t\rd s
 \\
 &\; = \frac{1}{T^2}\iint_{\cS} \Big\{\frac{1}{h^2} \int_0^h K\Big(\frac{s}{h}\Big) H(\tau_1-\tau_2-s)\rd s\Big\} 
  \Big\{\frac{1}{h^2} \int_0^h K\Big(\frac{t}{h}\Big) H(\tau_1-\tau_2+t)\rd t\Big\}
  \rd\tau_1 \rd\tau_2.
\end{align*}
On the set $\cS$ functions $H(\tau_1-\tau_2 -\cdot)$ and $H(\tau_1-\tau_2 +\cdot)$
can be expanded to Taylor's series up to the second order; therefore for 
$(\tau_1, \tau_2) \in \cS$ one has 
\begin{align*}
\bigg| \frac{1}{h^2} \int_0^h K\Big(\frac{s}{h}\Big) H(\tau_1-\tau_2-s)\rd s - H^\prime(\tau_1-\tau_2)\bigg| 
= \tfrac{1}{2} h\bigg|\int_0^1 y^2 K(y) H^{\prime\prime}(\tau_1-\tau_2-\vartheta_1yh) \rd y\bigg|
\\
\leq \tfrac{1}{2}C_K h
\max_{y\in [0,h]} \big|H^{\prime\prime}(\tau_1-\tau_2-y)\big|
\leq  \frac{3dC_KLh}{2r^2B(\frac{d+1}{2}, \frac{1}{2})}= c_{11}Lh r^{-2}
 \\
 \bigg|\frac{1}{h^2} \int_0^h K\Big(\frac{t}{h}\Big) H(\tau_1-\tau_2+t)\rd s - 
 H^\prime(\tau_1-\tau_2)\bigg| = \tfrac{1}{2} h\bigg|\int_0^1 y^2 K(y) 
 H^{\prime\prime}(\tau_1-\tau_2+\vartheta_2yh) \rd y\bigg|
 \\
 \leq 
 \tfrac{1}{2} C_Kh \max_{y\in [0,h]} \big|H^{\prime\prime}(\tau_1-\tau_2+y)\big|
 \leq \frac{3dC_KLh}{2r^2B(\frac{d+1}{2}, \frac{1}{2})}=c_{11}Lh r^{-2},
\end{align*}
where $\vartheta_1, \vartheta_2\in [0,1]$ and we have used Lemma~\ref{lem:H-diff}.
This yields
\begin{align*}
|V_{1,1}(\cS)|& \leq \frac{1}{T^2}\iint_{\cS} \big[|H^{\prime}(\tau_1-\tau_2)|+c_{11}Lhr^{-2}\big]^2
\rd \tau_1\rd \tau_2
\\
&\leq \frac{2}{T^2} \iint_{\cS} |H^\prime(\tau_1-\tau_2)|^2 \rd\tau_1\rd\tau_2
\;+ c_{12} h^2 L^2 r^{-4}
 \leq  c_{13}\Big\{\frac{L^2}{rT}  +  L^2h^{2}r^{-4}\Big\},
\end{align*}
where in the last line we have used (\ref{eq:Hprime-bound}) and (\ref{eq:Hprime-int}).
Combining the obtained inequalities for 
$|V_{1,1}(\cS_1)|, |V_{1,1}(\cS_2)|$ and $|V_{1,1}(\cS)|$  we obtain 
\begin{align}\label{eq:V1-final}
 |V_{1}| \leq c_{13}\Big\{\frac{L}{rT}\big[1+ L +  hr^{-1}\big] +  L^2 h^2 r^{-4}\Big\}.
\end{align}
\par 
4$^0$.
Now we consider  
\begin{align*}
 V_{2}:= \frac{1}{h^4}\int_0^h \int_0^h K\Big(\frac{t}{h}\Big)K\Big(\frac{s}{h}\Big)
 \frac{1}{T^2}\int_0^T\int_0^T A_2(\tau_1-\tau_2, s, t) \rd \tau_1\rd \tau_2\rd s \rd t.
\end{align*}
Write for brevity 
$\tau_{\max} :=\max\{0, \tau_1-\tau_2, \tau_1-\tau_2+s, t\}$ and 
$\tau_{\min} := \min\{0, \tau_1-\tau_2, \tau_1-\tau_2+s, t\}$, and recall that 
\begin{align*}
 A_2(\tau_1-\tau_2, s, t)= \frac{1}{{\rm vol}(B)} \bE 
 {\rm vol}\{B\cap B(v(\tau_1-\tau_2))\cap B(v(\tau_1-\tau_2+t))\cap B(vs)\}
\\ =
 \frac{1}{{\rm vol}(B)} \bE {\rm vol}\{B(v(\tau_{\max}-\tau_{\min}))\}.
\end{align*}
Note that for $t, s\in [0,h]$ we have 
$\tau_{\max}=\max\{\tau_1-\tau_2+s, t\}$ and 
$\tau_{\min} = \min\{0, \tau_1-\tau_2\}$ so that 
\[
 A_2(\tau_1-\tau_2, s, t)= H (\tau_{\max}-\tau_{\min}).
\]
Now we proceed 
by partitioning $[0,T]^2$ in subsets $\cS_1$, $\cS_2$ and $\cS$ as defined in 
(\ref{eq:setS1})--(\ref{eq:setS2}).  
On the set $\cS_1$ we have $\tau_{\max}= (\tau_1-\tau_2+s)\vee t$ and $\tau_{\min}=0$ 
so that 
\begin{align}
 V_2(\cS_1) & := \frac{1}{h^2}\int_0^h K\Big(\frac{s}{h}\Big)
 \frac{1}{T^2}\iint_{\cS_1} \Big\{
 \frac{1}{h^2} \int_0^h K\Big(\frac{t}{h}\Big)
 H((\tau_1-\tau_2+s)\vee t) \rd t\Big\} \rd \tau_1\rd \tau_2 \rd s,
 \label{eq:V2S1}
 \end{align}
and 
\begin{align*}
 &\bigg|\frac{1}{h^2} \int_0^h K\Big(\frac{t}{h}\Big)
 H((\tau_1-\tau_2+s)\vee t) \rd t \bigg|
 \\
 & = \bigg|\frac{1}{h^2}
 \int_0^h K\Big(\frac{t}{h}\Big)
 H(t) {\bf 1}\{t\geq \tau_1-\tau_2+s\} \rd t  
 + \frac{1}{h^2}
 \int_0^h K\Big(\frac{t}{h}\Big)
 H(\tau_1-\tau_2+s) {\bf 1}\{t < \tau_1-\tau_2+s\} \rd t  \bigg|
 \\
 &\leq \bigg|\frac{1}{h^2}
 \int_0^h K\Big(\frac{t}{h}\Big)
 H(t)  \rd t\bigg|  
 + 
 \bigg|
 \frac{1}{h^2}
 \int_0^h K\Big(\frac{t}{h}\Big)
 \big[H(\tau_1-\tau_2+s) - H(t)\big] {\bf 1}\{t < \tau_1-\tau_2+s\} \rd t\bigg|
 \\
 & \leq 
 |H^\prime(0+)| +  c_{1}Lhr^{-2} +  c_{2}Lr^{-1} \leq c_{3} Lr^{-1}(1+hr^{-1}),
\end{align*}
where in the last line we have used bounds on $|H^\prime(t)|$ established in Lemma~\ref{lem:H-diff},
and the fact that $|\tau_1-\tau_2+s-t|\leq 2h$ on the set $\cS_1$. 
Substituting this bound in (\ref{eq:V2S1}) we obtain
\[
 |V_2(\cS_1)| \leq \frac{1}{h^2}\int_0^h \Big|K\Big(\frac{s}{h}\Big)\Big|\rd s
 \frac{1}{T^2} \iint_{\cS_1} c_{16}Lr^{-1}(1+hr^{-1})\rd \tau_1\rd \tau_2 \leq
 \frac{c_{4}L}{rT} (1+hr^{-1}).
\]
On the set $\cS_2$, $\tau_{\max}=(\tau_1-\tau_2+s)\vee t$
and $\tau_{\min}=\tau_1-\tau_2$; therefore 
\begin{align}
 V_2(\cS_2) & := \frac{1}{h^2}\int_0^h K\Big(\frac{t}{h}\Big)
 \frac{1}{T^2}\iint_{\cS_2} \Big\{
 \frac{1}{h^2} \int_0^h K\Big(\frac{s}{h}\Big)
 H((\tau_1-\tau_2+s)\vee t - (\tau_1-\tau_2)) \rd s\Big\} \rd \tau_1\rd \tau_2 \rd t,
 \label{eq:V2S2}
 \end{align}
and similarly to bounding $V_2(\cS_1)$ we have 
\begin{align*}
 &\bigg|\frac{1}{h^2} \int_0^h K\Big(\frac{s}{h}\Big)
 H\big((\tau_1-\tau_2+s)\vee t - (\tau_1-\tau_2)\big) \rd s \bigg|
 \\
 & = \bigg|\frac{1}{h^2}
 \int_0^h K\Big(\frac{s}{h}\Big)
 H(t-(\tau_1-\tau_2)) {\bf 1}\{t\geq \tau_1-\tau_2+s\} \rd s  
 + \frac{1}{h^2}
 \int_0^h K\Big(\frac{s}{h}\Big)
 H(s) {\bf 1}\{t < \tau_1-\tau_2+s\} \rd t  \bigg|
 \\
 &\leq \bigg|\frac{1}{h^2}
 \int_0^h K\Big(\frac{s}{h}\Big)
 H(s)  \rd s\bigg|  
 + 
 \bigg|
 \frac{1}{h^2}
 \int_0^h K\Big(\frac{s}{h}\Big)
 \big[H(t-(\tau_1-\tau_2)) - H(s)\big] {\bf 1}\{t < \tau_1-\tau_2+s\} \rd s\bigg|
 \\
 & \leq 
 |H^\prime(0+)| +  c_{1}Lhr^{-2} +  c_{2}Lr^{-1} \leq c_{3} Lr^{-1}(1+hr^{-1}),
\end{align*}
so that we also have 
\[
 |V_2(\cS_2)|\leq \frac{c_{4}L}{rT}(1+hr^{-1}).
\]
On the set $\cS\setminus\cS_1$ we have $\tau_{\max}=\tau_1-\tau_2+s$ and 
$\tau_{\min}=0$ so that 
\begin{align*}
 V_2(\cS\setminus \cS_1) = \frac{1}{h^4}\int_0^h \int_0^h 
 K\Big(\frac{s}{h}\Big) K\Big(\frac{t}{h}\Big)
 \iint_{\cS\setminus \cS_1} H(\tau_1-\tau_2+s) \rd \tau_1\rd \tau_2\rd t\rd s =0.
\end{align*}
On the set $\cS\setminus \cS_2$ we have  $\tau_{\max}=t$, $\tau_{\min}=\tau_1-\tau_2$
so that 
\begin{align*}
 V_2(\cS\setminus \cS_2) = \frac{1}{h^4}\int_0^h \int_0^h 
 K\Big(\frac{s}{h}\Big) K\Big(\frac{t}{h}\Big)
 \iint_{\cS\setminus \cS_2} H\big(t-(\tau_1-\tau_2)\big) \rd \tau_1\rd \tau_2\rd t\rd s =0.
\end{align*}
Finally, combining all these bounds we obtain
\begin{align}\label{eq:V2-final}
 |V_2| \leq \frac{c_{5}L}{rT} (1+ hr^{-1}). 
\end{align}
\par 
4$^0$. Now we are in a position to complete the theorem proof. 
Combining (\ref{eq:V1-final}), (\ref{eq:V2-final}), (\ref{eq:V-bar}) with (\ref{eq:var-psi-1})
we obtain
\begin{align*}
 V \leq \frac{c_1(\rho^2+\rho)L}{rT}
 \Big[1+ L  + hr^{-1}\Big]+ c_2\rho^2 L^2r^{-4} h^2 
 +c_3(\rho^2+\rho)\bigg(\frac{1}{hT^2}+ 
  \frac{h}{(T-h)^2}\bigg)\int_0^T H(t)\rd t.
\end{align*}
Combining this inequality with (\ref{eq:bias-variance}) we obtain 
\begin{align*}
 \bE |\hat{\mu}_h-\mu|^2 \leq & c_4\Big(1+\frac{1}{\rho}\Big)\frac{Lr}{T} \Big[1+L+ hr^{-1}\Big]+ 
 c_5 L^2r^{-2} h^2
 \\
&\;\;\;+ c_6 r^2\Big(1+\frac{1}{\rho}\Big)\bigg(\frac{1}{hT^2}+ 
  \frac{h}{(T-h)^2}\bigg)\int_0^T H(t)\rd t.
\end{align*}
The result of the theorem follows from the above inequality by selecting $h_*$
as stated in the premise of the theorem. Under this choice  as $T\to\infty$ the first term on the right hand side dominates 
the other two terms and leads to the announced result.
We also took into account that 
for bounded densities 
$\int_0^T H(t)\rd t \leq  O(\ln T)$ as $T\to\infty$.
\epr
\subsection{Proof of Lemma~\ref{lem:properties-varphi}}
By condition~(K),  
$K(x/h)$ is supported on $[0,h]$ and $K(\ln (x/x_0)/h)$ is supported on $[x_0, x_0e^h]$.
Therefore  we have for $h\leq t\leq x_0$  
\begin{align*}
 \varphi_{x_0,h}(t)= \frac{1}{h} \int_0^t K\Big(\frac{t}{h}\Big)\rd x =  \frac{1}{h} \int_0^1 K\Big(\frac{t}{h}\Big)\rd x=1,\;\;\;\hbox{for}\;\;h\leq t\leq x_0,
\end{align*}
and for $t\geq x_0e^h$ 
\begin{align*}
 \varphi_{x_0,h}(t)=1  - \int_0^t \frac{1}{xh} K\Big(\frac{\ln (x/x_0)}{h}\Big)\rd x=
 1- \int_0^{\ln (t/x_0)/h}  K(u)\rd u =0.
\end{align*}
This proves the first statement. 
\par 
Since $K$ is compactly supported, $\widehat{K}$ is an entire function. 
Also, $\widetilde{K}$ is entire because $K$ is infinitely differentiable at zero and compactly 
supported. 
Now we compute the Mellin transform of $\varphi_{x_0,h}$:
\begin{align*}
 \widetilde{\varphi}_{x_0,h}(z) &= \int_0^{x_0e^h} t^{z-1} \varphi_{x_0}(t)\rd t = \int_0^{x_0e^h} 
 \int_0^{x_0e^h} t^{z-1}{\bf 1}(x\leq t) \Big[\frac{1}{h} K\Big(\frac{x}{h}\Big) - \frac{1}{xh}K\Big(\frac{\ln(x/x_0)}{h}\Big) \Big] \rd x \rd t
 \\
 &=\int_0^{x_0e^h} \frac{1}{z}\big[x_0^z e^{zh}- x^z\big] \Big[\frac{1}{h} K\Big(\frac{x}{h}\Big) - \frac{1}{xh}K\Big(\frac{\ln(x/x_0)}{h}\Big) \Big] \rd x
 \\
 &= - \frac{1}{z}\int_0^{x_0e^h} x^z \Big[\frac{1}{h} K\Big(\frac{x}{h}\Big) - \frac{1}{xh}K\Big(\frac{\ln(x/x_0)}{h}\Big) \Big] \rd x= - \frac{1}{z}\big[h^z\widetilde{K}(z+1) - x_0^z 
 \widehat{K}(-zh)\big],
\end{align*}
where in the last equality we took into account that $x_0e^h>h$. The lemma is proved.
\epr 
\subsection{Proof of Lemma~\ref{lem:psi}}
We begin with the proof of (\ref{eq:abs-conv}). This condition ensures that the integral in the definition 
of function $\psi_{x_0, h}$ is  absolutely convergent so that $\psi_{x_0,h}$ is well defined.
\par 
1$^0$. We note that 
$|\psi_{x_0,h}(t)| \leq I_1(t)+ I_2(t)$,
where 
\begin{align*}
  &I_1(t) = \frac{t^{-s}x_0^{1-s}}{2\pi}
 \int_{-\infty}^\infty \bigg| \frac{\widehat{K}((s-1+i\omega)h)}{\tilde{w}(1-s-i\omega)}\bigg|
 [(1-s)^2+\omega^2]^{-1/2}
 \rd\omega,\;\;\;s<1,
 \\*[2mm]
 & I_2(t) = \frac{t^{-s} h^{1-s}}{2\pi}
 \int_{-\infty}^\infty \bigg| \frac{\widetilde{K}(2-s-i\omega)}{\widetilde{w}(1-s-i\omega)}\bigg|
 [(1-s)^2+\omega^2]^{-1/2}
 \rd\omega,\;\;\;s<1,
\end{align*}
It suffices to show that the integrals appearing in definitions of $I_1(t)$ and $I_2(t)$ are finite for $s<1$.
\par
First consider the integral in the definition of $I_1(t)$. 
Since 
\begin{equation}\label{eq:widehatKK}
\widehat{K}(\sigma+i\omega) = \int_0^1 e^{(\sigma+i\omega)t} K(t)\rd t,
\end{equation}
in view of condition~(K)
we have 
for any positive integer $\ell$
\begin{equation}\label{eq:hatKbound}
 |\widehat{K}(\sigma+i\omega)| \leq \min\Big\{c_1 e^{|\sigma|},\, c_2(\ell) 
 e^{|\sigma|}(\sigma^2+\omega^2)^{-\ell/2}\Big\},\;\;\;\forall \sigma, \omega,
\end{equation}
where the second inequality follows from the repeated integration by parts in 
(\ref{eq:widehatKK}).
Therefore using Lemma~\ref{lem:Mellin-w} with $\sigma=1-s$ we obtain for any $-d\leq s<1$  
\begin{align*}
 &\int_{-\infty}^\infty \bigg| \frac{\widehat{K}((s-1+i\omega)h)}{\widetilde{w}(1-s-i\omega)}\bigg|
 [(1-s)^2+\omega^2]^{-1/2}
 \rd\omega
 \\
 & \;\;\;\leq \frac{c_3}{(2r)^{1-s}}\bigg\{ \int_{|\omega|\leq 2}   |\widehat{K}((s-1)h+i\omega h)|\rd \omega  
 + \int_{|\omega|\geq 2} |\widehat{K}((s-1)h+i\omega h)| |\omega|^{(d+1)/2}\rd \omega\bigg\}
 \\
  &\;\;\; \leq \frac{c_3}{(2r)^{1-s}}\bigg\{
  4e^{|s-1|h} + h^{-(d+3)/2}\int_{|\xi|\geq 2h} 
  |\widehat{K}((s-1)h+i\xi)|\, |\xi|^{(d+1)/2}\rd \xi\bigg\}
  \\
  &\;\;\;\leq \frac{c_4 e^{|s-1|h}}{(2r)^{1-s}}
  \bigg\{ 
  1+ 
  h^{-(d+3)/2}\bigg(1+\int_{|\xi|\geq 1} 
  |\xi|^{(d+1)/2-\ell} \rd \xi \bigg)\bigg\}
  \leq \frac{c_5 e^{|s-1|h}}{(2r)^{1-s}}h^{-(d+3)/2},
\end{align*}
where in order to get the penultimate inequality we split the integral into
the sets $2h\leq |\xi|\leq 1$ and $|\xi|\geq 1$, and on the first set we use the first 
inequality in (\ref{eq:hatKbound}),  while  on the second set, the second inequality 
 in (\ref{eq:hatKbound})
is used 
with $\ell>(d+3)/2$. 
\par 
Now consider the integral appearing in the definition of $I_2(t)$. Similarly to 
(\ref{eq:hatKbound}) it follows from  
$\widetilde{K}(\sigma+i\omega)= \int_0^1 t^{\sigma-1+i\omega} K(t)\rd t$ that for $\sigma> 0$ 
\begin{equation}\label{eq:tildeKbound}
 |\widetilde{K}(\sigma+i\omega)| \leq c_6(\ell) \min\{1, (\sigma^2+\omega^2)^{-\ell/2}\}. 
\end{equation}
Therefore by Lemma~\ref{lem:Mellin-w} with $\sigma=1-s$ and (\ref{eq:tildeKbound}) 
we have for $-d\leq s<1$
\begin{align*}
 &\int_{-\infty}^\infty 
 \bigg|\frac{\widetilde{K}(2-s-i\omega)}{\widetilde{w}(1-s-i\omega)}\bigg| [(1-s)^2+\omega^2]^{-1/2}\rd \omega
\\
 &\;\;\leq \frac{c_7}{(2r)^{1-s}} \bigg\{\int_{|\omega|\leq 2}  |\widetilde{K}(2-s-i\omega)| \rd \omega +
 \int_{|\omega|\geq 2} |\widetilde{K}(2-s-i\omega)| |\omega|^{(d+1)/2}\rd\omega \bigg\} \leq 
 \frac{c_8}{(2r)^{1-s}}.
\end{align*}
Combining these inequalities with definitions of $I_1(t)$ and $I_2(t)$
we finally obtain for any $-d\leq s<1$ 
\begin{equation}\label{eq:psi-bound}
 |\psi_{x_0,h}(t)|\leq  c_9(2r)^{s-1} t^{-s} 
 \Big(h^{1-s}+ x_0^{1-s}e^{|s-1|h} h^{-(d+3)/2}\Big),\;\;\forall t\geq 0.
\end{equation}
\par 
2$^0$.  Now we prove (\ref{eq:lin-strategy}). In view of (\ref{eq:w})
\begin{align*}
 \int_0^\infty \psi_{x_0,h}(t) H(t)\rd t = \int_0^\infty \psi_{x_0,h}(t) \int_0^\infty w(tx) \rd F(x)\rd t=
 \int_0^\infty \bigg[\int_0^\infty \psi_{x_0,h}(t)w(tx) \rd t\bigg] \rd F(x). 
\end{align*}
Therefore in order to prove (\ref{eq:lin-strategy}) it suffices to show that 
\begin{equation}\label{eq:lin-strategy-1}
 \int_0^\infty \psi_{x_0,h}(t)w(tx)\rd t = \varphi_{x_0,h}(x),\;\;\forall x\geq 0.
\end{equation}
Talking the Mellin transform of the left hand side we obtain 
\begin{align*}
 \int_0^\infty x^{z-1}  \int_0^\infty \psi_{x_0,h}(t)w(tx)\rd t \rd x= 
 \tilde{w}(z) \int_0^\infty t^{-z} \psi_{x_0,h}(t)\rd t = 
 \widetilde{w}(z) \widetilde{\psi}_{x_0,h}(1-z). 
\end{align*}
Note that the convergence region of $\widetilde{w}(z)$ is $\{z: {\rm Re}(z)>0\}$. By definition of 
$\psi_{x_0,h}$ in (\ref{eq:psi}) the Mellin transform of $\psi_{x_0,h}$ is 
$\widetilde{\psi}_{x_0,h}(z)= \widetilde{\varphi}_{x_0}(1-z)/\tilde{w}(1-z)$, 
and its convergence region is
$\{z: {\rm Re}(z)<1\}$. Therefore the function 
$\widetilde{w}(z)\widetilde{\psi}_{x_0,h}(1-z)$ is analytic in $\{z: {\rm Re}(z)>0\}$, 
and in this region 
\[
 \widetilde{w}(z)\widetilde{\psi}_{x_0,h}(1-z) = \widetilde{\varphi}_{x_0,h}(z),\;\;\;\forall z: {\rm Re}(z)>0.
\]
Therefore, by uniqueness of the Mellin transform, (\ref{eq:lin-strategy-1}) is fulfilled for all $x\geq 0$ which complete the proof of the lemma.
\epr  

\subsection{Proof of Theorem~\ref{th:distribution}}\label{sec:proof-of-distr}
First we state a result on 
the rate of decay of the Mellin transform of $w$  on vertical lines 
in the convergence region.
\begin{lemma}\label{lem:Mellin-w}
 The Mellin transform of function $w$ defined in (\ref{eq:w}) is given by 
 \[
  \widetilde{w}(z)=\frac{(2r)^z}{z} \frac{B(\frac{d+1}{2}, \frac{z+1}{2})}{B(\frac{d+1}{2}, \frac{1}{2})},\;\;\; {\rm Re}(z)>0.
 \]
Moreover, for $\sigma>0$ and $\omega\in \bR$ one has  
\begin{eqnarray*}
 |\widetilde{w}(\sigma+i\omega)| &\geq& \frac{C_1 (2r)^\sigma}{\sqrt{\sigma^2+\omega^2}}
 \;\frac{
 \Gamma(\frac{\sigma+1}{2})}{
 \Gamma(\frac{d+\sigma+2}{2})},\;\;\;\;\;\forall |\omega|\leq 2,
 \\*[2mm]
 |\widetilde{w}(\sigma+i\omega)|  &\geq& 
\frac{C_2 (2r)^\sigma}{\sqrt{\sigma^2+\omega^2}} |\omega|^{-(d+1)/2}, \;\;\;\;\;\forall |\omega|\geq 2,
\end{eqnarray*}
where constants $C_1$ and $C_2$ depend on $d$ only.
 \end{lemma}
\pr 
We have 
\begin{align*}
 \widetilde{w}(z) &= \frac{1}{B(\frac{d+1}{2}, \frac{1}{2})} \int_0^\infty t^{z-1}
 \int_0^1 {\bf 1}\{t\leq 2r\sqrt{1-y}\} y^{\frac{1}{2}(d-1)}(1-y)^{-1/2}\rd y \rd t
 \\
 & = \frac{(2r)^z}{z B(\frac{d+1}{2},\frac{1}{2})}\int_0^1   (1-y)^{z/2} y^{\frac{1}{2}(d-1)}(1-y)^{-1/2}\rd y  = \frac{(2r)^z B(\frac{d+1}{2}, \frac{z+1}{2})}{z B(\frac{d+1}{2}, \frac{1}{2})},
\end{align*}
where the second equality holds only if ${\rm Re}(z)>0$ and the last equality holds for 
${\rm Re}(z)>-1$; thus the equality holds for all $z$ such that ${\rm Re}(z)>0$.
Furthermore, for $z=\sigma+i\omega$, $\sigma>0$, $\omega\in \bR$ we have 
\begin{align*}
 |\widetilde{w}(\sigma+i\omega)| = \frac{ (2r)^\sigma}{\sqrt{\sigma^2+\omega^2}}\;
 \frac{\Gamma(\frac{d}{2}+1)}{\Gamma(\frac{1}{2})}\;\;
 \bigg|\frac{\Gamma(\frac{\sigma+1}{2}+\frac{i\omega}{2})}{\Gamma(\frac{d+\sigma+2}{2}+\frac{i\omega}{2})}\bigg|.
\end{align*}
\par 
We use the following well known properties of the Gamma function.
[see, e.g., \citeasnoun{Carlson} and \citeasnoun{Andrews}]:
\begin{itemize}
\item[(i)] function $\Gamma(z)$ does not have zeros on $\bC$, and it is analytic in 
$\bar{\bC}:=\bC \setminus \{0, -1, -2, \ldots\}$;
 \item[(ii)]
 $|\Gamma(x+iy)|\leq \Gamma(x)$ for all $x+iy\in \bar{\bC}$, and 
 $|\Gamma(x+iy)| \geq \Gamma(x) e^{-\pi |y|/2}$, $\forall x\geq 1/2$, $\forall y\in \bR$; 
\item[(iii)] for  all $x_1\leq x \leq x_2$ and  $|y|\geq 2$ 
there exist constants $c_1\leq c_2$
depending on $x_1$ and $x_2$ such that 
\[
c_1 |y|^{x-1/2} e^{-\pi |y|/2}  \leq |\Gamma(x+iy)| \leq c_2 |y|^{x-1/2} e^{-\pi |y|/2}
\]
\end{itemize}
\par 
By property~(ii)
since $(\sigma+1)/2 >1/2$,
\[
 \bigg|\Gamma\Big(\frac{\sigma+1}{2}+\frac{i\omega}{2}\Big)\bigg| \geq 
 \Gamma\Big(\frac{\sigma+1}{2}\Big) e^{-\pi |\omega|/4},\;\;\;\forall \omega,
\]
which yields 
\[
 |\widetilde{w}(\sigma+i\omega)| \geq \frac{(2r)^\sigma}{\sqrt{\sigma^2+\omega^2}}
 \;\frac{
 \Gamma(\frac{d}{2}+1)\Gamma(\frac{\sigma+1}{2})}{\Gamma(\frac{1}{2})
 \Gamma(\frac{d+\sigma+2}{2})}\; e^{-\pi/2},\;\;\;\forall |\omega|\leq 2.
\]
If $|\omega|\geq 2$ then we use property~(iii):
\[
 |\widetilde{w}(\sigma+i\omega)| \geq \frac{(2r)^\sigma}{\sqrt{\sigma^2+\omega^2}}
\frac{c_1 |\frac{1}{2}\omega|^{(\sigma+1)/2}}{c_2 |\frac{1}{2}\omega|^{(d+\sigma+2)/2}}  \geq 
\frac{c_3 (2r)^\sigma}{\sqrt{\sigma^2+\omega^2}} |\omega|^{-(d+1)/2}, \;\;\;\forall |\omega|\geq 2.
 \]
where $c_3$ depends on $d$ only.
This completes the proof.
\epr
\par 
Now we proceed with the proof  of the theorem.
\paragraph{Proof of Theorem~\ref{th:distribution}.}
Throughout the proof $c_1, c_2,\ldots$ stand for positive constants  that may 
depend on $d$, $\beta$ and $\alpha$ only.
\par 
It follows from (\ref{eq:psi-bound}) and the bound on $H(t)$ in (\ref{eq:H(t)<=}) 
that  
$\int_0^\infty |\psi_{x_0,h}(t)| H(t)\rd t<\infty$. Next, 
we have 
\begin{align}
 &|\hat{F}(x_0)- F(x_0)| =  \bigg|\int_0^{T/2} \psi_{x_0,h}(t) [\hat{H}(t)-H(t)]\rd t +
 \int_0^{T/2} \psi_{x_0,h}(t) H(t)\rd t - F(x_0)\bigg|
 \nonumber
 \\
 &\;\;\leq  \bigg| \int_0^{T/2} \psi_{x_0,h}(t) [\hat{H}(t)-H(t)]\rd t\bigg| +
 \bigg|\int_{T/2}^\infty \psi_{x_0,h}(t) H(t)\rd t\bigg|+
 \bigg|\int_0^{\infty} \psi_{x_0,h}(t) H(t)\rd t - F(x_0)\bigg| 
\nonumber
 \\
 & 
 \;\;= \bigg| \int_0^{T/2} \psi_{x_0,h}(t) [\hat{H}(t)-H(t)]\rd t\bigg| +
  \bigg|\int_{T/2}^\infty \psi_{x_0,h}(t) H(t)\rd t\bigg|
 +
 \bigg|\int_0^{\infty} \varphi_{x_0,h}(x) \rd F(x)  - F(x_0)\bigg|,
\label{eq:error}
 \end{align}
where in the last line we have used (\ref{eq:lin-strategy}). Our goal is to derive bounds on the expectation
of the squared terms on the right hand side of the above formula. 
\par 
Let $\epsilon>0$ be a small number to be specified. Throughout the proof  in the definition of 
$\psi_{x_0, h}$ 
in (\ref{eq:psi}) we put $s=1-\epsilon$.
\par 
1$^0$.
By the Cauchy--Schwarz inequality 
\begin{align*}
 \bE \bigg| \int_0^{T/2} \psi_{x_0,h}(t) [\hat{H}(t)-H(t)]\rd t\bigg|^2 
 \leq
 \Big(\int_0^{T/2} |\psi_{x_0,h}(t)|^2 t^{1-2\epsilon} \rd t\Big) \Big(\int_0^{T/2}
 \frac{1}{t^{1-2\epsilon}} \bE |\hat{H}(t)-H(t)|^2 \rd t\Big).
\end{align*}
In view of  Theorem~\ref{lem:R-Rhat},
\begin{equation}\label{eq:Hhat-H}
 \int_0^{T/2}
 \frac{1}{t^{1-2\epsilon}} \bE |\hat{H}(t)-H(t)|^2 \rd t \leq \frac{c_1}{2\epsilon T^{1-2\epsilon}}\Big(1+\frac{1}{\rho}\Big)
 \int_0^T H(t)\rd t,
\end{equation}
which provides an upper bound on the second integral in previous display formula.
Moreover,  by (\ref{eq:Parseval-Mellin}),
\begin{align*}
 &\int_0^{T/2} |\psi_{x_0,h}(t)|^2 t^{-2\epsilon+1}\rd t  
  = \int_0^{T/2} 
 |\psi_{x_0,h}(t)|^2 t^{2(1-\epsilon) - 1}\rd t
 \\
 &\;\;\leq 
 \frac{1}{2\pi}\int_{-\infty}^\infty |\widetilde{\psi}_{x_0,h}(1-\epsilon + i\omega)|^2\rd \omega
 =
 \frac{1}{2\pi}\int_{-\infty}^\infty 
 \bigg|\frac{\widetilde{\varphi}_{x_0,h}(\epsilon-i\omega)}{\widetilde{w}(\epsilon-i\omega)}\bigg|^2\rd \omega
  \\
   &\;\; \leq \frac{1}{\pi} \int_{-\infty}^\infty
 \bigg| \frac{x_0^{\epsilon} \widehat{K}((-\epsilon + i\omega)h)}{[\epsilon^2+\omega^2]
 \widetilde{w}(\epsilon-i\omega)}
 \bigg|^2 \rd \omega  
 + \frac{1}{\pi}
 \int_{-\infty}^\infty
 \bigg| \frac{h^{\epsilon} \widetilde{K}(1+\epsilon - i\omega)}
 {[\epsilon^2+\omega^2]\widetilde{w}(\epsilon-i\omega)}
 \bigg|^2 \rd \omega 
 =: I_1+I_2.
\end{align*}
We proceed with bounding the integrals $I_1$, $I_2$ on the right hand side. 
\par 
By Lemma~\ref{lem:Mellin-w} and (\ref{eq:hatKbound}) and using the same reasoning as 
in the proof of Lemma~\ref{lem:psi}, we obtain
\begin{align*}
 I_1  &\leq c_1 \Big(\frac{x_0}{2r}\Big)^{2\epsilon} \bigg[
 \int_{|\omega|\leq 2} |\widehat{K}((-\epsilon+ i\omega)h)|^2 \rd \omega
 + \int_{|\omega|\geq 2} |\widehat{K}((-\epsilon+ i\omega)h)|^2 |\omega|^{d+1}\rd \omega
\bigg]
\\*[2mm]
& \leq c_2 \Big(\frac{x_0}{2r}\Big)^{2\epsilon} \bigg[
e^{2\epsilon h} + h^{-d-2}\int_{|\xi|\geq 2h} |\widehat{K}(-\epsilon h + i\xi)|^2 
|\xi|^{d+1}\rd \xi\bigg]
\\*[2mm]
& \leq c_3 e^{2 \epsilon h}\Big(\frac{x_0}{2r}\Big)^{2 \epsilon}
\bigg[1+  h^{-d-2}\bigg(1+ \int_{|\xi|\geq 1} \frac{|\xi|^{d+1}}{[\epsilon^2h^2+ \xi^2]^\ell}
\rd \xi \bigg)\bigg]\leq c_4 e^{2 \epsilon h}\Big(\frac{x_0}{2r}\Big)^{2\epsilon} h^{-d-2},
 \end{align*}
where we have used inequality (\ref{eq:hatKbound}) with  $2\ell> d+2$.
Similarly, 
\begin{align*}
 I_2 \leq c_5 \Big(\frac{h}{2r}\Big)^{2\epsilon}
  \bigg\{\int_{|\omega|\leq 2}  |\widetilde{K}(1+\epsilon-i\omega)|^2 \rd \omega +
 \int_{|\omega|\geq 2} |\widetilde{K}(1+\epsilon-i\omega)|^2 |\omega|^{d+1}\rd\omega \bigg\} \leq 
 c_6 \Big(\frac{h}{2r}\Big)^{2\epsilon}.
\end{align*}
Thus 
\[
 \int_0^{T/2} |\psi_{x_0,h}(t)|^2 t^{-2\epsilon+1}\rd t \leq  c_4 e^{2 \epsilon h}\Big(\frac{x_0}{2r}\Big)^{2\epsilon} h^{-d-2} + c_6 \Big(\frac{h}{2r}\Big)^{2\epsilon} \leq c_7 e^{2 \epsilon h}\Big(\frac{x_0}{2r}\Big)^{2\epsilon} h^{-d-2}.
\]
Combining this inequality with (\ref{eq:Hhat-H}) we obtain 
\begin{equation}\label{eq:var-term}
\bE \bigg| \int_0^{T/2} \psi_{x_0,h}(t) [\hat{H}(t)-H(t)]\rd t\bigg|^2 
 \leq
 c_8 \Big(\frac{x_0}{2r}\Big)^{2\epsilon} \frac{e^{2 \epsilon h}(1+\frac{1}{\rho})}{\epsilon h^{d+2}
 T^{1-2\epsilon}} \int_0^T H(t)\rd t.
\end{equation}
\par 
2$^0$.  Now we bound the second term on the right hand side of (\ref{eq:error}).  
It follows from (\ref{eq:psi-bound}) applied with $s=1-\epsilon$ that 
\[
 |\psi_{x_0,h}(t)| \leq c_1 (2r)^\epsilon t^{-1+\epsilon}\Big[h^\epsilon + x_0^\epsilon e^{\epsilon h} 
 h^{-(d+3)/2}\Big].
\]
Therefore using a bound on $H(t)$ in (\ref{eq:H(t)<=}) for sufficiently large $T$ we obtain
\begin{align}
\int_{T/2}^\infty |\psi_{x_0,h}(t)| H(t)\rd t \leq 
c_1 (2r)^\epsilon \Big[h^\epsilon + x_0^\epsilon e^{\epsilon h} 
 h^{-(d+3)/2}\Big] \int_{T/2}^\infty H(t) t^{-1+\epsilon} \rd t
\nonumber
 \\
 \leq \frac{c_2M}{1+\alpha} (2r)^{1+\alpha+\epsilon} \Big[h^\epsilon + x_0^\epsilon e^{\epsilon h} 
 h^{-(d+3)/2}\Big]
 \int_{T/2}^\infty t^{-2-\alpha+\epsilon}\rd t
 \nonumber
 \\
 \leq \frac{c_3 Mx_0^\epsilon e^{\epsilon h}} 
 {(1+\alpha)(1+\alpha-\epsilon)} \bigg(\frac{1}{h^{(d+3)/2}T^{1+\alpha-\epsilon}}\bigg),
\label{eq:bias-first-term}
 \end{align}
provided that $1+\alpha-\epsilon>0$.
\par 
3$^0$. Now we work with the thrid term on the right hand side of (\ref{eq:error}). 
By definition of $\varphi_{x_0,h}$ we have 
\begin{align}
& \int_0^{\infty} \varphi_{x_0,h}(x)\rd F(x) = \int_0^h \varphi_{x_0,h}(x)\rd F(x) + F(x_0)-F(h) + 
 \int_{x_0}^{x_0e^h} \varphi_{x_0,h}(x)\rd F(x),
 \label{eq:111}
 \\
& \int_0^h \varphi_{x_0,h}(x)\rd F(x) = \int_0^h \int_0^h {\bf 1}(t\leq x) \frac{1}{h}K\Big(\frac{t}{h}\Big)
 \rd t \rd F(x)= F(h) - \frac{1}{h}\int_0^h K\Big(\frac{t}{h}\Big)F(t)\rd t
\label{eq:222}
 \end{align}
and 
\begin{align}
 &\int_{x_0}^{x_0e^h} \varphi_{x_0}(x)\rd F(x)= \int_{x_0}^{x_0e^h} 
 \Big[1- \int_{x_0}^x \frac{1}{th} K\Big(\frac{\ln (t/x_0)}{h}\Big)\Big] \rd F(x)
 \nonumber
 \\
 &= F(x_0e^h)- F(x_0) - \int_{x_0}^{x_0e^h}\int_{x_0}^{x_0e^h} {\bf 1}(t\leq x) \frac{t}{th}
 K\Big(\frac{\ln (t/x_0)}{h}\Big)\rd t\rd F(x)
\nonumber
 \\
 & = -F(x_0) + \int_{x_0}^{x_0e^h} \frac{1}{th}K\Big(\frac{\ln (t/x_0)}{h}\Big)F(t)\rd t = 
 \int_0^1 K(y) F(x_0e^{yh}) \rd y - F(x_0).
\label{eq:333}
 \end{align}
Combining (\ref{eq:111}), (\ref{eq:222}) and (\ref{eq:333}) we obtain 
\begin{align*}
 \bigg|
 \int_0^{\infty} \varphi_{x_0}(x)\rd F(x) - F(x_0)\bigg| \leq 
 \bigg|\frac{1}{h}\int_0^h K\Big(\frac{t}{h}\Big)F(t)\rd t\bigg| + 
 \bigg|\int_0^1 K(y) F(x_0e^{yh}) \rd y - F(x_0)\bigg|.
\end{align*}
\par 
Because $F\in \sH_{\beta}(A)$ and $F(0)=0$, expanding in Taylor's series we have for some 
$\xi \in (0, h)$
\begin{align*}
 \bigg|\frac{1}{h}\int_0^h K\Big(\frac{t}{h}\Big)F(t)\rd t\bigg| \leq 
 \frac{h^\ell}{\ell !}\int_0^1 |K(y)|\,|F^{(\ell)}(\xi) - F^{(\ell)}(0+)| y^\ell \rd y\leq 
 c_1 Ah^{\beta}.
\end{align*}
To bound the second term we define function 
$R_{x_0}(t)= F(x_0e^t)$; with this notation 
\begin{align*}
 \bigg|\int_0^1 K(y) F(x_0e^{yh}) \rd y - F(x_0)\bigg| = 
 \bigg|\int_0^1 K(y) [R_{x_0}(yh)- R_{x_0}(0)] \rd y\bigg|.
\end{align*}
Expanding function $R_{x_0}$ in Taylor's series around $0$ we have 
\[
 R_{x_0}(yh)= R_{x_0}(0) + \sum_{j=1}^{\ell-1} \frac{R_{x_0}^{(j)}(0)}{j!} (yh)^j + \frac{1}{\ell!} 
 R^{(\ell)}_{x_0}(\xi) (yh)^\ell,\;\;\;0\leq \xi \leq h,
\]
and therefore 
\begin{equation}\label{eq:bias-11}
 \bigg|\int_0^1 K(y) [R_{x_0}(yh)- R_{x_0}(0)] \rd y\bigg| \leq 
 \frac{h^\ell}{\ell!}\int_0^1 |K(y)|  |R^{(\ell)}_{x_0}(\xi) - R_{x_0}^{(\ell)}(0)| |y|^\ell \rd y.
\end{equation}
By the Fa\'a~di~Bruno formula 
\[
 R_{x_0}^{(\ell)}(t)= \frac{d^\ell}{dt^\ell} F(x_0e^t)= \sum_{k=1}^\ell F^{(k)}(x_0e^t) B_{\ell, k}(x_0e^t,\ldots, x_0e^t),
\]
where $B_{\ell, k}$ is the Bell polynomial of degree $\ell$ in $\ell-k+1$ variables given by the formula 
\[
 B_{\ell, k}(x_1, \ldots, x_{\ell-k+1})=
 \sum \frac{\ell!}{j_1!\cdots j_{\ell-k+1}!}\Big(\frac{x_1}{1!}\Big)^{j_1} \cdots 
 \Big(\frac{x_{\ell-k+1}}{(\ell-k+1)!}\Big)^{j_{\ell-k+1}};
\]
the sum is taken over all over all subsets $j_1, ..., j_{\ell-k+1}$ 
of non--negative integers such that 
$j_1+\cdots+j_{\ell-k+1}=k$ and 
$j_1+2j_2+\cdots+(\ell-k+1)j_{\ell-k+1}=\ell$. 
In our  specific case the Fa\'a~di~Bruno formula takes the form  
\[
 R_{x_0}^{(\ell)}(t)= \sum_{k=1}^\ell c_{k,l} F^{(k)}(x_0e^t) (x_0e^t)^k,
\]
where $c_{k,l}$ are coefficients depending on $k$ and $l$ only.
Therefore 
\begin{align*}
 &|R_{x_0}^{(\ell)}(\xi)- R_{x_0}^{(\ell)}(0)| = 
 \Big|\sum_{k=1}^\ell c_{k,l}x_0^k \big[F^{(k)}(x_0 e^\xi) e^{k\xi} - F^{(k)}(x_0)\big]\Big|
 \\
 &\leq \sum_{k=1}^\ell c_{k,l}x_0^k \Big[e^{k\xi} |F^{(k)}(x_0e^\xi)- F^{(k)}(x_0)|
 + |F^{(k)}(x_0)| |e^{k\xi}-1|\Big]
 \\
 &\leq \sum_{k=1}^{\ell-1} c_{k,l}x_0^k \Big[ e^{kh} Ax_0 |e^h-1| +A|e^{kh}-1|\Big]+ c_{\ell,\ell}x_0^\ell
 \Big[Ax_0^{\beta-\ell}|e^h-1|^{\beta-\ell} + A|e^{\ell h}-1|\Big]
 \\
 &\leq c_2 A  \Big[ x_0^\beta h^{\beta-\ell} + h\sum_{k=1}^{\ell-1} x_0^{k+1}\Big],
\end{align*}
where constant $c_2$ depends on $\beta$ only.
To obtain the last formula we have used the 
elementary inequality $e^x-1\leq x e^x$, $x\geq 0$,  and the fact that $h<1/2$. 
Combining this inequality with (\ref{eq:bias-11}) we obtain 
\begin{align*}
  \bigg|\int_0^1 K(y) [R_{x_0}(yh)- R_{x_0}(0)] \rd y\bigg| \leq c_3 A 
  \Big[ x_0^\beta h^{\beta} + h^{\ell+1}\sum_{k=1}^{\ell-1} x_0^{k+1}\Big], 
\end{align*}
and finally
\begin{align}\label{eq:bias-second-term}
 \bigg|
 \int_0^{\infty} \varphi_{x_0,h}(x)\rd F(x) - F(x_0)\bigg| \leq  c_4 Ah^\beta(x_0^\beta+1) +   
 c_3Ah^{\ell+1}\sum_{k=1}^{\ell-1} x_0^{k+1}.
\end{align}
\par 
4$^0$. Now we are in a position to complete the proof of the theorem.
We combine bounds (\ref{eq:var-term}), (\ref{eq:bias-first-term}) and (\ref{eq:bias-second-term})
and set $\epsilon=1/\ln T$. 
Then for $T$ large enough we obtain 
\begin{align*}
 \bE |\hat{F}(x_0)-F(x_0)|^2 \leq  \frac{c_1(1+\frac{1}{\rho})\ln T}{h^{d+2}
 T} \int_0^T H(t)\rd t +  \frac{c_2M}{h^{d+3}T^{2}}
+ c_3A^2 h^{2\beta} (x_0^\beta+1)^2
 \end{align*}
In view of (\ref{eq:int-H}) we have  
 \begin{align*}
  \int_0^T H(t)\rd t \leq c_4\tilde{\eta}_T,\;\;\;\tilde{\eta}_T=[r + Mr^{1\vee (1+\alpha)} \eta_T],\;\;\;
 \end{align*}
 and we recall that $\eta_T$ is $1$ for $\alpha>0$; $\ln T$ for $\alpha=0$, and $T^{-\alpha}$ for 
 $-1<\alpha<0$ [see (\ref{eq:int-H})]. 
Then the  choice $h=h_*$ as in (\ref{eq:h-*}) leads to the announced result. 
\epr

\section{Proofs for Section~\ref{sec:Brownian}}
\label{sec:proofs-Brownian}
\subsection{Proof of Lemma~\ref{lem:cov-properties}}
\par\medskip
(a).~{\em  Behavior at zero.} If $d=1$ then by straightforward algebra
\begin{align*}
 H(t)&= \frac{1}{B(1, \frac{1}{2}) \Gamma(\frac{1}{2})} 
 \int_0^1 \int_0^\infty x^{-1/2} e^{-x} {\bf 1}\Big(x\leq \frac{2r^2y}{\sigma^2 t}\Big) y^{-1/2}\rd x
 \rd y
 \\
 &= \frac{2\Gamma(\frac{3}{2})}{\Gamma^2(\frac{1}{2})} 
 \int_0^\infty
x^{-1/2} e^{-x}\Big(1-\frac{\sigma \sqrt{tx}}{\sqrt{2}r}\Big) {\bf 1}
 \Big(x\leq \frac{2r^2}{\sigma^2 t}\Big) \rd x
 \\
 &= \frac{2\Gamma(\frac{3}{2})}{\Gamma^2(\frac{1}{2})} 
 \bigg[\int_0^\infty x^{-1/2} e^{-x}{\bf 1}
 \Big(x\leq \frac{2r^2}{\sigma^2 t}\Big) \rd x + \frac{\sigma \sqrt{t}}{\sqrt{2} r}
 \exp\Big\{-\frac{2r^2}{\sigma^2 t}\Big\} - \frac{\sigma \sqrt{t}}{\sqrt{2} r}
 \bigg]
\end{align*}
so that 
\[
 1- H(t)  \sim  \frac{\sigma\sqrt{t}}{\sqrt{2}\Gamma(\frac{1}{2}) r}=\frac{\sigma\sqrt{t}}{\sqrt{2\pi} r}\;\;\;\hbox{as}\;\;t\to 0.
\]
\par 
If $d\geq 2$ then the asymptotic approximation of the upper incomplete Gamma function
[see, e.g., \cite[6.5.32]{Abr-Ste}] yields:
 \begin{equation}\label{eq:1-H(t)}
  1-H(t) \sim  \frac{1}{B(\tfrac{d+1}{2}, \tfrac{1}{2}) \Gamma(\tfrac{d}{2})} 
  \Big(\frac{2r^2}{\sigma^2 t}\Big)^{\frac{d}{2}-1} 
  \int_0^1 \exp\Big(-\frac{2r^2 y}{\sigma^2 t}\Big) 
  (1-y)^{\frac{1}{2}(d-1)} y^{\frac{1}{2}(d-3)}\rd y,\;\;\;
  t\to 0.
 \end{equation}
The  integral on the right hand side is expressed in terms of  Kummer's function 
[cf.~\cite[Chapter~13]{Abr-Ste}] that is defined as follows:
for $a$ and $b$ satisfying ${\rm Re}(b)>{\rm Re}(a)>0$
\[
 M(a,b; z):= 1+ \sum_{k=1}^\infty \bigg(\prod_{j=0}^{k-1} \frac{a+j}{b+j}\bigg) \frac{z^k}{k!}
 = \frac{\Gamma(b)}{\Gamma(b-a) \Gamma(a)} \int_0^1 e^{zu} u^{a-1} (1-u)^{b-a-1} \rd u.
\]
With this notation letting $a=(d-1)/2$ and $b=d$ we have
\begin{align}\label{eq:int=}
 \int_0^1 \exp\Big(-\frac{2r^2 y}{\sigma^2 t}\Big) (1-y)^{\frac{1}{2}(d-1)} y^{\frac{1}{2}(d-3)}\rd y = \frac{\Gamma(\frac{d+1}{2}) \Gamma(\frac{d-1}{2})}{\Gamma(d)}
 M\Big(\frac{d-1}{2}, d; -\frac{2r^2}{\sigma^2t}\Big).
\end{align}
By \cite[13.1.5]{Abr-Ste},  
$M(a, b; z)= \frac{\Gamma(b)}{\Gamma(b-a)} (-z)^{-a}\big[1+O(|z|^{-1})\big]$, ${\rm Re}(z)<0$ as $|z|\to\infty$.
Therefore combining (\ref{eq:int=}) and (\ref{eq:1-H(t)}) we obtain  
\begin{align*}
 1- H(t)\;&\sim\; \frac{\Gamma(\frac{d}{2}+1)}{\Gamma(\frac{1}{2})\Gamma(\frac{d+1}{2})\Gamma(\frac{d}{2})} \Big(\frac{2r^2}{\sigma^2t}\Big)^{\frac{d}{2}-1}
 \frac{\Gamma(\frac{d+1}{2}) \Gamma(\frac{d-1}{2})}{\Gamma(d)}\cdot 
 \frac{\Gamma(d)}{\Gamma(\frac{d+1}{2})}
 \Big(\frac{2r^2}{\sigma^2t}\Big)^{-\frac{d-1}{2}}
 \\
 & =  \frac{d}{(d-1) \Gamma(\frac{1}{2})}
 \Big(\frac{\sigma^2 t}{2r^2}\Big)^{1/2}= \frac{d}{d-1}
 \Big(\frac{\sigma \sqrt{t}}{\sqrt{2\pi}r}\Big)\;\;\;\;\hbox{as}\;\;t\to 0.
\end{align*}
These calculations show that   
  $H^\prime (0+)=\infty$. 
 \par 
 {\em Behavior at infinity.} As $t\to\infty$ we have
\[
 H(t) \sim \frac{1}{B(\tfrac{d+1}{2}, \tfrac{1}{2}) \Gamma(\tfrac{d}{2})}
 \int_0^1  \frac{2}{d}\Big(\frac{2r^2y}{\sigma^2 t}\Big)^{d/2}(1-y)^{(d-1)/2} y^{-1/2}\rd y = 
 \frac{\Gamma(\frac{d+1}{2})}{\Gamma(d+1)  \Gamma(\frac{1}{2})} \bigg(\frac{\sqrt{2}r}{\sigma\sqrt{t}}\bigg)^{d},
\]
as claimed.
\par\medskip 
(b). It follows from (\ref{eq:covariance}) that  
\begin{align*}
 H(t)= \frac{1}{B(\frac{d+1}{2},\frac{1}{2})\Gamma(\frac{d}{2})} 
 \int_0^1\int_0^\infty x^{\frac{d}{2}-1} e^{-x} {\bf 1}\Big(x\leq \frac{2r^2y}{\sigma^2 t}\Big)
 (1-y)^{\frac{1}{2}(d-1)} y^{-1/2}\rd x\rd y
 \\
 \leq \frac{1}{B(\frac{d+1}{2},\frac{1}{2})\Gamma(\frac{d}{2})} \int_0^1
 \frac{2}{d} \Big(\frac{2r^2y}{\sigma^2 t}\Big)^{d/2} (1-y)^{\frac{1}{2}(d-1)} y^{-1/2}\rd y 
 \\
 = \frac{\Gamma(\frac{d+1}{2})}{\Gamma(d+1) \Gamma(\frac{1}{2})} \Big(\frac{2r^2}{\sigma^2 t}\Big)^{d/2},\;\;\;\forall t>0.
\end{align*} 
 We also have  
\begin{align*}
\int_0^T H(t) \rd t =  \frac{1}{B(\frac{d+1}{2}, \frac{1}{2})\Gamma(\frac{d}{2})}
\int_0^1 \int_0^\infty x^{\frac{d}{2}-1}e^{-x} \min\Big\{ \frac{2r^2y}{\sigma^2 x}, T\Big\}
(1-y)^{\tfrac{1}{2}(d-1)} y^{-1/2} \rd x\rd y.
\end{align*}
If $d>2$ then 
\begin{align*}
 \int_0^T H(t) \rd t  \leq 
 \frac{1}{B(\frac{d+1}{2}, \frac{1}{2})\Gamma(\frac{d}{2})}
\int_0^1 \int_0^\infty x^{\frac{d}{2}-1}e^{-x}  \frac{2r^2y}{\sigma^2 x}
(1-y)^{\tfrac{1}{2}(d-1)} y^{-1/2} \rd x\rd y
 \\
  = 
 \frac{2r^2}{\sigma^2} \frac{\Gamma(\frac{d}{2}-1) B(\frac{d+1}{2}, \frac{3}{2})}{B(\frac{d+1}{2},\frac{1}{2}) \Gamma(\frac{d}{2})}= \frac{4r^2}{\sigma^2(d^2-4)}.
\end{align*}
If $d=2$ then  
\begin{align*}
 H(t) & = \frac{1}{B(\frac{3}{2}, \frac{1}{2})} \int_0^1 \int_0^\infty e^{-x} 
 {\bf 1} \Big(x\leq \frac{2r^2y}{\sigma^2 t}\Big) (1-y)^{1/2} y^{-1/2} \rd x\rd y 
 \\
 &= 
 \frac{1}{B(\frac{3}{2},\frac{1}{2})} \int_0^1 
 \Big[1-\exp\Big\{-\frac{2r^2y}{\sigma^2 t}\Big\}\Big]
 (1-y)^{1/2} y^{-1/2} \rd y 
 \leq 
 \frac{r^2}{2\sigma^2 t},\;\;\;\forall t>0,
\end{align*}
where we have used the elementary inequality $1-e^{-x}\leq x$. Therefore 
for $d=2$ 
\[
 \int_0^T H(t) \rd t \leq 1 + \frac{r^2}{2\sigma^2}\int_1^T\frac{\rd t}{t} = 1+ \frac{r^2}{2\sigma^2}
 \ln T.
\]
If $d=1$ then 
\begin{align*}
 H(t) &= \frac{1}{B(1, \frac{1}{2})\Gamma(\frac{1}{2})}
 \int_0^1 \int_0^\infty x^{-1/2} e^{-x}  {\bf 1}\Big(x\leq \frac{2r^2y}{\sigma^2t}\Big) 
 y^{-1/2}\rd x\rd y
\\
&\leq \frac{2}{B(1, \frac{1}{2})\Gamma(\frac{1}{2})} \int_0^1  
\Big(\frac{2r^2y}{\sigma^2t}\Big)^{1/2} y^{-1/2}\rd y =\sqrt{\frac{2}{\pi}} \Big(\frac{r^2}{\sigma^2 t}
\Big)^{1/2},\;\;\forall t>0.
 \end{align*}
 This completes the proof.
\epr

\subsection{Proof of Theorem~\ref{th:sigma-1}}
In the subsequent proof $c_1, c_2, \ldots$ stand for positive constant that may depend on $d$
only. The proof is divided in two steps.
\par\medskip
(a). 
First we establish an upper bound on accuracy of estimating the functional $\Psi_\alpha$. 
We have 
\[
 \hat{\Psi}_{\alpha,b}-\Psi_\alpha=  \int_0^b \frac{\hat{H}(t)-H(t)}{t^{1-\alpha}} \rd t + \int_b^\infty 
 \frac{H(t)}{t^{1-\alpha}}\rd t;
\]
hence 
\begin{align*}
 \bE |\hat{\Psi}_{\alpha, b} - \Psi_\alpha|^2 \leq 2 
 \int_0^b\int_0^b \frac{\bE [(\hat{H}(t)-H(t))(\hat{H}(s)-H(s))]}{t^{1-\alpha}s^{1-\alpha}}
 \rd t\rd s 
 \\
 + 2\bigg(\int_b^\infty \frac{H(t)}{t^{1-\alpha}}\rd t\bigg)^2
 =: 2I_1 +2 I_2.
\end{align*}
\par 
The bound on $I_2$ is readily obtained from (\ref{eq:H(t)leq}): 
\begin{align*}
 \int_b^\infty \frac{H(t)}{t^{1- \alpha}} \rd t \leq c_1 \Big(\frac{r}{\sigma}\Big)^d \int_b^\infty
 \frac{\rd t}{t^{1-\alpha+d/2}} = \Big(\frac{2c_1}{d-2\alpha}\Big) 
 \Big(\frac{r}{\sigma}\Big)^d b^{\alpha-d/2},
\end{align*}
so that 
\[
 I_2 \leq \Big(\frac{2c_1}{d-2\alpha}\Big)^2  \Big(\frac{r}{\sigma}\Big)^{2d} b^{-d+2\alpha}~.
\]
\par 
Now we bound $I_1$: by 
the Cauchy--Schwarz inequality, (\ref{eq:cov-err}) and (\ref{eq:H-int})
\begin{align*}
 I_1\leq \bigg[\int_0^b \frac{1}{t^{1-\alpha}}
 \Big(\bE |\hat{H}(t)-H(t)|^2\Big)^{1/2}\rd t\bigg]^2 \leq c_1 \Big(1+\frac{1}{\rho}\Big) \int_0^T H(t)\rd t 
  \bigg[\int_0^b \frac{\rd t}{t^{1-\alpha} (T-t)^{1/2}} \bigg]^2
  \\
  \leq 
  c_2\Big(1+\frac{1}{\rho}\Big)\frac{b^{2\alpha}}{\alpha^2T} \int_0^T H(t)\rd t. 
\end{align*}
Using (\ref{eq:H-int}) and combining the bounds for $I_1$ and $I_2$ we obtain the following results.
\par 
If $d=1$ then 
 \begin{equation*}
   \bE |\hat{\Psi}_{\alpha, b} - \Psi_\alpha|^2 \leq  c_3 \Big(1+\frac{1}{\rho}\Big) \Big(\frac{r}{\sigma}\Big)
 \frac{b^{2\alpha}}{\alpha^2 \sqrt{T}} + \Big(\frac{2c_1}{1-2\alpha}\Big)^2 
 \Big(\frac{r}{\sigma}\Big)^2 b^{-1+2\alpha}.
\end{equation*}
Letting $b_*= (1-2\alpha)^{-2}\alpha^2 \sqrt{T}$ we obtain 
\begin{equation}\label{eq:d=1}
 \bE |\hat{\Psi}_{\alpha, b_*} - \Psi_\alpha|^2 \leq  c_4 \Big[ \Big(1+\frac{1}{\rho}\Big)
 \Big(\frac{r}{\sigma}\Big) + \Big(\frac{r}{\sigma}\Big)^2\Big]
 \frac{\alpha^{-2+4\alpha}}{(1-2\alpha)^{4\alpha}T^{1/2-\alpha}}~.
\end{equation}
\par 
If $d=2$ then 
\begin{equation*}
   \bE |\hat{\Psi}_{\alpha, b} - \Psi_\alpha|^2 \leq  c_5 \Big(1+\frac{1}{\rho}\Big) 
 \frac{b^{2\alpha}}{\alpha^2 T} 
 \Big[1+ \Big(\frac{r}{\sigma}\Big)^2\ln T\Big]
 + \Big(\frac{2c_1}{2-2\alpha}\Big)^2 
 \Big(\frac{r}{\sigma}\Big)^4 b^{-2+2\alpha}.
\end{equation*}
Therefore letting $b_*=\alpha \sqrt{T/\ln T}$ we obtain
for sufficiently large $T$
\begin{equation}\label{eq:d=2}
  \bE |\hat{\Psi}_{\alpha, b_*} - \Psi_\alpha|^2  \leq 
  c_6 \Big[\Big(1+\frac{1}{\rho}\Big) \Big(\frac{r}{\sigma}\Big)^2 +\Big(\frac{r}{\sigma}\Big)^4
  \Big]
  \alpha^{-2+2\alpha} 
  \Big(\frac{\ln T}{T}\Big)^{1-\alpha}.
\end{equation}
\par 
Finally, if $d>2$ then 
\begin{equation*}
   \bE |\hat{\Psi}_{\alpha, b} - \Psi_\alpha|^2 \leq  c_7 \Big(1+\frac{1}{\rho}\Big)
   \Big(\frac{r}{\sigma}\Big)^2
 \frac{b^{2\alpha}}{\alpha^2 T} 
 + \Big(\frac{2c_1}{d-2\alpha}\Big)^2 
 \Big(\frac{r}{\sigma}\Big)^{2d} b^{-d+2\alpha},
\end{equation*}
and for $b_*=(\alpha^2 T)^{1/d}$ we get 
\begin{equation}\label{eq:d>2}
  \bE |\hat{\Psi}_{\alpha, b_*} - \Psi_\alpha|^2 \leq c_8  \Big[\Big(1+\frac{1}{\rho}\Big)
   \Big(\frac{r}{\sigma}\Big)^2 +\Big(\frac{r}{\sigma}\Big)^{2d} \Big]
   \Big(\frac{1}{\alpha^2 T}\Big)^{1 - \frac{2\alpha}{d}}.
\end{equation}
\par\medskip 
(b). Now we relate error in estimating $\sigma^2$ by $\hat{\sigma}_{\alpha, b}$ 
to the mean squared error of 
$\hat{\Psi}_{\alpha, b}$. 
 By definition 
\begin{align*}
 & |\hat{\sigma}_{\alpha, b}^2 - \sigma^2| = 
 \Big(\frac{J_{\alpha}}{\Psi_\alpha\hat{\Psi}_{\alpha,b}}\Big)^{1/\alpha}
 |\Psi_\alpha^{1/\alpha}-\hat{\Psi}_{\alpha,b}^{1/\alpha}| \leq  
  \Big(\frac{J_{\alpha}}{\Psi_\alpha \hat{\Psi}_{\alpha,b}}\Big)^{1/\alpha}\frac{1}{\alpha} 
  \big(\hat{\Psi}_{\alpha,b} \vee \Psi_\alpha\big)^{\frac{1}{\alpha}-1} 
  |\Psi_\alpha-\hat{\Psi}_{\alpha,b}|
 \\
 &=\Big(\frac{J_{\alpha}}{\Psi_\alpha}\Big)^{1/\alpha} {\bf 1}(\hat{\Psi}_{\alpha,b}\geq 
 \Psi_\alpha) 
 \frac{1}{\alpha\hat{\Psi}_{\alpha,b}} |\Psi_\alpha-\hat{\Psi}_{\alpha,b}| + 
 \Big(\frac{J_{\alpha}}{\hat{\Psi}_{\alpha, b}}\Big)^{1/\alpha} {\bf 1}
 (\hat{\Psi}_{\alpha, b} < \Psi_\alpha) 
 \frac{1}{\alpha\Psi_\alpha} |\Psi_\alpha-\hat{\Psi}_{\alpha,b}|
 \\
 &\leq \frac{1}{\alpha \Psi_\alpha} |\Psi_\alpha-\hat{\Psi}_{\alpha, b}|
 \Big[\sigma^2 {\bf 1}(\hat{\Psi}_{\alpha, b} \geq \Psi_\alpha)+
 \hat{\sigma}_{\alpha, b}^2 {\bf 1}(\hat{\Psi}_{\alpha,b}<\Psi_\alpha)\Big]
 \leq 
 \frac{1}{\alpha\Psi_\alpha} |\Psi_\alpha-\hat{\Psi}_{\alpha, b}|\big(\sigma^2 +
 \hat{\sigma}_{\alpha, b}^2\big),
  \end{align*}
 where in the first line we have used the elementary inequality 
 $|a^{1/\alpha}- b^{1/\alpha}|\leq \frac{1}{\alpha} (a\vee b)^{\frac{1}{\alpha}-1}|a-b|$ which holds 
 for all $a, b>0$ and $0<\alpha\leq 1$. 
Therefore  
\begin{equation}\label{eq:Delta-bound}
\bE \big[\Delta (\hat{\sigma}^2_{\alpha, b}, \sigma^2)\big]^2 = \bE \bigg|\frac{\hat{\sigma}_{\alpha, b}^2- \sigma^2}{\hat{\sigma}_{\alpha, b}^2 + \sigma^2}\bigg|^2
\leq \frac{\sigma^{4\alpha}}{\alpha^2J_{\alpha}^2} \bE |\Psi_\alpha-\hat{\Psi}_{\alpha,b}|^2
\leq c_1 \Big(\frac{\sigma}{r}\Big)^{4\alpha} \bE |\Psi_\alpha-\hat{\Psi}_{\alpha,b}|^2,
 \end{equation}
where we have taken into account that
 \[
  \alpha^2 J_\alpha^2 = 
  \bigg[\frac{(2r^2)^\alpha\Gamma(\frac{d}{2}-\alpha)B(\frac{1}{2}+\alpha, \frac{d+1}{2})}
 {\Gamma(\frac{d}{2})B(\frac{d+1}{2},\frac{1}{2})}\bigg]^2\geq c_2 r^{4\alpha},\;\;\;\forall \alpha\in (0, 1/2).
 \]
\par
To complete the proof we combine 
(\ref{eq:Delta-bound}) with (\ref{eq:d=1}), (\ref{eq:d=2}) and (\ref{eq:d>2}). 
and set $\alpha_*=1/\ln T$. 
\epr  
 \subsection{Proof of Lemma~\ref{lem:d=1}}
We have 
\begin{align*}
 H(t)&= \frac{1}{B(1, \frac{1}{2}) \Gamma(\frac{1}{2})} 
 \int_0^1 \int_0^\infty x^{-1/2} e^{-x} {\bf 1}\Big(x\leq \frac{2r^2y}{\sigma^2 t}\Big) y^{-1/2}\rd x
 \rd y
 \\
 &= \frac{2\Gamma(\frac{3}{2})}{\Gamma^2(\frac{1}{2})} 
 \int_0^\infty
x^{-1/2} e^{-x}\Big(1-\frac{\sigma \sqrt{tx}}{\sqrt{2}r}\Big) {\bf 1}
 \Big\{x\leq \frac{2r^2}{\sigma^2 t}\Big\} \rd x
 \\
 &= \frac{2\Gamma(\frac{3}{2})}{\Gamma^2(\frac{1}{2})} 
 \bigg[\int_0^\infty x^{-1/2} e^{-x}{\bf 1}
 \Big\{x\leq \frac{2r^2}{\sigma^2 t}\Big\} \rd x + \frac{\sigma \sqrt{t}}{\sqrt{2} r}
 \exp\Big\{-\frac{2r^2}{\sigma^2 t}\Big\} - \frac{\sigma \sqrt{t}}{\sqrt{2} r}
 \bigg]
 \\
 &=  1- \frac{1}{\sqrt{\pi}}  \int_{2r^2/(\sigma^2 t)}^\infty x^{-1/2}e^{-x}\rd x + \frac{\sigma \sqrt{t}}{\sqrt{2\pi} r}
 \exp\Big\{-\frac{2r^2}{\sigma^2 t}\Big\} - \frac{\sigma \sqrt{t}}{\sqrt{2\pi} r}.
\end{align*}
Then for any $t\leq \sigma^2/(2r^2)$ one has 
\[
 1-H(t)= \frac{\sigma\sqrt{t}}{\sqrt{2\pi} r}\Big[1- \exp\Big\{-\frac{2r^2}{\sigma^2 t}\Big\} \Big] + 
 \delta(t),\;\;\;
|\delta(t)| \leq  \exp\Big\{-\frac{2r^2}{\sigma^2 t}\Big\}.
\]
\epr 
 \subsection{Proof of Theorem~\ref{th:sigma-2}}
It follows from Lemma~\ref{lem:d=1} that 
\[
 \bigg|\sigma - \frac{\sqrt{2\pi} r}{\sqrt{\tau}} [1-H(\tau)]\bigg| \leq
 \exp\Big\{-\frac{r^2}{\sigma^2\tau}\Big\} \Big[\sigma + \frac{\sqrt{\tau}}{\sqrt{2\pi}r}\Big].
\]
Therefore
\begin{align*}
 |\hat{\sigma}_\tau - \sigma| \leq \frac{\sqrt{2\pi} r}{\sqrt{\tau}} |\hat{H}(\tau) - H(\tau)|+ 
 \exp\Big\{-\frac{r^2}{\sigma^2\tau}\Big\} \Big[\sigma + \frac{\sqrt{\tau}}{\sqrt{2\pi}r}\Big],
\end{align*}
and 
\begin{align*}
 \bE |\hat{\sigma}_\tau - \sigma|^2 \leq 
 c\Big\{\Big(1+\frac{1}{\rho}\Big) \frac{1}{T}\int_0^T H(t)\rd t
 +  \exp\Big\{-\frac{2r^2}{\sigma^2\tau}\Big\} \Big[\sigma^2 + \frac{\tau}{r^2}\Big]\Big\}
 \\
 \leq
 c\Big\{\Big(1+\frac{1}{\rho}\Big) \Big(\frac{r}{\sigma}\Big)^3 \frac{\sigma^2}{\tau \sqrt{T}}
 +  \exp\Big\{-\frac{2r^2}{\sigma^2\tau}\Big\} \Big[\sigma^2 + \frac{\tau}{r^2}\Big]\Big\},
\end{align*}
where we have used (\ref{eq:H-int}).
Setting for some $\bar{\sigma}>0$
\[
\tau=\tau_*=\frac{4r^2}{\bar{\sigma}^2\ln T}
\]
we obtain 
\[
 \bE \bigg|\frac{\hat{\sigma}_{\tau_*} - \sigma}{\sigma}\bigg|^2 \leq 
 c\bigg\{\Big(1+\frac{1}{\rho}\Big)\Big(\frac{r}{\sigma}\Big) \Big(
 \frac{\bar{\sigma}}{\sigma}\Big)^2 \frac{\ln T}{\sqrt{T}} + 
 \Big(\frac{1}{\sqrt{T}}\Big)^{\bar{\sigma}^2/\sigma^2}\Big[1+ 
 \frac{1}{\sigma^2\bar{\sigma}^2\ln T}
 \Big]\bigg\}.
\]
\epr

\bibliographystyle{agsm}
 
\end{document}